\documentstyle[epsf]{article}

\newcommand{\bigzerou}{%
\smash{\lower1.7ex\hbox{\bg 0}}}
\newtheorem{theorem}{Theorem}

\newtheorem{defi}{Definition}
\newtheorem{cor}{Corollary}
\newtheorem{conj}{Conjecture}

\newtheorem{Rem}{Remark}
\newtheorem{lem}{Lemma}

\newcommand{\ba}{\begin{eqnarray}}
\newcommand{\ea}{\end{eqnarray}}
\newcommand{\no}{\nonumber}

\newcommand{\dd}{{\bf d}}
\def\d{{\partial}}
\newcommand{\mapright}[1]{%
\smash{\mathop{%
\hbox to 1.0cm{\rightarrowfill}}\limits^{#1}}}
\newcommand{\mapleft}[1]{%
\smash{\mathop{%
\hbox to 1.3cm{\leftarrowfill}}\limits^{#1}}}

\begin{document}
\title{
\begin{flushright}
  \begin{minipage}[b]{5em}
    \normalsize
    ${}$      \\
  \end{minipage}
\end{flushright}
{\bf Mirror Map as Generating Function of Intersection Numbers: Toric Manifolds 
with Two K\"ahler Forms}}
\author{Masao Jinzenji\\
\\
\it Division of Mathematics, Graduate School of Science \\
\it Hokkaido University \\
\it  Kita-ku, Sapporo, 060-0810, Japan\\
{\it e-mail address: jin@math.sci.hokudai.ac.jp}}
\maketitle
\begin{abstract}
In this paper, we extend our geometrical derivation of the expansion coefficients of mirror maps by localization computation to the case of 
toric manifolds with two K\"ahler forms. In particular, we consider Hirzebruch surfaces $F_{0}$, $F_{3}$ and Calabi-Yau hypersurface 
in weighted projective space ${\bf P}(1,1,2,2,2)$ 
as examples. We expect that our results can be easily generalized to arbitrary toric manifolds.  
\end{abstract}
\section{Introduction}
In the study of mirror symmetry, gauged linear sigma model is expected to play an important role \cite{mp}. It has been considered 
to be slightly different from the topological (non-linear) sigma model, whose correlation function is nothing but the Gromov-Witten 
invariant. Let us restrict our attention to the genus $0$ Gromov-Witten invariants of toric manifolds.
The moduli space used in topological (non-linear) sigma model is the moduli space of stable maps, which 
is a compactification of the moduli space of holomorphic maps from $CP^1$ to toric manifolds by using stable maps.  
On the other hand, the moduli space used in gauged sigma model is another compactification (toric compactification)
of the moduli space of holomorphic maps from $CP^1$ to toric manifold. In this case, we use "rational maps" from $CP^1$ to toric manifolds to compactify the moduli space.
A rational map $f:X\rightarrow Y$ is the map which allows some Zariski-closed subset $U\subset X$ whose image is undefined by $f$.
Therefore,  a rational map is not an actual map in some cases. The merit of using toric compactification is that the boundary structure 
of toric compactification is simpler than the one of stable map compactification. Since the moduli space is different, the correlation 
functions of gauged linear sigma model do not always coincide with the corresponding Gromov-Witten invariants.
Motivated by these facts, our general conjecture is the following.\\
\\
{\bf  General Conjecture}\\
{\it The 2-point correlation functions computed by using the moduli space of gauged linear sigma model give us the information of the B-model used in the 
mirror computation of the Gromov-Witten invariants.  In particular, some 2-point correlation functions give us the expansion coefficients 
of the mirror map used in the mirror computation and the remaining 2-point functions are translated into 2-point Gromov-Witten invariants 
via the (generalized) mirror transformation caused by the mirror map. }  
\\
\\
Of course, the above conjecture is a little bit abstract. For example, we have to define the 2-point correlation function of gauged 
linear sigma model.
 We will give more explicit details in the following part of this section.
Before we turn into details, we remark here that this paper is a continuation of our previous work \cite{vs}, which is our first paper aiming at 
establishing the above conjecture when the toric manifold is  $CP^{N-1}$. 

In \cite{vs}, we proposed a residue integral representation 
of virtual structure constant $\tilde{L}_n^{N,k,d}$, which is a B-model analogue of genus $0$ Gromov-Witten invariants of 
the degree $k$ hypersurface in $CP^{N-1}$ (we denote this hypersurface by $M_{N}^{k}$). $\tilde{L}_n^{N,k,d}$ is our candidate 
of 2-point correlation function of gauged linear sigma model.  The  virtual structure constant $\tilde{L}_n^{N,k,d}$ is the 
rational number which is non-zero if and only if $0\leq n\leq N-1-(N-k)d$. It is defined by the initial condition 
\begin{equation}
\sum_{n=0}^{k-1}\tilde{L}_n^{N,k,1}w^n=k\prod_{j=1}^{k-1}((k-j)+jw),
\end{equation} 
and the recursive formulas that represent $\tilde{L}_n^{N,k,d}$ as a weighted homogeneous polynomial in $\tilde{L}_m^{N+1,k,d^{\prime}}$
($d^{\prime}\leq d$). We will show explicit form of the recursive formulas in Section 2.
Let us first review the main results on the virtual structure constants presented in \cite{gene,gm}. 
For this purpose, we introduce genus $0$ degree $d$ two-point Gromov-Witten invariant $\langle{\cal O}_{h^a}{\cal O}_{h^b}\rangle_{0,d}$
of $M_{N}^{k}$. Here,
$h$ is the cohomology class of $M_{N}^{k}$ induced from hyperplane class of $CP^{N-1}$, and $\langle{\cal O}_{h^a}{\cal O}_{h^b}\rangle_{0,d}$ is defined by the formula:
\begin{equation}
\langle{\cal O}_{h^a}{\cal O}_{h^b}\rangle_{0,d}=\int_{\overline{M}_{0,2}(CP^{N-1},d)}ev_{1}^{*}(h^a)\wedge
ev_{2}^{*}(h^b)\wedge c_{top}(R^{0}\pi_{*}ev_{3}^{*}{\cal O}_{CP^{N-1}}(k)). 
\label{gwdef}
\end{equation}
In (\ref{gwdef}), $\overline{M}_{0,n}(CP^{N-1},d)$ is the moduli space of stable maps of degree $d$ from genus $0$ stable curves with $n$ marked 
points to $CP^{N-1}$. $ev_{i}:\overline{M}_{0,n}(CP^{N-1},d)\rightarrow CP^{N-1}$ is the evaluation map at the i-th marked point.
$\pi:\overline{M}_{0,3}(CP^{N-1},d)\rightarrow \overline{M}_{0,2}(CP^{N-1},d)$ is the forgetful map that forgets the third marked point.
 
If $N-k\geq 1$, i.e., the hypersurface is a Fano manifold, 
we have the following equality:
\begin{equation}
\frac{\tilde{L}_{n}^{N,k,d}}{d}=\displaystyle{\frac{\langle{\cal O}_{h^{N-2-n}}{\cal O}_{h^{n-1+(N-k)d}}\rangle_{0,d}}{k}}, 
\end{equation}
except for $N-k=1$ and $d=1$ case. 

If $N-k=1$ and $d=1$, 
we have an equality:
\begin{eqnarray}
\tilde{L}_n^{k+1,k,1}-k!=\frac{\langle{\cal O}_{h^{k-1-n}}{\cal O}_{h^{n}}\rangle_{0,1}}{k}.
\label{spfano}
\end{eqnarray}
If $N-k\leq 0$, these two numbers differ from each other. In this case, 
$\tilde{L}_{n}^{N,k,d}$ appears as the matrix element of the connection matrix of the virtual Gauss-Manin system \cite{gm}
associated with the Picard-Fuchs differential equation used in 
the mirror computation:
\begin{equation}
\biggl((\d_{x})^{N-1}-k\cdot e^{x}\cdot (k\d_{x}+k-1)(k\d_{x}+k-2)\cdots(k\d_{x}+1)\biggr)w(x)=0.
\label{fun}
\end{equation}
Let us explain the relation between $\tilde{L}_{n}^{N,k,d}$ and (\ref{fun}) more explicitly when $N=k$, i.e., when the hypersurface 
is a Calabi-Yau manifold. A linearly independent basis of solutions of (\ref{fun}) around $x=-\infty$ is given by:
\begin{eqnarray}
u_{j}^{k,k}(x):=\left.\frac{1}{j!}(\d_z)^j\left(\sum_{d=0}^{\infty}\exp((d+z)x)\frac{(kd)!}{(d!)^{N}}\prod_{j=1}^{kd
}(1+\frac{k}{j}z)
\prod_{j=1}^{d}(1+\frac{1}{j}z)^{-N}\right)\right|_{z=0},\;\; (j=0,1,\cdots,k-2).
\label{psol}
\end{eqnarray}
On the other hand, we introduce a generating function of the $\tilde{L}_{n}^{k,k,d}$:
\begin{eqnarray}
\tilde{L}_{n}^{k,k}(e^x):=1+\sum_{d=1}^{\infty}\tilde{L}_{n}^{k,k,d}e^{dx}.
\label{lgen}
\end{eqnarray}
In \cite{gm}, we proved the following equality:
\begin{equation}
u_{j}^{k,k}(x):={\tilde{L}_{0}^{k,k}(e^x)}\int^{x}_{-\infty}dx_{1}
{\tilde{L}_{1}^{k,k}(e^{x_{1}})}\int^{x_{1}}_{-\infty}dx_{2}
{\tilde{L}_{2}^{k,k}(e^{x_{2}})}\cdots\int^{x_{j-1}}_{-\infty}dx_{j}
{\tilde{L}_{j}^{k,k}(e^{x_{j}})},
\label{cy}
\end{equation}
where we apply a formal rule: $\int^{x}_{-\infty}x^mdx=\frac{1}{m+1}x^{m+1}$ in integrating  
the top term of expansion (\ref{lgen}). By using (\ref{cy}), we can represent 
$\tilde{L}_{n}^{k,k}(e^x)$ in terms of $u_{j}^{k,k}(x)$'s. The most important relation derived from (\ref{cy})
is the following equality:
\begin{eqnarray}
x+\sum_{d=1}^{\infty}\frac{\tilde{L}_1^{k,k,d}}{d}e^{dx}=\frac{u_{1}^{k,k}(x)}{u_{0}^{k,k}(x)},
\label{mirrk}
\end{eqnarray}
where the r.h.s. gives us the celebrated mirror map: $t(x)=\frac{u_{1}^{k,k}(x)}{u_{0}^{k,k}(x)}$ used in the mirror computation. 
With this mirror map,
We can compute $\displaystyle{\frac{\langle{\cal O}_{h^{k-2-n}}{\cal O}_{h^{n-1}}\rangle_{0,d}}{k}}$
from the equality:
\begin{eqnarray}
t+\sum_{d=1}^{\infty}\displaystyle{\frac{\langle{\cal O}_{h^{k-2-n}}{\cal O}_{h^{n-1}}\rangle_{0,d}}{k}}e^{dt}
=x(t)+\sum_{d=1}^{\infty}\frac{\tilde{L}_n^{k,k,d}}{d}e^{dx(t)}.
\label{mirrckk}
\end{eqnarray}
This is the mirror transformation caused by the mirror map in (\ref{mirrk}).\\
If $N<k$, we can also compute $\displaystyle{\frac{\langle{\cal O}_{h^{N-2-n}}{\cal O}_{h^{n-1+(N-k)d}}\rangle_{0,d}}{k}}$
by using a generalization of (\ref{mirrk}) and (\ref{mirrckk}) \cite{givc,iri,gene}. In this case, 
$\frac{\tilde{L}_{1+(k-N)d}^{N,k,d}}{d}$ appears in the mirror map as follows:
\begin{equation}
t_{n}=x_{n}+\sum_{d=1}^{\infty}\frac{\tilde{L}_{1+(k-N)d}^{N,k,d}}{d}e^{dx_1}\cdot\delta_{1+(k-N)d,n},\;\; (n=0,1,\cdots,N-2).
\label{genmirmap}
\end{equation}
where $x_{n}$ (resp. $t_{n}$) is the B-model (resp. A-model) deformation parameter associated with ${\cal O}_{h^{n}}$. 
The two-point Gromov-Witten invariants of $M_{N}^{k}$ are obtained after operating the generalized mirror transformation 
caused by (\ref{genmirmap}) on $\tilde{L}_n^{N,k,d}$'s.
Here, we explicitly write down 
the formulas that represent $\displaystyle{\frac{\langle{\cal O}_{h^{N-2-n}}{\cal O}_{h^{n-1+(N-k)d}}\rangle_{0,d}}{k}}$ 
in terms of the virtual structure constants up to $d=3$. They were proved in \cite{prmir}.  
\begin{eqnarray}
\frac{1}{k} \langle{\cal O}_{h^{N-2-n}}{\cal O}_{h^{n-1+N-k}}\rangle_{0,1} &=&\tilde{L}_{n}^{N,k,1}-\tilde{L}_{1+(k-N)}^{N,k,1},\no\\
\frac{1}{k}\langle{\cal O}_{h^{N-2-n}}{\cal O}_{h^{n-1+2(N-k)}}\rangle_{0,2}&=&\frac{1}{2}(\tilde{L}_{n}^{N,k,2}-\tilde{L}_{1+2(k-N)}^{N,k,2})
-\tilde{L}_{1+(k-N)}^{N,k,1}(\sum_{j=0}^{k-N}(\tilde{L}_{n-j}^{N,k,1} 
- \tilde{L}_{1+2(k-N)-j}^{N,k,1})),\no\\
\frac{1}{k}\langle{\cal O}_{h^{N-2-n}}{\cal O}_{h^{n-1+3(N-k)}}\rangle_{0,3} &=&\frac{1}{3}(\tilde{L}_{n}^{N,k,3}-\tilde{L}_{1+3(k-N)}^{N,k,3})
-\tilde{L}_{1+(k-N)}^{N,k,1}(\sum_{j=0}^{k-N}
(\tilde{L}_{n-j}^{N,k,2}-\tilde{L}_{1+3(k-N)-j}^{N,k,2})
+C_{1,1}^{N,k,3}(n))\nonumber\\
&&-\frac{1}{2}\tilde{L}_{1+2(k-N)}^{N,k,2}
(\sum_{j=0}^{2(k-N)}(\tilde{L}_{n-j}^{N,k,1}
-\tilde{L}_{1+3(k-N)-j}^{N,k,1}))\no\\
&&+\frac{3}{2}(\tilde{L}_{1+(k-N)}^{N,k,1})^{2}
(\sum_{j=0}^{2(k-N)}A_{j}(\tilde{L}_{n-j}^{N,k,1}
-\tilde{L}_{1+3(k-N)-j}^{N,k,1})),
\label{th3}
\end{eqnarray}
where
\begin{eqnarray}
A_{j}&:=&j+1,\;\;\mbox{if}\;\;\;(0\leq j\leq k-N),\;\;\;
A_{j}:=1+2(k-N)-j,\;\;\mbox{if}\;\;\; (k-N\leq j\leq 2(k-N)),\no\\
C_{1,1}^{N,k,3}(n)&=&\sum_{j=0}^{(k-N)-1}\bigl(\sum_{m=0}^{j}\tilde{L}_{n-m}^{N,k,1}
\tilde{L}_{n-2(k-N)+j-m}^{N,k,1}-\tilde{L}_{(k-N)+2+j}^{N,k,1}
(\sum_{m=0}^{2(k-N)}\tilde{L}_{n-m}^{N,k,1})\no\\
&&+\tilde{L}_{1+(k-N)}^{N,k,1}
(\sum_{m=j+1}^{2(k-N)-j-1}\tilde{L}_{n-m}^{N,k,1})\bigr)\no\\
&&-\sum_{j=0}^{(k-N)-1}\bigl(\sum_{m=0}^{j}\tilde{L}_{1+3(k-N)-m}^{N,k,1}
\tilde{L}_{1+(k-N)+j-m}^{N,k,1}-\tilde{L}_{(k-N)+2+j}^{N,k,1}
(\sum_{m=0}^{2(k-N)}\tilde{L}_{1+3(k-N)-m}^{N,k,1})\no\\
&&+\tilde{L}_{1+(k-N)}^{N,k,1}
(\sum_{m=j+1}^{2(k-N)-j-1}\tilde{L}_{1+3(k-N)-m}^{N,k,1})\bigr).
\label{fini}
\end{eqnarray}
Now, we go back to the argument given in \cite{vs}. In \cite{vs}, our conjectural residue integral representation 
of $\frac{\tilde{L}_n^{N,k,d}}{d}$ leads us to speculate that if $N-2-n\geq 0$ and $n-1+(N-k)d\geq 0$, 
$\frac{\tilde{L}_n^{N,k,d}}{d}$ can be interpreted as an intersection number on the moduli space of polynomial maps 
with two marked points. 
Let $\widetilde{Mp}_{0,2}(N,d)$ be the compactified moduli space of polynomial maps from $CP^{1}$ to 
$CP^{N-1}$ of degree $d$ with two marked points, which was introduced in \cite{vs} and will be explicitly defined in Section 2 of this paper. 
This space is the moduli space that corresponds to gauged linear sigma model.
We defined an intersection 
number:
\begin{eqnarray}
w({\cal O}_{h^{\alpha}}{\cal O}_{h^{\beta}})_{0,d}:=
\int_{\widetilde{Mp}_{0,2}(N,d)}ev_{1}^{*}(h^{\alpha})\wedge ev_{2}^{*}(h^{\beta})\wedge
c_{top}({\cal E}_{d}^{k}),
\label{defwN}
\end{eqnarray} 
where ${\cal E}_{d}^{k}$ is a rank $kd+1$ orbi-bundle on $\widetilde{Mp}_{0,2}(N,d)$ that corresponds to $R^{0}\pi_{*}ev_{3}^{*}{\cal O}_{CP^{N-1}}(k)$ 
on $\overline{M}_{0,2}(CP^{N-1},d)$: the corresponding moduli space of non-linear sigma model. 
In (\ref{defwN}), $ev_{i}:\widetilde{Mp}_{0,2}(N,d)\rightarrow CP^{N-1}$ is the evaluation map at
 the $i$-th marked point. We computed 
$w({\cal O}_{h^{\alpha}}{\cal O}_{h^{\beta}})_{0,d}$ by localization techniques and concluded that our residue integral 
representation suggests,
\begin{eqnarray}
k\cdot\frac{\tilde{L}_n^{N,k,d}}{d}= w({\cal O}_{h^{N-2-n}}{\cal O}_{h^{n-1+(N-k)d}})_{0,d}.
\label{keyvs}
\end{eqnarray}
In Section 2 of this paper, we prove,
\begin{theorem}
(\ref{keyvs}) is true for arbitrary $M_{N}^{k}$ if $0\leq N-2-n\leq N-2$ and if $0\leq n-1+(N-k)d\leq N-2$.
\end{theorem}

At this stage, we go back to the equality (\ref{mirrk}). We introduce here the classical three-point function and metric,
\begin{eqnarray}
w({\cal O}_{h^{\alpha}}{\cal O}_{h^{\beta}}{\cal O}_{h^{\gamma}})_{0,0}&:=&\int_{CP^{k-1}}kh\wedge h^{\alpha}\wedge 
h^{\beta}\wedge h^{\gamma}=k\cdot\delta_{\alpha+\beta+\gamma,k-2},\no\\
\eta_{\alpha\beta}&:=&w({\cal O}_{h^{\alpha}}{\cal O}_{h^{\beta}}{\cal O}_{h^{0}})_{0,0}=
k\cdot\delta_{\alpha+\beta,k-2},\no\\
\eta^{\alpha\beta}&:=&\frac{1}{k}\cdot\delta_{\alpha+\beta,k-2}.
\end{eqnarray}
We also introduce perturbed two-point functions:
\begin{eqnarray}
w({\cal O}_{h^{\alpha}}{\cal O}_{h^{\beta}})_{0,0}(x)&:=&w({\cal O}_{h^{\alpha}}{\cal O}_{h^{\beta}}{\cal O}_{h})_{0,0}\cdot x,\no\\
w({\cal O}_{h^{\alpha}}{\cal O}_{h^{\beta}})_{0,d}(x)&:=&w({\cal O}_{h^{\alpha}}{\cal O}_{h^{\beta}})_{0,d}\cdot e^{dx}\;\;(d\geq 1).
\label{pertN}
\end{eqnarray}
With this setup, we can conclude from (\ref{mirrk}) and (\ref{keyvs}) that the equality:
\begin{eqnarray}
t(x)=\eta^{1\alpha}\left(\sum_{d=0}^{\infty}w({\cal O}_{h^{\alpha}}{\cal O}_{h^{0}})_{0,d}(x)\right),
\label{genemirr}
\end{eqnarray}
gives us the mirror map used in the mirror computation. One of our motivations in this paper is to generalize (\ref{genemirr}) 
to the mirror computation of toric manifolds with two K\"ahler forms.
In this paper, we consider Hirzebruch surfaces $F_{0}$, $F_{3}$ and resolution of weighted projective space 
${\bf P}(1,1,2,2,2)$ (we denote it by $WP_1$) as examples. These toric manifolds have two K\"ahler forms.
Let $z$ and $w$ be these two K\"ahler forms. Polynomial maps from $CP^1$ to these toric manifolds are classified 
by bi-degree:  
\begin{equation}
\dd:=(d_a,d_b),\;\;(\dd\neq (0,0)),
\label{bid}
\end{equation}
where $d_a$ and $d_b$ are non-negative integers.
Let $\widetilde{Mp}_{0,2}(X,\dd)$ be  the compactified moduli space of polynomial maps from $CP^1$ to $X$ with two marked 
points of degree $\dd$, which can be constructed by generalizing the construction of $\widetilde{Mp}_{0,2}(N,d)$.
Of course, $X$ considered here is  $F_{0}$ or $F_{3}$ or $WP_1$. Then we consider the following intersection numbers 
on $\widetilde{Mp}_{0,2}(X,\dd)$:
\begin{eqnarray}
w({\cal O}_{\alpha}{\cal O}_{\beta})_{0,\dd}&:=&
\int_{\widetilde{Mp}_{0,2}(F_{0},\dd)}ev_{1}^{*}(\alpha)\wedge ev_{2}^{*}(\beta)\wedge
c_{top}({\cal E}_{\dd}),\no\\
w({\cal O}_{\alpha}{\cal O}_{\beta})_{0,\dd}&:=&
\int_{[\widetilde{Mp}_{0,2}(F_{3},\dd)]_{ver.}}ev_{1}^{*}(\alpha)\wedge ev_{2}^{*}(\beta),\no\\
w({\cal O}_{\alpha}{\cal O}_{\beta})_{0,\dd}&:=&
\int_{[\widetilde{Mp}_{0,2}(WP_1,\dd)]_{ver.}}ev_{1}^{*}(\alpha)\wedge ev_{2}^{*}(\beta)\wedge
c_{top}({\cal E}_{\dd}),
\label{wdef}
\end{eqnarray}
where $\alpha$ and $\beta$ are elements of $H^{*}(X,{\bf C})$. ${\cal E}_{\dd}$ in the first (resp. the third) line of 
(\ref{wdef}) is an orbi-bundle on $\widetilde{Mp}_{0,2}(F_{0},\dd)$ (resp. $\widetilde{Mp}_{0,2}(WP_1,\dd)$) that corresponds 
to $R^{1}\pi_{*}ev_{3}^{*}K_{F_{0}}$ (resp. $R^{0}\pi_{*}ev_{3}^{*}K_{WP_1}^{*}$) on $\overline{M}_{0,2}(F_{0},\dd)$
 (resp. $\overline{M}_{0,2}(WP_1,\dd)$).  
These intersection numbers are analogues of two-point Gromov-Witten invariants of 
$K_{F_0}$, $F_{3}$ and the Calabi-Yau hypersurface in $WP_1$ respectively. 
In this paper, we derive closed formulas to compute these intersection numbers by applying the localization theorem.
The resulting formulas are written as a sum of contributions from connected components of fixed point sets labeled 
by ordered partitions $\sigma_{\dd}$ of bi-degree $\dd$:
\begin{equation}
\sigma_{\dd}=(\dd_{1},\dd_{2},\cdots,\dd_{l(\sigma_{\dd})}),
\:\;(\sum_{j=1}^{l(\sigma_{d})}\dd_{j}=\dd\;\;,\;\;\dd_{j}=(d_{a,j},0)\;\mbox{or}\;(0,d_{b,j}),\;\;d_{a,j},d_{b,j}>0).
\label{part0} 
\end{equation}
This structure can be regarded as a natural generalization of the $CP^{N-1}$ case, because in the $CP^{N-1}$ case, 
$w({\cal O}_{h^{\alpha}}{\cal O}_{h^{\beta}})_{0,d}$ is written as sum of contributions labeled by ordered 
partitions of positive integers $d$. With these formulas, we numerically compute $w({\cal O}_{\alpha}{\cal O}_{\beta})_{0,\dd}$
by using MAPLE for low degrees. For the special cases $F_{0}$ and $WP_1$, we also compute classical intersection numbers, 
metrics and perturbed two-point 
functions by introducing deformation parameters $x_1$ and $x_2$ associated with $z$ and $w$ respectively. 
\begin{eqnarray}
w({\cal O}_{\alpha}{\cal O}_{\beta})_{0,(0,0)}(x_1,x_2)&=&w({\cal O}_{\alpha}{\cal O}_{\beta}{\cal O}_{z})_{0,(0,0)}\cdot x_1+
w({\cal O}_{\alpha}{\cal O}_{\beta}{\cal O}_{w})_{0,(0,0)}\cdot x_2,\no\\
w({\cal O}_{\alpha}{\cal O}_{\beta})_{0,\dd}(x_1,x_2)&=&w({\cal O}_{\alpha}{\cal O}_{\beta})_{0,\dd}e^{d_{a}x_{1}+d_{2}x_{2}},\;\;
({\bf d}\neq(0,0)).
\label{pu2}
\end{eqnarray}
With
this setup, we test whether the equalities:
\begin{eqnarray}
t_1(x_1,x_2)&=&\eta^{z\alpha}\left(\sum_{\dd\geq{\bf 0}}^{\infty}w({\cal O}_{\alpha}{\cal O}_{1})_{0,\dd}(x_1,x_2)\right),\no\\
t_2(x_1,x_2)&=&\eta^{w\alpha}\left(\sum_{\dd\geq{\bf 0}}^{\infty}w({\cal O}_{\alpha}{\cal O}_{1})_{0,\dd}(x_1,x_2)\right),
\label{intromirr}
\end{eqnarray}
give us the mirror map of $K_{F_0}$ and the Calabi-Yau hypersurface in $WP_1$. The numerical results confirm our speculation. 
Therefore, we conjecture that (\ref{intromirr}) indeed gives us the mirror map used in the mirror computation. This conjecture explains 
the meaning of the title of this paper. As in the $CP^{N-1}$ case, we can also compute the standard two-point Gromov-Witten invariants 
$\langle{\cal O}_{\alpha}{\cal O}_{\beta}\rangle_{0,\dd}$ by using $w({\cal O}_{\alpha}{\cal O}_{\beta})_{0,\dd}$ 
 by generalizing the equality (\ref{mirrckk}).  
 In sum, we propose the following conjecture:
 \begin{conj}
 In the case of $K_{F_{0}}$ (resp.  Calabi-Yau hypersurface of $WP_1$), (\ref{intromirr}) gives us the mirror map used in the mirror computation 
 of Gromov-Witten invariants, and 
 \begin{eqnarray}
 &&\sum_{\dd\geq{\bf 0}}^{\infty}w({\cal O}_{\alpha}{\cal O}_{\beta})_{0,\dd}(x_1(t_1,t_2),x_2(t_1,t_2))=\no\\
 &&\langle{\cal O}_{\alpha}{\cal O}_{\beta}{\cal O}_{z}\rangle_{0,(0,0)}t_1+
 \langle{\cal O}_{\alpha}{\cal O}_{\beta}{\cal O}_{w}\rangle_{0,(0,0)}t_2+
\sum_{\dd>{\bf 0}}^{\infty}\langle{\cal O}_{\alpha}{\cal O}_{\beta}\rangle_{0,\dd}e^{d_a t_1+d_b t_2}
 \end{eqnarray}
 where $\langle{\cal O}_{\alpha}{\cal O}_{\beta}\rangle_{0,{\bf d}}$ is the two-point Gromov-Witten invariant of $K_{F_{0}}$ (resp.  Calabi-Yau hypersurface of $WP_1$).
 \end{conj}

Let us turn into the non-nef example $F_{3}$. In this case, we first review Givental-Coates-Guest-Iritani's approach \cite{givc,guest,iri} of the mirror computation of 
Gromov-Witten invariants of non-nef toric manifolds by taking $F_3$ as an eaxample. In this approach, we 
start from the Givental's $I$-function: 
\begin{eqnarray}
	I_{F_3}=e^{(zx_1 + wx_2)/\hbar}\sum_{\dd}\frac{\prod_{m=-\infty}^0(-3z+w+m\hbar)}{\prod_{m=-\infty}^{-3d_a+d_b}(-3z+w+m\hbar)\prod_{m=1}^{d_a}(z+m\hbar)^2\prod_{m=1}^{d_b}(w+m\hbar)}e^{d_ax_1+d_bx_2},
\label{If3intro}
\end{eqnarray}
where $I_{F_3}$ is the cohomology-valued function.
Note that $I_{F_3}$ contains the parameter $\hbar$, which plays a central role in Givental-Coates-Guest-Iritani's approach.
We take $1,z,w,w^2$ as the basis of $H^{*}(F_{3},{\bf C})$ and expand cohomology-valued function $F$ into the form:
\begin{equation}
F=F(1)\cdot 1+F(z)\cdot z+F(w)\cdot w+F(w^2)\cdot w^2.
\end{equation}  
Next, we define the $4\times 4$ matrix $S$ given by,
\begin{eqnarray}
{}^{t}S=\left(\begin{array}{cccc} I_{F_{3}}(1)&\bigl((\hbar\d_{x_1}+z )I_{F_{3}}\bigr)(1)&\bigl((\hbar\d_{x_2}+w )I_{F_{3}}\bigr)(1)
&\bigl((\hbar\d_{x_2}+w )^2I_{F_{3}}\bigr)(1)\\
 I_{F_{3}}(z)&\bigl((\hbar\d_{x_1}+z )I_{F_{3}}\bigr)(z)&\bigl((\hbar\d_{x_2}+w )I_{F_{3}}\bigr)(z)
&\bigl((\hbar\d_{x_2}+w )^2I_{F_{3}}\bigr)(z)\\
 I_{F_{3}}(w)&\bigl((\hbar\d_{x_1}+z )I_{F_{3}}\bigr)(w)&\bigl((\hbar\d_{x_2}+w )I_{F_{3}}\bigr)(w)
&\bigl((\hbar\d_{x_2}+w )^2I_{F_{3}}\bigr)(w)\\
 I_{F_{3}}(w^2)&\bigl((\hbar\d_{x_1}+z )I_{F_{3}}\bigr)(w^2)&\bigl((\hbar\d_{x_2}+w )I_{F_{3}}\bigr)(w^2)
&\bigl((\hbar\d_{x_2}+w )^2I_{F_{3}}\bigr)(w^2)
\label{solm}
 \end{array} \right),
\end{eqnarray}
and the connection matrices $\Omega_{z}\;\;\Omega_{w}$:
\begin{equation}
\Omega_{z}=\bigl((\hbar\d_{x_1}+z)S\bigr)S^{-1},\;\;\Omega_{w}=\bigl((\hbar\d_{x_2}+w)S\bigr)S^{-1},\;\;
\end{equation} 
Since $F_3$ is a non-nef manifold, expansion of $S$ around $\hbar=0$ includes both negative and positive powers of $\hbar$. We then 
take Birkhoff factorization of $S=S(\hbar,\hbar^{-1})$ with respect to $\hbar$:
\begin{equation}
S(\hbar,\hbar^{-1})=Q(\hbar)R(\hbar^{-1}).
\end{equation} 
The positive part $Q(\hbar)$ provides a gauge transformation which converts $\Omega_z$ and $\Omega_w$ into $\hbar$ independent 
matrices $B_{z}$ and $B_{w}$:
\begin{equation}
B_z=Q^{-1}(\hbar)\Omega_z Q(\hbar)+\hbar(\d_{x_1}Q^{-1}(\hbar))Q(\hbar),\;\;
B_w=Q^{-1}(\hbar)\Omega_w Q(\hbar)+\hbar(\d_{x_2}Q^{-1}(\hbar))Q(\hbar).
\label{bmat}
\end{equation}
Let us identify the subscript $1,2,3,4$ of $4\times4$ matirices with the cohomology elements $1,z,w,w^2$ respectively. 
We then introduce the clasical metric $\eta=(\eta_{ij})$ of $H^{*}(F_{3},{\bf C})$ and define the matrices: 
\begin{equation}
C_{z}=B_z\eta,\;\;C_w=B_w\eta.
\label{cmat}
\end{equation}  
These are the intermediate results in the mirror computation of the Gromov-Witten invariants of $F_{3}$. In order 
to obtain the three-point Gromov-Witten invariants of $F_{3}$, we have to operate the generalized mirror transformation induced 
from the mirror map: $t_{\alpha}=t_{\alpha}(x_1,x_{2}) $ to $C_{z}$ and $C_{w}$. The mirror map is determined from the matrix 
elements of the $\hbar$
independent connection matrices in (\ref{cmat}) via the relation: 
\begin{equation}
\frac{\d t_{\alpha}}{\d x_{1}}=(C_{z})_{1\gamma}\eta^{\gamma\alpha},\;\;
\frac{\d t_{\alpha}}{\d x_{2}}=(C_{w})_{1\gamma}\eta^{\gamma\alpha}, \;\;(\alpha,\beta,\gamma\in \{1,z,w,w^2\},\;\;
\eta_{\alpha\beta}\eta^{\beta\gamma}=\delta_{\alpha}^{\gamma}\;).
\label{jacobi}
\end{equation}
This final step requires a lot of computations and results in a generalization of the formula (\ref{th3}) to the case of $F_{3}$. 
See \cite{fj1} for details. Let us remark one important point on Givental-Coates-Guest-Iritani's approach applied to $M_{N}^{k}$
with $k>N$.
In this case, the virtual structure constant $\tilde{L}_{n}^{N,k,d}$ is nothing but the matrix element of the $\hbar$ independent 
connection matrix $B_{h}$ that appears in the step (\ref{bmat}).

With these results in mind, we look back at the $w({\cal O}_{\alpha}{\cal O}_{\beta})_{0,\dd}$'s for $F_{3}$. 
The matrix element $(C_{\alpha})_{\beta\gamma}$ in (\ref{cmat}) is a power series in $e^{x_1}$ and $e^{x_{2}}$, 
and we denote the coefficient of $e^{d_a x_1+d_b x_2}$ of $(C_{\alpha})_{\beta\gamma}$ by $(C_{\alpha})_{\beta\gamma}({\bf d})$.
Our conjecture in the case of $F_{3}$ is the following:
\begin{conj}
\begin{equation}
(C_{z})_{\alpha\beta}({\bf d})=d_a\cdot w({\cal O}_{\alpha}{\cal O}_{\beta})_{0,\dd},\;\;
 (C_{w})_{\alpha\beta}({\bf d})=d_b\cdot w({\cal O}_{\alpha}{\cal O}_{\beta})_{0,\dd},
\end{equation}
and the mirror map used in the generalized mirror transformation is given by,
\begin{equation}
t_{\alpha}=x_{\alpha}+\sum_{{\bf d}\neq (0,0)}w({\cal O}_{1}{\cal O}_{\beta})_{0,\dd}\eta^{\beta\alpha}e^{d_a x_1+d_b x_2},
\end{equation}
where $x_{\alpha}$ (resp. $t_{\alpha}$) is the B-model (resp. A-model) deformation parameter associated with ${\cal O}_{\alpha}$
and $x_z$ (resp. $x_{w}$) corresponds to $x_1$ (resp. $x_{2}$).
\end{conj}
In this paper, we compute $w({\cal O}_{\alpha}{\cal O}_{\beta})_{0,\dd}$ for lower $\dd$ by using the definition (\ref{wdef}) 
and the localization theorem.
Our numerical results agree with the numerical data of $C_{z}$ and $C_{w}$ computed in \cite{fj1}. 
Since these connection matrices are enough for the mirror 
computation of Gromov-Witten invariants of $F_{3}$, our formula to compute $w({\cal O}_{\alpha}{\cal O}_{\beta})_{0,\dd}$
gives us another way of carrying out the mirror computation without using Birkhoff factorization.

 Our results in this paper compute 
nothing new from the point of view of the mirror computation, but our construction gives concrete geometrical 
footing to the B-model data as intersection numbers on the moduli space of polynomial maps (or gauged linear sigma model), which can be regarded  as an alternate 
compactification of the moduli space of holomorphic maps from $CP^1$ to a toric manifold. The examples treated in this 
paper imply that our construction can be generalized to arbitrary toric manifolds.  

In this paper, we also give supplemental discussions on our arguments given in \cite{vs}. 
Especially, we present the explicit construction of $\widetilde{Mp}_{0,2}(N,d)$, which was briefly outlined in our previous paper \cite{vs}.
We propose that $\widetilde{Mp}_{0,2}(N,d)$ is given as a toric variety whose weight matrix of ${\bf C}^{\times}$ actions 
includes the $A_{d-1}$ Cartan matrix. This construction explains not only the structure of the boundary components of 
$\widetilde{Mp}_{0,2}(N,d)$, 
but also the reason why expressions associated with the $A_{d-1}$ Cartan matrix appear in the definition of the virtual
structure constant $\tilde{L}_n^{N,k,d}$. We also give a detailed construction of $\widetilde{Mp}_{0,2}(F_0,\dd)$ as a
toric variety. It plays an important role in carrying out the localization computation of $w({\cal O}_{\alpha}{\cal O}_{\beta})_{0,\dd}$
for $F_0$, $F_3$ and $WP_1$. 

In the last part of this paper, we extend our construction to the mirror computation of the K3 surface in the weighted projective
space ${\bf P}(1,1,1,3)$. It is well-known that the mirror map in this example is written by using the elliptic $j$-function. 
Combining this fact with our conjecture, we propose a formula that expresses Fourier expansion coefficients of the $j$-function 
in terms of intersection numbers on $\widetilde{Mp}_{0,2}({\bf P}(1,1,1,3),d)$.

This paper is organized as follows.  

In Section 2, 
we reconsider the argument given in our previous paper \cite{vs} and discuss problems that remained 
unsolved. First, we explicitly construct $\widetilde{Mp}_{0,2}(N,d)$ used in \cite{vs} as a toric  variety.
In this construction, we emphasize that it is obtained from compactifying the moduli space of polynomial 
maps from $CP^1$ to $CP^{N-1}$ with two marked points.  
 We also discuss a problem that is related to the so-called point-instanton, which is included in $\widetilde{Mp}_{0,2}(N,d)$
but excluded in the moduli space of stable maps.  Next, we define the intersection number on $\widetilde{Mp}_{0,2}(N,d)$
that corresponds to the two-point Gromov-Witten invariant of $M_{N}^{k}$ and compute it explicitly by the localization theorem. 
Lastly, we prove Theorem 1 by deriving explicitly the residue integral representation of $\frac{\tilde{L}_n^{N,k,d}}{d}$.

In Section 3, we generalize the localization computation of 
intersection numbers on polynomial maps to toric manifolds with two K\"ahler forms. First, we take the Hirzebruch 
surface $F_{0}={\bf P}^1\times{\bf P}^1$ and construct $\widetilde{Mp}_{0,2}(F_0,\dd)$ as a toric 
variety. Next, we define intersection numbers on $\widetilde{Mp}_{0,2}(F_0,\dd)$ that correspond to local 
Gromov-Witten invariants of $K_{F_0}$ and derive closed formulas for them by using the localization theorem. 
We then give some numerical results on these intersection numbers and use these to carry out the mirror computation of $K_{F_{0}}$. We take the non-nef Hirzebruch surface $F_3$ as our next example. We assume that 
$\widetilde{Mp}_{0,2}(F_3,\dd)$ has the same boundary structure as $\widetilde{Mp}_{0,2}(F_0,\dd)$ and compute 
intersection numbers on $\widetilde{Mp}_{0,2}(F_3,\dd)$ that correspond to Gromov-Witten invariants of $F_{3}$. We 
show that our numerical results coincide with the expansion coefficients of matrix elements of connection
matrices obtained from  Birkhoff factorization of the Givental $I$-function of $F_3$. Our last example in this section is the 
resolution of weighted projective space ${\bf P}(1,1,2,2,2)$, which we call $WP_1$. We define intersection numbers on 
$\widetilde{Mp}_{0,2}(WP_1,\dd)$ that correspond to Gromov-Witten invariants of Calabi-Yau hypersurface in $WP_1$ 
and compute them by the localization theorem. We end this section by demonstrating the mirror computation of the Calabi-Yau 
hypersurface by using numerical data of the intersection numbers.

In Section 4, we extend our computation to the K3 surface in weighted projective space ${\bf P}(1,1,1,3)$. 
This example is well-known because the mirror map of it is closely related with the elliptic $j$-function. We show numerically that the expansion coefficients of the mirror map are given by intersection numbers on 
 $\widetilde{Mp}_{0,2}({\bf P}(1,1,1,3),d)$. Next, we present a formula which expresses the Fourier coefficients 
of the elliptic $j$-function in terms of these intersection numbers. Lastly, we mention a resolution of the weighted 
projective space ${\bf P}(1,1,2,2,6)$, which can be regarded as a ${\bf P}(1,1,1,3)$ bundle over ${\bf P}^1$.
\\
\\ 
{\bf Notation} Throughout this paper, we denote by $\frac{1}{2\pi\sqrt{-1}}\oint_{C_a}dz$ the operation of taking residue at $z=a$. 
If we write $$\frac{1}{2\pi\sqrt{-1}}\oint_{C_{(a_1,a_2,\cdots,a_m)}}dz$$, it means taking residues at $z=a_j,\;(j=1,2,\cdots,m)$.    
\\
\\
{\bf Acknowledgment} The author would like to thank Prof. T.Eguchi and Dr. B.Forbes for valuable discussions.
He especially would like to thank Dr. B.Forbes for correcting English of this paper.  
He would also like to thank organizers of the workshop "Branes, Strings and Black Holes" at the Yukawa Institute of 
Theoretical Physics, during which part of this work was done.
Lastly, he would like to thank Miruko Jinzenji for kind encouragement. His research is partially supported 
by JSPS grant No. 22540061.  
\newpage
\section{$CP^{N-1}$ case Revisited}
\subsection{Review of the Results in the case of $CP^{N-1}$}
\subsubsection{Toric Compactification of the Moduli Space of Degree $d$ Polynomial Maps with Two Marked Points}
Let ${\bf a}_{j},\;(j=0,1,\cdots,d)$ be vectors in ${\bf C}^{N}$ 
and let $\pi_{N}: {\bf C}^{N}\rightarrow CP^{N-1}$ be a projection map. 
In this paper, we define a degree $d$ polynomial map $p$ from ${\bf C}^{2}$ to ${\bf C}^{N}$ as a map that consists of 
${\bf C}^{N}$vector-valued degree $d$ homogeneous polynomials in two coordinates $s,t$ of ${\bf C}^{2}$:    
\begin{eqnarray}
&&p:{\bf C}^{2}\rightarrow {\bf C}^{N}\no\\
&&p(s,t)={\bf a}_{0}s^{d}+{\bf a}_{1}s^{d-1}t+{\bf a}_{2}s^{d-2}t^{2}+\cdots+{\bf a}_{d}t^{d}.
\label{polyp1}
\end{eqnarray}
The parameter space of polynomial maps is given by ${\bf C}^{N(d+1)}=\{ ({\bf a}_{0},{\bf a}_{1},\cdots,{\bf a}_{d}) \}$.
We denote by $Mp_{0,2}(N,d)$ the space obtained from dividing $\{({\bf a}_{0},\cdots,{\bf a}_{d})\in
{\bf C}^{N(d+1)}|\;{\bf a}_{0}\neq {\bf 0},\;{\bf a}_{d}\neq {\bf 0}\}$ by two ${\bf C}^{\times}$ actions induced from the following 
two ${\bf C}^{\times}$ actions on ${\bf C}^{2}$ via the map $p$ in (\ref{polyp1}). 
\begin{equation}
(s,t)\rightarrow ( \mu s,\mu t),\;\;(s,t)\rightarrow (s,\nu t). 
\label{torus1}
\end{equation}
With the above two torus actions, $Mp_{0,2}(N,d)$ can be regarded as the parameter space of degree $d$ rational maps from $CP^{1}$ to $CP^{N-1}$ with
two marked points in $CP^{1}$: $0(=(1:0))$ and $\infty(=(0:1))$ . Set theoretically, it is given as follows:
\begin{eqnarray}
Mp_{0,2}(N,d) = \{ ( {\bf a}_{0},{\bf a}_{1},\cdots,{\bf a}_{d}) \in CP^{N(d+1)}\;|\;{\bf a}_{0},{\bf a}_{d}\neq {\bf 0}\}/({\bf C}^{\times})^2,
\end{eqnarray}
where the two ${\bf C}^{\times}$actions are given by, 
\begin{eqnarray}
&&( {\bf a}_{0},{\bf a}_{1},\cdots,{\bf a}_{d}) \rightarrow (\mu{\bf a}_{0},\mu{\bf a}_{1},\cdots,\mu{\bf a}_{d-1},\mu{\bf a}_{d})\no\\
&&( {\bf a}_{0},{\bf a}_{1},\cdots,{\bf a}_{d}) \rightarrow ({\bf a}_{0},\nu{\bf a}_{1},\cdots,\nu^{d-1}{\bf a}_{d-1},\nu^{d}{\bf a}_{d})
\label{two22}
\end{eqnarray} 
The condition ${\bf a}_{0}, {\bf a}_{d}\neq {\bf 0}$ assures that the images of $0$ and $\infty$ are well-defined in $CP^{N-1}$. 

At this stage, we have to note the difference between the moduli space of holomorphic maps from $CP^{1}$ to $CP^{N-1}$ and the moduli
space of polynomial maps from $CP^{1}$ to $CP^{N-1}$. In short, the latter includes the points that are not the actual maps from 
$CP^{1}$ to $CP^{N-1}$ but the rational maps from $CP^{1}$ to $CP^{N-1}$. These points are called point instantons by physicists. 
More explicitly, a point instanton is a polynomial map 
$\sum_{j=0}^{d}{\bf a}_js^{j}t^{d-j}$ which can be factorized as 
\begin{eqnarray}
\sum_{j=0}^{d}{\bf a}_js^{j}t^{d-j}=p_{d-d_1}(s,t)\cdot(\sum_{j=0}^{d_1}{\bf b}_js^{j}t^{d_1-j}),
\label{pfacin}
\end{eqnarray}
where $p_{d-d_1}(s,t)$ is a homogeneous polynomial of degree $d-d_1(>0)$. If we consider $\sum_{j=0}^{d}{\bf a}_js^{j}t^{d-j}$
as a map from $CP^{1}$ to $CP^{N-1}$, it should be regarded as a rational map whose images of the zero points of $p_{d-d_1}$ is 
undefined.　Moreover, the closure of the image of this map is a rational curve of degree $d_1(<d)$ in $CP^{N-1}$. The reason why 
we include point instantons is that we can obtain simpler compactification of the moduli space than the moduli space of the stable maps 
$\overline{M}_{0,2}(CP^{N-1},d)$, the standard moduli space used to define the two-point Gromov-Witten invariants.  

Now, let us turn into the problem of compactification of $Mp_{0,2}(N,d)$. If $d=1$, $Mp_{0,2}(N,1)$ is given by,  
\begin{eqnarray}
Mp_{0,2}(N,1) = \{ ( {\bf a}_{0},{\bf a}_{1}) \in CP^{N(d+1)}\;|\;{\bf a}_{0},{\bf a}_{1}\neq {\bf 0}\}/({\bf C}^{\times})^2,
\end{eqnarray}
where $({\bf C}^{\times})^2$ action is given as follows.
\begin{eqnarray}
&&( {\bf a}_{0},{\bf a}_{1}) \rightarrow (\mu{\bf a}_{0},\mu{\bf a}_{1})\no\\
&&( {\bf a}_{0},{\bf a}_{1}) \rightarrow ({\bf a}_{0},\nu{\bf a}_{1}).
\label{two221}
\end{eqnarray} 
Therefore, $Mp_{0,2}(N,1)$ is nothing but $CP^{N-1}\times CP^{N-1}$ and is already compact. 
If $d\geq 2$, we have to use the two ${\bf C}^{\times}$ actions in (\ref{two22}) to turn ${\bf a}_0$ and ${\bf a}_d$ into the 
points in $CP^{N-1}$, $[{\bf a}_0]$ and $[{\bf a}_d]$. Therefore, we can easily see,
\begin{eqnarray}
Mp_{0,2}(N,d) = \{ ( [{\bf a}_{0}],{\bf a}_{1},\cdots,{\bf a}_{d-1},[{\bf a}_{d}]) \in CP^{N-1}\times {\bf C}^{N(d-1)}\times CP^{N-1}\;|\}/{\bf Z}_{d}.
\label{openmp}
\end{eqnarray} 
In (\ref{openmp}), the ${\bf Z}_d$ acts on ${\bf C}^{N(d-1)}$ as follows.
\begin{equation}
({\bf a}_{1},{\bf a}_{2}\cdots,{\bf a}_{d-1})\rightarrow ((\zeta_{d})^j{\bf a}_{1},(\zeta_{d})^{2j}{\bf a}_{2}\cdots,
(\zeta_{d})^{(d-1)j}{\bf a}_{d-1}), 
\end{equation} 
where $\zeta_d$ is the d-th primitive root of unity.  In this way, we can see that $Mp_{0,2}(N,d)$ is not compact if $d\geq 2$.
In order to compactify $Mp_{0,2}(N,d)$, we imitate the stable map compactification and add the following chains of polynomial maps 
\begin{equation}
\cup_{j=1}^{l(\sigma_{d})}\bigl(\sum_{m_{j}=0}^{d_{j}-d_{j-1}} {\bf a}_{d_{j-1}+m_{j}}(s_{j})^{m_{j}}(t_{j})^{d_{j}-d_{j-1}-m_{j}}\bigr),
\;\;\bigl({\bf a}_{d_{j}}\neq {\bf 0},\;\;j=0,1,\cdots,l(\sigma_{d})\bigr),
\label{chain00}
\end{equation}
at the infinity locus of $Mp_{0,2}(N,d)$. In (\ref{chain00}), $d_{j}$'s are integers that satisfy,
\begin{equation}
1\leq d_1<d_2<\cdots<d_{l(\sigma_{d})}\leq d-1.
\end{equation}
We denote by $\widetilde{Mp}_{0,2}(N,d)$ the space 
obtained after this compactification. 
This $\widetilde{Mp}_{0,2}(N,d)$ is the moduli space we use in this paper.
It is explicitly constructed as a toric orbifold by introducing boundary divisor coordinates $u_1,u_2,\cdots 
u_{d-1}$ as follows. 
\begin{eqnarray}
&&\widetilde{Mp}_{0,2}(N,d) = \no\\
&&\{ ( {\bf a}_{0},{\bf a}_{1},\cdots,{\bf a}_{d},u_1,u_2,\cdots,u_{d-1}) \in CP^{N(d+1)+d-1}\;|\;{\bf a}_{0},
({\bf a}_{1},u_{1}),\cdots,({\bf a}_{d-1},u_{d-1}),{\bf a}_{d}\neq {\bf 0}\}/({\bf C}^{\times})^{d+1},\no\\
\end{eqnarray}
where the (d+1) ${\bf C}^{\times}$actions are given by,
\begin{eqnarray}
&&( {\bf a}_{0},{\bf a}_{1},\cdots,{\bf a}_{d},u_1,\cdots,u_{d-1}) \rightarrow (\mu_{0}{\bf a}_{0},\cdots,\mu_{0}^{-1}u_{1},
\cdots),\no\\
&& ( {\bf a}_{0},{\bf a}_{1},\cdots,{\bf a}_{d},u_1,\cdots,u_{d-1}) \rightarrow (\cdots,\mu_{1}{\bf a}_{1},\cdots,\mu_{1}^{2}
u_{1},\mu_{1}^{-1}u_{2},\cdots), \no\\
&& ( {\bf a}_{0},{\bf a}_{1},\cdots,{\bf a}_{d},u_1,\cdots,u_{d-1}) \rightarrow (\cdots,\mu_{i}{\bf a}_{i},\cdots,\mu_{i}^{-1}u_{i-1},\mu_{i}^{2}u_{i},
\mu_{i}^{-1}u_{i+1},\cdots),\;(i=2,\cdots,d-1),\no\\ 
&& ( {\bf a}_{0},{\bf a}_{1},\cdots,{\bf a}_{d},u_1,\cdots,u_{d-1}) \rightarrow (\cdots,\mu_{d-1}{\bf a}_{d-1},\cdots,\mu_{d-1}^{-1}u_{d-2},
\mu_{d-1}^{2}u_{d-1}
),\no\\ 
&& ( {\bf a}_{0},{\bf a}_{1},\cdots,{\bf a}_{d},u_1,\cdots,u_{d-1}) \rightarrow (\cdots,\mu_{d}{\bf a}_{d},\cdots,\mu_{d}^{-1}u_{d-1}). 
\label{action1}
\end{eqnarray}
In(\ref{action1}), "$\cdots$" in the r.h.s indicates that the ${\bf C}^{\times}$ actions are trivial.   
These torus actions are represented by a $(d+1)\times 2d$ weight matrix $W_{d}$:
\begin{eqnarray} 
W_{d}:=\bordermatrix{                    &{\bf a}_{0}&{\bf a}_{1}&{\bf a}_{2}&\cdots&{\bf a}_{d-3}&{\bf a}_{d-2}&{\bf a}_{d-1}&{\bf a}_{d}&u_{1}&
                                  u_{2}&u_{3}&\cdots&u_{d-2}&u_{d-1}\cr
                                  h_{0}&1&0&0&\cdots&0&0&0&0&-1&0&0&\cdots&0&0\cr 
                                  h_{1}&0&1&0&\cdots&0&0&0&0&2&-1&0&\cdots&0&0\cr  
                                  h_{2}&0&0&1&\ddots&0&\vdots&0&0&-1&2&-1&\ddots&0&0\cr   
                                  \vdots&\vdots&\vdots&\ddots&\ddots&\ddots&\vdots&\vdots&\vdots&\vdots&\ddots&\ddots&\ddots&\vdots&\vdots\cr 
                                  \vdots&\vdots&\vdots&0&\ddots&1&0&0&0&0&0&\ddots&\ddots&\ddots&0\cr  
                                  \vdots&\vdots&\vdots&\vdots&\ddots&0&1&0&0&0&0&\ddots&-1&2&-1\cr    
                                  h_{d-1}&0&0&0&\cdots&0&0&1&0&0&0&\cdots&0&-1&2\cr
                                  h_{d}&0&0&0&\cdots&0&0&0&1&0&0&\cdots&0&0&-1\cr}
\label{toric1}
\end{eqnarray}

Notice that the $A_{d-1}$ Cartan matrix appears in $W_{d}$.
If $u_{1},u_{2},\cdots,u_{d-1}\neq 0$, we can set all the $u_{i}$'s to 1 by using the $(d+1)$ torus actions. The remaining  two torus actions 
that leave them invariant are nothing but the ones given in (\ref{two22}). Therefore, the subspace given by the condition 
$u_{1},u_{2},\cdots,u_{d-1}\neq 0$ corresponds to $Mp_{0,2}(N,d)$.
If $u_{d_1}=0, u_{j}\neq 0\;\;(j\neq d_{1})$, we have to delete the $u_{d_1}$ column of matrix $W_{d}$. 
This operation turns the $A_{d-1}$ Cartan matrix into the $A_{d_{1}-1}\times A_{d-d_{1}-1}$ Cartan matrix and results in 
chains of two polynomial maps:
\begin{equation}
(\sum_{j=0}^{d_{1}}{\bf a}_{j}s_{1}^{j}t_{1}^{d_{1}-j})\cup(\sum_{j=0}^{d-d_{1}}{\bf a}_{j+d_{1}}s_{2}^{j}t_{2}^{d-d_{1}-j}),\;\;
({\bf a}_{0},{\bf a}_{d_{1}},{\bf a}_{d}\neq {\bf 0}).
\end{equation}
Therefore, the corresponding boundary locus is given by $\displaystyle{Mp_{0,2}(N,d_1)\mathop{\times}_{CP^{N-1}}Mp_{0,2}(N,d-d_1)}$,
where $\displaystyle{\mathop{\times}_{CP^{N-1}}}$ is the fiber product with respect to the following projection maps:
\begin{eqnarray}
&&\pi_{1}:Mp_{0,2}(N,d_1)\rightarrow CP^{N-1},\;\;\pi_{1}({\bf a}_{0},\cdots,{\bf a}_{d_1})=[{\bf a}_{d_1}]\no\\ 
&&\pi_{2}:Mp_{0,2}(N,d-d_1)\rightarrow CP^{N-1},\;\;\pi_{2}({\bf a}_{d_1},\cdots,{\bf a}_{d})=[{\bf a}_{d_1}] 
\end{eqnarray}
In general, the subspace given by the condition 
\begin{equation}
u_{d_{i}}=0,\;(1\leq d_{1}<d_{2}<\cdots<d_{l(\sigma_{d})-1}\leq d-1), u_{j}\neq 0,\;\; (j\notin \{d_{1},d_{2},\cdots
,d_{l(\sigma_{d})-1}\}),
\end{equation} 
corresponds to 
chains of polynomial maps labeled by ordered partition $\sigma_{d}=(d_{1}-d_{0},d_{2}-d_{1},d_{3}-d_{2},\cdots,d_{l(\sigma_{d})}-
d_{l(\sigma_{d})-1})$:
\begin{equation}
\cup_{j=1}^{l(\sigma_{d})}\bigl(\sum_{m_{j}=0}^{d_{j}-d_{j-1}} {\bf a}_{d_{j-1}+m_{j}}(s_{j})^{m_{j}}(t_{j})^{d_{j}-d_{j-1}-m_{j}}\bigr),
\;\;\bigl({\bf a}_{d_{j}}\neq {\bf 0},\;\;j=0,1,\cdots,l(\sigma_{d})\bigr),
\label{chain0}
\end{equation}
where we set $d_{0}=0,d_{l(\sigma_{d})}=d$. In this case, the corresponding boundary locus is,
\begin{equation}
Mp_{0,2}(N,d_1-d_0)\mathop{\times}_{CP^{N-1}}Mp_{0,2}(N,d_2-d_1)\mathop{\times}_{CP^{N-1}}\cdots\mathop{\times}_{CP^{N-1}}
Mp_{0,2}(N,d_{l(\sigma_d)}-d_{l(\sigma_d)-1}).
\end{equation} 
Since the lowest dimensional boundary:
\begin{equation}
Mp_{0,2}(N,1)\mathop{\times}_{CP^{N-1}}Mp_{0,2}(N,1)\mathop{\times}_{CP^{N-1}}\cdots\mathop{\times}_{CP^{N-1}}
Mp_{0,2}(N,1),
\end{equation} 
is identified with the compact space $(CP^{N-1})^{d+1}$, we can conclude that $\widetilde{Mp}_{0,2}(N,d)$ is compact.

Next, we discuss the structure of the cohomology ring $H^{*}( \widetilde{Mp}_{0,2}(N,d) )$. In (\ref{toric1}), we labeled row vectors of 
$W_{d}$ by $h_{i}\;(i=0,1,\cdots,d)$, which represents K\"ahler forms of $\widetilde{Mp}_{0,2}(N,d)$ associated with the torus 
action of $\mu_i$ in (\ref{action1}). By using standard results on toric varieties, we can see that these $h_{i}$'s are generators of $H^{*}( \widetilde{Mp}_{0,2}(N,d) )$
and that relations between the generators are given by the data of elements of $W_{d}$ as follows:    
\begin{eqnarray}
&&(h_{0})^{N}=0,\;(h_{d})^{N}=0,\no\\
&&(h_{i})^{N}(2h_{i}-h_{i-1}-h_{i+1})=0,\;\;(i=1,2,\cdots,d-1).
\label{rel}
\end{eqnarray}

\subsubsection{Construction of Two Point Intersection Numbers on $\widetilde{Mp}_{0,2}(N,d)$ }
In this section, we define the following intersection number on $ \widetilde{Mp}_{0,2}(N,d) $, which is an analogue of 
a two point Gromov-Witten invariant of 
the degree $k$ hypersurface in $CP^{N-1}$:
\begin{eqnarray}
w({\cal O}_{h^{a}}{\cal O}_{h^{b}})_{0,d}:=
\int_{ \widetilde{Mp}_{0,2}(N,d) }ev_{1}^{*}( h^{a})\wedge ev_{2}^{*}( h^{b})\wedge c_{top}({\cal E}_d^k).\no\\
\label{gw1} 
\end{eqnarray}  
In (\ref{gw1}), $h$ is the hyperplane class of $CP^{N-1}$, and $ev_{1}: \widetilde{Mp}_{0,2}(N,d)\rightarrow CP^{N-1} $ (resp. 
$ev_{2}: \widetilde{Mp}_{0,2}(N,d)\rightarrow CP^{N-1}$) is the evaluation map at the first (resp. second) marked point. 
These maps are easily constructed as follows:
\begin{eqnarray}
ev_1([({\bf a}_{0},\cdots,{\bf a}_{d},u_1,\cdots,u_{d-1})]):=[{\bf a}_{0}]\in CP^{N-1},\no\\
ev_2([({\bf a}_{0},\cdots,{\bf a}_{d},u_1,\cdots,u_{d-1})]):=[{\bf a}_{d}]\in CP^{N-1}.
\label{evc}
\end{eqnarray}
We also 
have to construct a rank $(kd+1)$ orbi-bundle ${\cal E}_d^k$ on $\widetilde{Mp}_{0,2}(N,d)$
that corresponds to $R^{0}\pi_{*}ev_{3}^{*}({\cal O}_{CP^{N-1}}(k))$ on 
the moduli space of stable maps $\overline{M}_{0,2}(CP^{N-1},d)$. 
In this step, we need to consider the problem of point instantons that were introduced in the previous section.  In the case of 
$\widetilde{Mp}_{0,2}(N,d)$, we include point instantons to compactify the moduli space. On the other hand, these are prohibited in the case of 
$\overline{M}_{0,2}(CP^{N-1},d)$ because they are not actual maps.
This difference can be considered as the origin of the (generalized) mirror transformation. 
Therefore, our problem here is how to define an orbi-bundle corresponding to $R^{0}\pi_{*}ev_{3}^{*}({\cal O}_{CP^{N-1}}(k))$
for point instantons. Our approach to this task is quite naive. Let $s_0$ be a global holomorphic section ${\cal O}_{CP^{N-1}}(k)$.
It is well-known that $s_0$ is identified with a homogeneous polynomial of degree $k$ in homogeneous coordinates $X_1,X_2,\cdots,X_N$ 
of $CP^{N-1}$. Therefore, we can take 
\begin{eqnarray}
s_{0}=(X_1)^k+(X_2)^k+\cdots+(X_N)^k,
\label{section}
\end{eqnarray}
for example. Let us regard $\sum_{j=0}^{d}{\bf a}_js^{j}t^{d-j}\;({\bf a}_0,{\bf a}_d\neq {\bf 0})$ as a map $\varphi$ from
${\bf C}^{2}$ to ${\bf C}^{N}$. Of course, $[({\bf a}_0,{\bf a}_1,\cdots,{\bf a}_d)]$ represents a point in $Mp_{0,2}(N,d)$.
Then we can consider,
\begin{eqnarray}
\varphi^{*}s_{0}=\sum_{j=0}^{kd}\varphi_{j}^{k}({\bf a}_0,\cdots,{\bf a}_d)s^j t^{kd-j},
\label{pbse}
\end{eqnarray} 
where $\varphi_{j}^{k}({\bf a}_0,\cdots,{\bf a}_d)$ is a homogeneous polynomial of degree $k$ 
in  $a_j^i$ $({\bf a}_j=(a_j^1,a_j^2,\cdots,a_j^N))$. If we set
\begin{eqnarray}
\tilde{s}_{0}({\bf a}_0,\cdots,{\bf a}_d):=(\varphi_{0}^{k}({\bf a}_0,\cdots,{\bf a}_d),\varphi_{1}^{k}({\bf a}_0,\cdots,{\bf a}_d),\cdots
,\varphi_{kd}^{k}({\bf a}_0,\cdots,{\bf a}_d)), 
\label{sec2}
\end{eqnarray}
we can easily see that the image of the corresponding polynomial map lies inside the hypersurface defined by (\ref{section}) 
if and only if $\tilde{s}_{0}({\bf a}_0,\cdots,{\bf a}_d)={\bf 0}$. Moreover, we can derive the following relations:
\begin{eqnarray}
\tilde{s}_{0}(\mu{\bf a}_0,\cdots,\mu{\bf a}_d)&=&(\mu^{k}\varphi_{0}^{k}({\bf a}_0,\cdots,{\bf a}_d),\mu^{k}\varphi_{1}^{k}({\bf a}_0,
\cdots,{\bf a}_d),\cdots,\mu^{k}\varphi_{kd}^{k}({\bf a}_0,\cdots,{\bf a}_d)),\no\\
\tilde{s}_{0}({\bf a}_0,\nu{\bf a}_1,\nu^2{\bf a}_2\cdots,\nu^{d-1}{\bf a}_d)&=&
(\varphi_{0}^{k}({\bf a}_0,\cdots,{\bf a}_d),\nu\varphi_{1}^{k}({\bf a}_0,\cdots,{\bf a}_d),
\nu^2\varphi_{2}^{k}({\bf a}_0,\cdots,{\bf a}_d),\cdots
,\nu^{kd}\varphi_{kd}^{k}({\bf a}_0,\cdots,{\bf a}_d)).\no\\
\label{homrel}
\end{eqnarray}
These relations tells us that $\tilde{s}_0$ defines a section of a rank $kd+1$ orbi-bundle on $Mp_{0,2}(N,d)$, 
because we can compute transition functions of the bundle by using (\ref{homrel}). Let us discuss this argument more explicitly. 
Since $Mp_{0,2}(N,d)=(CP^{N-1}\times{\bf C}^{N(d-1)}\times CP^{N-1})/{\bf Z}_d$, we can take the following local coordinate system 
$U_{ij}$. 
\begin{eqnarray}
&&\phi_{ij}:U_{ij}\subset{\bf C}^{N(d+1)-2}\rightarrow Mp_{0,2}(N,d),\no\\
&&\phi_{ij}(x_1,x_2,\cdots,x_{N-1},{\bf y}_{1},{\bf y}_2,\cdots{\bf y}_{d-1},z_1,z_2.\cdots,z_{N-1})=\no\\
&&[(x_1,\cdots,x_{i-1},1,x_i,\cdots,x_{N-1},{\bf y}_{1},\cdots,{\bf y}_{d-1},z_1,\cdots,z_{j-1},1,z_{j},\cdots,z_{N-1})],
\label{local}
\end{eqnarray}
where ${\bf y}_i\in{\bf C}^{N}$. Let $(\tilde{x}_{*},\tilde{\bf y}_{*},\tilde{z}_{*})\in U_{kl}$. We assume that $i<k$ and $j<l$ for 
simplicity. The coordinate transformation between $U_{ij}$ and $U_{kl}$ is given by,
\begin{eqnarray}
&&x_{m}=\frac{\tilde{x}_m}{\tilde{x}_i} \;\;(m\leq i-1),\;\;x_{m}=\frac{\tilde{x}_{m+1}}{\tilde{x}_i}\;\;(i\leq m\leq k-2),\;\; 
x_{k-1}=\frac{1}{\tilde{x}_i},
\;\;x_{m}=\frac{\tilde{x}_{m}}{\tilde{x}_i}\;\;(k\leq m\leq N-1),\no\\
&&z_{m}=\frac{\tilde{z}_m}{\tilde{z}_j} \;\;(m\leq j-1),\;\;z_{m}=\frac{\tilde{z}_{m+1}}{\tilde{z}_j}\;\;(j\leq m\leq l-2),\;\; 
z_{l-1}=\frac{1}{\tilde{z}_j},
\;\;z_{m}=\frac{\tilde{z}_{m}}{\tilde{z}_j}\;\;(l\leq m\leq N-1),\no\\
&&{\bf y}_{m}=\frac{1}{(\tilde{x}_i)^{\frac{d-m}{d}}(\tilde{z}_j)^{\frac{m}{d}}}\tilde{\bf y}_m.
\end{eqnarray}
If we represent the section $s_{0}$ on $U_{ij}$ by $\tilde{s}_{0}(\phi_{ij}(x_{*},{\bf y}_{*},z_{*}))=(\varphi_{0}(x_{*},{\bf y}_{*},z_{*}),\varphi_{1}(x_{*},{\bf y}_{*},z_{*}),\cdots,\varphi_{kd}(x_{*},{\bf y}_{*},z_{*}))$, we obtain the following 
relation:
\begin{eqnarray}
\frac{1}{(\tilde{x}_i)^{\frac{kd-m}{d}}(\tilde{z}_j)^{\frac{m}{d}}}\varphi_{m}(\tilde{x}_{*},\tilde{\bf y}_{*},\tilde{z}_{*})=
\varphi_{m}(x_{*},{\bf y}_{*},{z}_{*}),\;\;(m=0,1,2,\cdots,kd).
\end{eqnarray} 
Therefore, we can regard $\tilde{s}_{0}$ as a section of the rank $kd+1$ bundle whose transition function is given by,
\begin{eqnarray}
(\tilde{x}_i)^{\frac{kd-m}{d}}(\tilde{z}_j)^{\frac{m}{d}}\tilde{e}_m=e_m, \;\;(m=0,1,2,\cdots,kd), 
\label{trans}
\end{eqnarray}
where $e_m$ (resp. $\tilde{e}_m$) is the base of trivialization on $U_{ij}$ (resp. $U_{kl}$).
We denote this orbi-bundle on $Mp_{0,2}(N,d)$ 
by ${\cal E}_d^k$. From (\ref{trans}), we can see that ${\cal E}_d^k\simeq\mathop{\oplus}_{m=0}^{kd}\bigl({\cal O}_{CP^{N-1}}(\frac{kd-m}{d})
\otimes{\cal O}_{CP^{N-1}}(\frac{m}{d})\bigr)$ as an orbi-bundle on $Mp_{0,2}(N,d)$.
Note that we can define ${\cal E}_d^k$ on whole $Mp_{0,2}(N,d)$ whether $\sum_{j=0}^{d}{\bf a}_js^{j}t^{d-j}$
is a point instanton or not.
Next, we extend ${\cal E}_d^k$ to $\widetilde{Mp}_{0,2}(N,d)$. Let us consider the locus in $\widetilde{Mp}_{0,2}(N,d)$
where 
\begin{eqnarray}
u_{d_j}&=&0,\;\;(1\leq d_1<d_2<\cdots<d_{l-1}\leq d-1),\no\\
u_{j}&\neq&0,\;\;(j\notin\{d_1,d_2,\cdots,d_{l-1}\}).
\label{loc1}
\end{eqnarray}
We denote this locus by $U_{(d_0,d_1,\cdots,d_l)}$ $(d_0:=0,\;d_l:=d)$. As was discussed in the previous  section, $U_{(d_0,d_1,\cdots,d_l)}$ is identified with, 
\begin{equation}
Mp_{0,2}(N,d_1-d_0)\mathop{\times}_{CP^{N-1}}Mp_{0,2}(N,d_2-d_1)\mathop{\times}_{CP^{N-1}}\cdots\mathop{\times}_{CP^{N-1}}
Mp_{0,2}(N,d_{l}-d_{l-1}),
\end{equation}
and its point is represented by a chain of polynomial maps:
\begin{eqnarray}
\mathop{\cup}_{j=1}^{l}(\sum_{h=0}^{d_j-d_{j-1}}{\bf a}_{d_{j-1}+h}(s_j)^h(t_j)^{d_j-d_{j-1}-h}).
\end{eqnarray}
For each $Mp_{0,2}(N,d_j-d_{j-1})$, we have $k(d_{j}-d_{j-1})+1$ 
dimensional orbi-bundle ${\cal E}_{d_j-d_{j-1}}^k$. We then introduce
a map $p_j:U_{(d_0,d_1,\cdots,d_l)}\rightarrow CP^{N-1}$ $(j=1,2,\cdots,l-1)$ defined by, 
\begin{eqnarray}
p_j(\mathop{\cup}_{j=1}^{l}(\sum_{h=0}^{d_j-d_{j-1}}{\bf a}_{d_{j-1}+h}(s_j)^h(t_j)^{d_j-d_{j-1}-h}))=[{\bf a}_{d_j}]\in CP^{N-1}.
\label{eproj}
\end{eqnarray}
With this setup, we define ${\cal E}_d^{k}|_{U_{(d_0,d_1,\cdots,d_l)}}$ by the following exact sequence:
\begin{eqnarray}
0\rightarrow {\cal E}_d^{k}|_{U_{(d_0,d_1,\cdots,d_l)}}\rightarrow 
\mathop{\oplus}_{j=1}^{l}{\cal E}_{d_j-d_{j-1}}^{k}\rightarrow \mathop{\oplus}_{j=1}^{l-1}p_{j}^{*}{\cal O}_{CP^{N-1}}(k)\rightarrow 0. 
\label{exseq1} 
\end{eqnarray}  
$ {\cal E}_d^{k}|_{U_{(d_0,d_1,\cdots,d_l)}}$ also has rank $kd+1$. In this way, we extend ${\cal E}_{d}^{k}$ to whole 
$\widetilde{Mp}_{0,2}(N,d)$.

We can also construct a rank $kd-1$ orbi-bundle ${\cal E}_{d}^{-k}$ on $Mp_{0,2}(N,d)$ that is isomorphic to $\mathop{\oplus}_{m=1}^{kd-1}\bigl({\cal O}_{CP^{N-1}}(\frac{m-kd}{d})
\otimes{\cal O}_{CP^{N-1}}(\frac{-m}{d})\bigr)$ as an orbi-bundle on $Mp_{0,2}(N,d)$.
We can also extend ${\cal E}_{d}^{-k}$ to the whole $\widetilde{Mp}_{0,2}(N,d)$
by using the exact sequence:
\begin{eqnarray}
0\rightarrow \mathop{\oplus}_{j=1}^{l-1}p_{j}^{*}{\cal O}_{CP^{N-1}}(-k)\rightarrow {\cal E}_d^{-k}|_{U_{(d_0,d_1,\cdots,d_l)}}\rightarrow 
\mathop{\oplus}_{j=1}^{l}{\cal E}_{d_j-d_{j-1}}^{-k}\rightarrow 0. 
\label{exseq2} 
\end{eqnarray}  
This bundle corresponds to $R^{1}\pi_{*}ev_{3}^{*}{\cal O}_{CP^{N-1}}(-k)$ on $\overline{M}_{0,2}(CP^{N-1},d)$ by Kodaira-Serre duality.

\subsubsection{Localization Computation of $w({\cal O}_{h^{a}}{\cal O}_{h^{b}})_{0,d}$}
In this section, we compute the intersection number $w({\cal O}_{h^{a}}{\cal O}_{h^{b}})_{0,d}$ by using the localization theorem.
For this purpose, we introduce the following ${\bf C}^{\times}$ action on $\widetilde{Mp}_{0,2}(N,d)$.
\begin{equation}
[(e^{\lambda_0 t}{\bf a}_0,e^{\lambda_1 t}{\bf a}_1,\cdots,e^{\lambda_{d-1} t}{\bf a}_{d-1},e^{\lambda_d t}{\bf a}_d,u_1,u_2,\cdots,u_{d-1})].
\label{tflow}
\end{equation}
In (\ref{tflow}), $\lambda_i\;\;(i=0,1,\cdots,d)$ is the equivariant parameter for the flow. In \cite{vs}, we took non-equivariant limit 
$\lambda_i \rightarrow 0$ from the start, but in this section, we perform the computation under non-zero equivariant parameters. The fixed point
sets of $\widetilde{Mp}_{0,2}(N,d)$ consist of connected components, each of which come from $U_{(d_0,d_1,\cdots,d_l)}$ defined in the previous section. We denote the connected component that comes from $U_{(d_0,d_1,\cdots,d_l)}$ by $F_{(d_0,d_1,\cdots,d_l)}$. Explicitly, a point in 
$F_{(d_0,d_1,\cdots,d_l)}$ is represented by the following chain of polynomial maps.
\begin{equation} 
\mathop{\cup}_{j=1}^{l}({\bf a}_{d_{j-1}}(s_j)^{d_j-d_{j-1}}+{\bf a}_{d_{j}}(t_j)^{d_j-d_{j-1}}).
\label{fixchain}
\end{equation}
Note here that $({\bf a}_{d_{j-1}}(s_j)^{d_j-d_{j-1}}+{\bf a}_{d_{j}}(t_j)^{d_j-d_{j-1}})$ is the ${\bf Z}_{d_{j}-d_{j-1}}$ singularity 
in $Mp_{0,2}(N,d_j-d_{j-1})$. We can easily see from (\ref{fixchain}) that $F_{(d_0,d_1,\cdots,d_l)}$ is set-theoretically isomorphic to 
$\prod_{j=0}^{l}(CP^{N-1})_{d_j}$ where $(CP^{N-1})_{d_j}$ is the $CP^{N-1}$ whose point is given by $[{\bf a}_{d_j}]$. 

Let us consider the contribution to $w({\cal O}_{h^{a}}{\cal O}_{h^{b}})_{0,d}$ from $F_{(d_0,d_1,\cdots,d_l)}$. We start from the case of 
$F_{(0,d)}\subset U_{(0,d)}=Mp_{0,2}(N,d)$. First, we have to determine the normal bundle of $F_{(0,d)}$ in $Mp_{0,2}(N,d)$. We already know 
from the previous discussion that,
\begin{equation}
Mp_{0,2}(N,d)=\{([{\bf a}_{0}],{\bf y}_{1},\cdots,{\bf y}_{d-1},[{\bf a}_{d}])\;|\;[{\bf a}_0],[{\bf a}_d]\in CP^{N-1},\;{\bf y}_i\in C^{N}\}
/{\bf Z}_d.
\label{Mpiso}
\end{equation}
Therefore, the normal bundle is given by $\displaystyle{\mathop{\oplus}_{i=1}^{d-1}\mathop{\oplus}_{j=1}^{N}\frac{\partial}{\partial y_{i}^{j}}}$. From the discussion of the previous section, we can see that $\frac{\partial}{\partial y_{i}^{j}}$ is isomorphic to 
${\cal O}_{(CP^{N-1})_{0}}(\frac{d-i}{d})\otimes{\cal O}_{(CP^{N-1})_{d}}(\frac{i}{d})$ as an orbi-bundle on $F_{(0,d)}$ and its first Chern class is 
given by,
\begin{equation}
\frac{d-i}{d}h_{0}+\frac{i}{d}h_d,
\end{equation}
where $h_{d_i}$ is the hyperplane class of $(CP^{N-1})_{d_i}$. On the other hand, the flow in (\ref{tflow}) acts on $y_{i}^{j}$ as 
$y_{i}^{j}\rightarrow e^{\bigl(\lambda_i-(\frac{d-i}{d}\lambda_{0}+\frac{i}{d}\lambda_d)\bigr)t}y_i^j$, and the character of the flow on $\frac{\partial}{\partial y_{i}^{j}}$ is given by, 
\begin{eqnarray}
\frac{d-i}{d}\lambda_{0}+\frac{i}{d}\lambda_d-\lambda_i.
\end{eqnarray}  
Next, we consider equivariant top Chern class of ${\cal E}_d^k$ on $F_{(0,d)}$. Since ${\cal E}_d^k$ is identified with $\mathop{\oplus}_{m=0}^{kd}\bigl({\cal O}_{(CP^{N-1})_0}(\frac{kd-m}{d})
\otimes{\cal O}_{(CP^{N-1})_d}(\frac{m}{d})\bigr)$ as an orbi-bundle on $Mp_{0,2}(N,d)$, its equivariant top Chern class on  $F_{(0,d)}$
is given by,
\begin{eqnarray}
\prod_{m=0}^{kd}\bigl(\frac{(kd-m)(h_0+\lambda_0)+m(h_d+\lambda_d)}{d}\bigr).
\end{eqnarray}
From the definition of the evaluation map for $Mp_{0,2}(N,d)$ in (\ref{evc}), we can easily see that equivariant representation of $ev_{1}^{*}(h^a)$ (resp. $ev_{2}^{*}(h^b)$) on $F_{(0,d)}$ is given by $(h_{0}+\lambda_0)^a$ (resp. $(h_d+\lambda_d)^b$).　
Finally, we have to remember that $F_{(0,d)}$ is also the singular locus on which ${\bf Z}_d$ acts. Therefore, we have to divide 
the results of integration on $F_{(0,d)}$ by $d$. Putting these results altogether, the contribution from $F_{(0,d)}$ becomes,
\begin{equation}
\frac{1}{d}\int_{(CP^{N-1})_0}\int_{(CP^{N-1})_d}(h_{0}+\lambda_0)^a \frac{\prod_{m=0}^{kd}\bigl(\frac{(kd-m)(h_0+\lambda_0)+m(h_d+\lambda_d)}{d}\bigr)}
{\prod_{i=1}^{d-1}\bigl(\frac{(d-i)(h_0+\lambda_0)+i(h_d+\lambda_d)}{d}-\lambda_i\bigr)^N}(h_d+\lambda_d)^b.
\label{c0d}
\end{equation} 
We then consider the contribution from $F_{(d_0,d_1.\cdots,d_l)}\;\;(l\geq 2)$. As for the normal bundle, we have additional factors coming from 
"smoothing the nodal singularities" of the image of the chain of polynomial maps, that are  given by $[{\bf a}_{d_j}]\;\;(j=1,2,\cdots,l-1)$. This factor 
is identified with the orbi-bundle $\frac{d}{d(\frac{s_j}{t_j})}\otimes\frac{d}{d(\frac{t_{j+1}}{s_{j+1}})}$ and its equivariant first Chern 
class is given by, 
\begin{eqnarray}
\frac{h_{d_j}+\lambda_{d_j}-h_{d_{j-1}}-\lambda_{d_{j-1}}}{d_j-d_{j-1}}+\frac{h_{d_j}+\lambda_{d_j}-h_{d_{j+1}}-\lambda_{d_{j+1}}}{d_{j+1}-d_{j}}.   
\end{eqnarray}
Equivariant top Chern class of ${\cal E}_d^k$ on $F_{(d_0,d_1,\cdots,d_l)}$ can be read off from the exact sequence in (\ref{exseq1}) 
as follows.
\begin{equation}
\frac{\prod_{j=1}^{l}\prod_{m=0}^{k(d_j-d_{j-1})}\bigl(\frac{(k(d_j-d_{j-1})-m)(h_{d_{j-1}}+\lambda_{d_{j-1}})+m(h_{d_j}+\lambda_{d_j})}{d_j-d_{j-1}}\bigr)}{\prod_{j=1}^{l-1}k(h_{d_j}+\lambda_{d_j})}.
\label{gentop}
\end{equation}
Combining these addtional factors with the consideration in the case of $F_{(0,d)}$, we can write down the contribution that comes from 
$F_{(d_0,d_1.\cdots,d_l)}$.
\begin{eqnarray}
&&\frac{1}{\prod_{j=1}^{l}(d_j-d_{j-1})}\int_{(CP^{N-1})_{d_0}}\cdots\int_{(CP^{N-1})_{d_l}}
(h_{d_0}+\lambda_{d_0})^a \times\no\\
&&\frac{\prod_{j=1}^{l}\prod_{m=0}^{k(d_j-d_{j-1})}\bigl(\frac{(k(d_j-d_{j-1})-m)(h_{d_{j-1}}+\lambda_{d_{j-1}})+m(h_{d_j}+\lambda_{d_j})}{d_j-d_{j-1}}\bigr)}
{\prod_{j=1}^{l}\prod_{i=1}^{d_j-d_{j-1}-1}\bigl(\frac{(d_j-d_{j-1}-i)(h_{d_{j-1}}+\lambda_{d_{j-1}})+i(h_{d_j}+\lambda_{d_j})}{d_j-d_{j-1}}-\lambda_{d_{j-1}+i}\bigr)^N}\times\no\\
&&\frac{1}{\prod_{j=1}^{l-1}\bigl(\frac{h_{d_j}+\lambda_{d_j}-h_{d_{j-1}}-\lambda_{d_{j-1}}}{d_j-d_{j-1}}+\frac{h_{d_j}+\lambda_{d_j}-h_{d_{j+1}}-\lambda_{d_{j+1}}}{d_{j+1}-d_{j}}\bigr)\bigl(k(h_{d_j}+\lambda_{d_j})\bigr)}(h_{d_l}+\lambda_{d_l})^b.
\label{cgend}
\end{eqnarray}
Since $\int_{CP^{N-1}}h^a=\frac{1}{2\pi\sqrt{-1}}\oint_{C_{(0)}}\frac{dz}{z^{N}}z^a$, we obtain the following closed formula for $w({\cal O}_{h^{a}}{\cal O}_{h^{b}})_{0,d}$.
\begin{eqnarray}
w({\cal O}_{h^{a}}{\cal O}_{h^{b}})_{0,d}&=&\sum_{0=d_0<d_1<\cdots<d_{l-1}<d_l=d}\frac{1}{\prod_{j=1}^{l}(d_j-d_{j-1})}\frac{1}{(2\pi\sqrt{-1})^{l+1}}\oint_{C_{(0)}}\frac{dz_{d_0}}{(z_{d_0})^{N}}\cdots\oint_{C_{(0)}}\frac{dz_{d_l}}{(z_{d_l})^{N}}\times\no\\
&&(z_{d_0}+\lambda_{d_0})^a\frac{\prod_{j=1}^{l}\prod_{m=0}^{k(d_j-d_{j-1})}\bigl(\frac{(k(d_j-d_{j-1})-m)(z_{d_{j-1}}+\lambda_{d_{j-1}})+m(z_{d_j}+\lambda_{d_j})}{d_j-d_{j-1}}\bigr)}
{\prod_{j=1}^{l}\prod_{i=1}^{d_j-d_{j-1}-1}\bigl(\frac{(d_j-d_{j-1}-i)(z_{d_{j-1}}+\lambda_{d_{j-1}})+i(z_{d_j}+\lambda_{d_j})}{d_j-d_{j-1}}-\lambda_{d_{j-1}+i}\bigr)^N}\times\no\\
&&\frac{1}{\prod_{j=1}^{l-1}\bigl(\frac{z_{d_j}+\lambda_{d_j}-z_{d_{j-1}}-\lambda_{d_{j-1}}}{d_j-d_{j-1}}+\frac{z_{d_j}+\lambda_{d_j}-z_{d_{j+1}}-\lambda_{d_{j+1}}}{d_{j+1}-d_{j}}\bigr)\bigl(k(z_{d_j}+\lambda_{d_j})\bigr)}(z_{d_l}+\lambda_{d_l})^b.
\label{wlint}
\end{eqnarray}
In the above formula, we can integrate the variable $z_{d_j}$ in arbitrary order. The formula (\ref{wlint}) has the form of residue integral
and we can take non-equivariant limit $\lambda_j\rightarrow 0$.  This operation makes the formula simpler. For simplicity, we introduce the following notations. We define the following two polynomials in $z$ and $w$:
\begin{eqnarray}
e(k,d;z,w)&:=&\prod_{j=0}^{kd}\bigl(\frac{jz+(kd-j)w}{d}\bigr)\no\\
t(N,d;z,w)&:=&\prod_{j=1}^{d-1}\bigl(\frac{jz+(d-j)w}{d}\bigr)^{N}.
\end{eqnarray}
We also introduce the ordered partition of a positive integer $d$: 
\begin{defi}
Let $OP_{d}$ be the set of ordered partitions of a positive integer $d$:
\begin{equation}
OP_{d}=\{\sigma_{d}=(d_{1},d_{2},\cdots,d_{l(\sigma_{d})})\;\;|\;\;
\sum_{j=1}^{l(\sigma_{d})}d_{j}=d\;\;,\;\;d_{j}\in{\bf N}\}.
\label{oparti} 
\end{equation}
In (\ref{oparti}), we denoted 
the length of the ordered partition $\sigma_{d}$ by $l(\sigma_{d})$.
\end{defi}
The increasing sequence of integer $(d_0,d_1,\cdots,d_l)\;\;(0=d_0<d_1<\cdots<d_{l-1}<d_l=d)$ used in (\ref{wlint}) can be 
replaced by the ordered partition $\sigma_d=(\tilde{d}_1,\tilde{d}_2,\cdots,\tilde{d}_l)\in OP_{d}$ if we use the following 
correspondence:
\begin{equation}
\tilde{d}_j=d_j-d_{j-1},\;\;(j=1,2,\cdots,l).
\end{equation}
With this setup, we can simplify the formula for $w({\cal O}_{h^{a}}{\cal O}_{h^{b}})_{0,d}$ after taking the non-equivariant limit, 
by relabeling the subscript of $z_{*}'s$ as follows.
\begin{eqnarray}
w({\cal O}_{h^{a}}{\cal O}_{h^{b}})_{0,d}&=&\sum_{\sigma_{d}\in OP_{d}}\frac{1}{(2\pi\sqrt{-1})^{l(\sigma_{d})+1}\prod_{j=0}^{l(\sigma_{d})}d_{j}}\oint_{C_{0}} \frac{dz_{0}}{(z_{0})^N}\cdots
\oint_{C_{0}} \frac{d z_{l(\sigma_{d})} }{(z_{l(\sigma_{d})})^N}(z_{0})^{a}\times\no\\
&&\times\prod_{j=1}^{l(\sigma_{d})-1}\frac{1}
{\biggl( \frac{z_{j}-z_{j-1}}{d_{j}}+\frac{z_{j}-z_{j+1}}{d_{j+1}}\biggr)kz_{j}}\prod_{j=1}^{l(\sigma_{d})}
\frac{e(k,d_j;z_{j-1},z_{j})}{t(N,d_j;z_{j-1},z_{j})}(z_{l(\sigma_{d})})^{b}.
\label{wnonequi}
\end{eqnarray}
\begin{Rem}
After taking non-equivariant limit, we have to take care of the order of integration of $z_j$'s. In (\ref{wnonequi}), we have to 
integrate $z_{j}'s$ in all the summands of the formula in descending (or ascending) order of the subscript $j$.   
\end{Rem}
\subsubsection{Numerical Results}
{\bf $N=7, k=5$ case }\\
In this case, $M_{7}^{5}$ is a Fano hypersurface and $N-k=2$. From (\ref{wnonequi}), we obtain the following $w({\cal O}_{h^{a}}{\cal O}_{h^{b}})_{0,d}$'s.
\begin{eqnarray}
&&w({\cal O}_{h^{1}}{\cal O}_{h^{5}})_{0,1}=600,\;\; w({\cal O}_{h^{2}}{\cal O}_{h^{4}})_{0,1}=3850,\;\; 
w({\cal O}_{h^{3}}{\cal O}_{h^{3}})_{0,1}=6725,\no\\
&& w({\cal O}_{h^{3}}{\cal O}_{h^{5}})_{0,2}=528000,\;\;
w({\cal O}_{h^{4}}{\cal O}_{h^{4}})_{0,2}=1731250,\;\;w({\cal O}_{h^{5}}{\cal O}_{h^{5}})_{0,3}=52200000.\;\;
\label{w75}
\end{eqnarray} 
On the other hand, we can evaluate the corresponding Gromov-Witten invariants by localization computation or mirror computation. 
The results are given as follows.
\begin{eqnarray}
&& \langle{\cal O}_{h^{1}}{\cal O}_{h^{5}}\rangle_{0,1}=600,\;\;\langle{\cal O}_{h^{2}}{\cal O}_{h^{4}}\rangle_{0,1}=3850 ,\;\;
 \langle{\cal O}_{h^{3}}{\cal O}_{h^{3}}\rangle_{0,1}=6725,\;\;\no\\
&& \langle{\cal O}_{h^{3}}{\cal O}_{h^{5}}\rangle_{0,2}=528000 ,\;\;\langle{\cal O}_{h^{4}}{\cal O}_{h^{4}}\rangle_{0,2}=1731250 ,\;\;
 \langle{\cal O}_{h^{5}}{\cal O}_{h^{5}}\rangle_{0,3}=52200000.
 \label{gw75}
\end{eqnarray} 
Therefore, we have $w({\cal O}_{h^{a}}{\cal O}_{h^{b}})_{0,d}=\langle{\cal O}_{h^{a}}{\cal O}_{h^{b}}\rangle_{0,d}$ in this case.
\\
\\
{\bf $N=5, k=5$ case }\\
Since $M_{5}^{5}$ is the celebrated quintic 3-fold, we have the following data of 2-point Gromov-Witten invariants.
\begin{eqnarray}
&& \langle{\cal O}_{h^{0}}{\cal O}_{h^{2}}\rangle_{0,1}=0,\;\;\langle{\cal O}_{h^{0}}{\cal O}_{h^{2}}\rangle_{0,2}=0,\;\;
 \langle{\cal O}_{h^{0}}{\cal O}_{h^{2}}\rangle_{0,3}=0,\cdots,\no\\
&&  \langle{\cal O}_{h^{1}}{\cal O}_{h^{1}}\rangle_{0,1}=2875,\;\;\langle{\cal O}_{h^{1}}{\cal O}_{h^{1}}\rangle_{0,2}=\frac{4876875}{2},\;\;
 \langle{\cal O}_{h^{1}}{\cal O}_{h^{1}}\rangle_{0,3}=\frac{8564575000}{3},\cdots.
 \label{gw55}
\end{eqnarray} 
The fact that $\langle{\cal O}_{h^{0}}{\cal O}_{h^{2}}\rangle_{0,d}=0$ follows from the puncture axiom of Gromov-Witten invariants.
On the other hand, the corresponding $w({\cal O}_{h^{a}}{\cal O}_{h^{b}})_{0,d}$'s are given as follows.
\begin{eqnarray}
&& w({\cal O}_{h^{0}}{\cal O}_{h^{2}})_{0,1}=3850,\;\;w({\cal O}_{h^{0}}{\cal O}_{h^{2}})_{0,2}=3589125,\;\;
 w({\cal O}_{h^{0}}{\cal O}_{h^{2}})_{0,3}=\frac{16126540000}{3},\cdots,\no\\
&&  w({\cal O}_{h^{1}}{\cal O}_{h^{1}})_{0,1}=6725,\;\;w({\cal O}_{h^{1}}{\cal O}_{h^{1}})_{0,2}=\frac{16482625}{2},\;\;
 w({\cal O}_{h^{1}}{\cal O}_{h^{1}})_{0,3}=\frac{44704818125}{3},\cdots.
 \label{w55}
\end{eqnarray} 
In this case, $w({\cal O}_{h^{a}}{\cal O}_{h^{b}})_{0,d}$ and  $\langle{\cal O}_{h^{a}}{\cal O}_{h^{b}}\rangle_{0,d}$ differ from 
each other. Let us consider here the generating function:
\begin{equation}
t(x):=x+\sum_{d=1}^{\infty}\frac{w({\cal O}_{h^{0}}{\cal O}_{h^{2}})_{0,d}}{5}e^{dx}=x+770e^x+717825e^{2x}+\frac{3225308000}{3}e^{3x}+\cdots.
\label{55mirror}
\end{equation}
This is nothing but the mirror map used in the mirror computation of the quintic 3-fold!
If we introduce another generating function:
\begin{equation}
F(x):=5x+\sum_{d=1}^{\infty}w({\cal O}_{h}{\cal O}_{h})_{0,d}e^{dx}=5x+6725e^x+\frac{16482625}{2}e^{2x}+\frac{44704818125}{3}e^{3x}+\cdots,
\label{55-2ptw}
\end{equation}
$F(x(t))$ gives us the generating function of $\langle{\cal O}_{h^{1}}{\cal O}_{h^{1}}\rangle_{0,d}$.
\begin{equation}
F(x(t))=5t+2875e^t+\frac{4876875}{2}e^{2t}+\frac{8564575000}{3}e^{3t}+\cdots.
\label{55-2pt}
\end{equation}
In section 3, we generalize these results to the case of some Calabi-Yau 3-folds with two K\"ahler forms.  
\\
\\
{\bf $N=8, k=9$ case }\\
In this case, $M_{8}^{9}$ is non-nef. The non-zero $w({\cal O}_{h^{a}}{\cal O}_{h^{b}})_{0,d}$'s up to $d=3$ are evaluated as follows.
\begin{eqnarray}
&&w({\cal O}_{h^{0}}{\cal O}_{h^{4}})_{0,1}=307250172 ,\;\;w({\cal O}_{h^{1}}{\cal O}_{h^{3}})_{0,1}=817713468 ,\;\;
w({\cal O}_{h^{2}}{\cal O}_{h^{2}})_{0,1}=1122806529 ,\;\;\no\\
&&w({\cal O}_{h^{0}}{\cal O}_{h^{3}})_{0,2}=75644409992388462,\;\;w({\cal O}_{h^{1}}{\cal O}_{h^{2}})_{0,2}=\frac{733562379269675757}{4} ,\;\;\no\\
&&w({\cal O}_{h^{0}}{\cal O}_{h^{2}})_{0,3}=34343397483304162555939158,\;\;w({\cal O}_{h^{1}}{\cal O}_{h^{1}})_{0,3}=56677396498174471672277559.
\label{w89}
\end{eqnarray}
On the other hand, the corresponding $\langle{\cal O}_{h^{a}}{\cal O}_{h^{b}}\rangle_{0,d}$'s are evaluated as follows.
\begin{eqnarray}
&&\langle{\cal O}_{h^{0}}{\cal O}_{h^{4}}\rangle_{0,1}=0 ,\;\;\langle{\cal O}_{h^{1}}{\cal O}_{h^{3}}\rangle_{0,1}=510463296 ,\;\;
\langle{\cal O}_{h^{2}}{\cal O}_{h^{2}}\rangle_{0,1}=815556357 ,\;\;\no\\
&&\langle{\cal O}_{h^{0}}{\cal O}_{h^{3}}\rangle_{0,2}=0,\;\;\langle{\cal O}_{h^{1}}{\cal O}_{h^{2}}\rangle_{0,2}=\frac{319615925538369285}{4},\;\;\no\\
&&\langle{\cal O}_{h^{0}}{\cal O}_{h^{2}}\rangle_{0,3}=0,\;\;\langle{\cal O}_{h^{1}}{\cal O}_{h^{1}}\rangle_{0,3}=12112667926597160835676659.
\label{gw89}
\end{eqnarray}
From the numerical data in \cite{gene0}, we can observe that $w({\cal O}_{h^{a}}{\cal O}_{h^{b}})_{0,d}$ is related to the virtual structure constants $\tilde{L}_n^{N,k,d}$ by the following equality:
\begin{equation}
k\tilde{L}_{n}^{N,k,d}=d\cdot w({\cal O}_{h^{6-n}}{\cal O}_{h^{n-1-d}})_{0,d}.
\end{equation}
Therefore, $w({\cal O}_{h^{a}}{\cal O}_{h^{b}})_{0,d}$'s are translated into $\langle{\cal O}_{h^{a}}{\cal O}_{h^{b}}\rangle_{0,d}$'s 
via the relations given in (\ref{th3}).\\
\\
The results in this section is the examples of Theorem 1, that will be proved in the next section.

\subsection{Proof of Theorem 1}
\subsubsection{Definition of the Virtual Structure Constants} 
In this subsection, we prove the conjecture proposed in \cite{vs} that represents the virtual structure constant $\tilde{L}_{n}^{N,k,d}$ 
for the degree $k$ hypersurface in $CP^{N-1}$ as a residue integral. We first write down the definition of $\tilde{L}_{n}^{N,k,d}$ 
given in our early papers \cite{fano,gm}. 
We introduce here a polynomial $Poly_{d}$ in 
$x,y,z_{1},z_{2},\cdots,z_{d-1}$ defined by the formula: 
\begin{eqnarray}
&&Poly_{d}(x,y,z_{1},z_{2},\cdots,z_{d-1})\no\\
&&:=\frac{d}{(2\pi\sqrt{-1})^{d-1}}
\oint_{D_{1}}{du_{1}}\cdots
\oint_{D_{d-1}}{du_{d-1}}\prod_{j=1}^{d-1}\biggl(
\frac{(u_{j})^{2}}{(2u_{j}-u_{j-1}-u_{j+1})(u_{j}-z_{j})}\biggr),
\label{trial2} 
\end{eqnarray}
where we denote $x$ (resp. $y$ ) by $u_{0}$ (resp. $u_{d}$) in the second line.
In (\ref{trial2}), $\frac{1}{2\pi\sqrt{-1}}\oint_{D_{j}}du_j$ represents, 
$$\frac{1}{2\pi\sqrt{-1}}\oint_{C_{(z_j,\frac{u_{j-1}+u_{j+1}}{2})}}du_j.$$ 
Let us consider the following "comb type" of a positive integer $d$ 
:
\begin{equation}
0=i_{0}<i_{1}<i_{2}<\cdots<i_{l-1}<i_{l}=d.
\end{equation}
The monomials that appear in $Poly_d$ are represented by,  
$$x^{m_{i_{0}}}z_{i_{1}}^{m_{i_{1}}}\cdots z_{i_{l-1}}^{m_{i_{l-1}}}
y^{m_{i_{l}}}, \;\;\;\;(\sum_{j=0}^{l}m_{i_{j}}=d-1).$$
We list some elements in ${\bf Z}^{l}$, which are 
determined for each comb type
as follows:
\begin{eqnarray}
\alpha&:=&(l-d,l-d,\cdots,l-d),\no\\
\beta&:=&(0,i_{1}-1,i_{2}-2,\cdots,i_{l-1}-l+1),\no\\
\gamma&:=&(0,i_{1}(N-k),i_{2}(N-k),\cdots,i_{l-1}(N-k)),\no\\
\epsilon_{1}&:=&(1,0,0,0,\cdots,0),\no\\
\epsilon_{2}&:=&(1,1,0,0,\cdots,0),\no\\
\epsilon_{3}&:=&(1,1,1,0,\cdots,0),\no\\
&&\cdots\no\\
\epsilon_{l}&:=&(1,1,1,1,\cdots,1).
\end{eqnarray}
Now we define $\delta=(\delta_{1},\cdots,\delta_{l})\in{\bf Z}^{l}$ 
by the formula: 
\begin{equation}
\delta:=\alpha+\beta+\gamma+\sum_{j=1}^{l-1}(m_{i_{j}}-1)\epsilon_{j}+
m_{i_{l}}\epsilon_{l}.
\label{delta}
\end{equation}
With this setup, we state the definition of $\tilde{L}_n^{N,k,d}$:
\begin{defi}
The virtual structure constant $\tilde{L}_{n}^{N,k,d}$ is a rational number which 
is non-zero only if $0\leq n\leq N-1-(N-k)d$. 
It is uniquely 
determined by the initial condition:
\begin{eqnarray}
&&\sum_{n=0}^{k-1}\tilde{L}_{n}^{N,k,1}w^{n}=
k\cdot\prod_{j=1}^{k-1}(jw+(k-j)),   \;\;(N\geq 2k),\no\\
&& \tilde{L}_{n}^{N,k,d}=0,  \;\;(d\geq 2,\; N\geq 2k),
\label{ini1}
\end{eqnarray}
and by the recursive formula:
\begin{equation}
\tilde{L}_{n}^{N,k,d}=\phi(Poly_{d}).
\label{mr}
\end{equation}
In (\ref{mr}), $\phi$ is a ${\bf Q}$-linear map from the ${\bf Q}$-vector
space of the homogeneous polynomials of degree $d-1$ in $x,y,z_{1},\cdots
,z_{d-1}$ to the ${\bf Q}$-vector
space of the weighted homogeneous polynomials of degree $d$ in 
$L_{m}^{N+1,k,d'}$. It is  defined on the basis by: 
\begin{equation}
\phi(x^{m_{0}}y^{m_{d}}z_{i_{1}}^{m_{i_{1}}}\cdots z_{i_{l-1}}^{m_{i_{l-1}}})
=\prod_{j=1}^{l}\tilde{L}_{n+\delta_{j}}^{N+1,k,i_{j}-i_{j-1}}.
\end{equation}
\end{defi}
\subsubsection{Proof}
In order to prove Theorem 1, it is enough for us to prove the following equality.
\begin{eqnarray}
\frac{\tilde{L}_{n}^{N,k,d}}{d}&=&\frac{1}{k}
\sum_{\sigma_{d}\in OP_{d}}\frac{1}{(2\pi\sqrt{-1})^{l(\sigma_{d})+1}\prod_{j=0}^{l(\sigma_{d})}d_{j}}\oint_{C_{0}} \frac{dz_{0}}{(z_{0})^N}\cdots
\oint_{C_{0}} \frac{d z_{l(\sigma_{d})} }{(z_{l(\sigma_{d})})^N}(z_{0})^{N-2-n}\times\no\\
&&\times\prod_{j=1}^{l(\sigma_{d})-1}\frac{1}
{\biggl( \frac{z_{j}-z_{j-1}}{d_{j}}+\frac{z_{j}-z_{j+1}}{d_{j+1}}\biggr)kz_{j}}\prod_{j=1}^{l(\sigma_{d})}
\frac{e(k,d_j;z_{j-1},z_{j})}{t(N,d_j;z_{j-1},z_{j})}(z_{l(\sigma_{d})})^{n-1+(N-k)d}.\no\\
\label{int}
\end{eqnarray}
As we have remarked in Section 2.1.3, the residue integral in (\ref{int}) strongly depends on the order of integration, and we have to take the residues of 
$z_{j}$'s in descending (or ascending) order of subscript $j$.

We prove the above theorem by showing that the r.h.s. of (\ref{int}) satisfies the initial 
condition (\ref{ini1}) and the recursion relation (\ref{mr}). For this purpose, 
we introduce the following lemma:
\begin{lem}
\begin{eqnarray}
&&\frac{1}{k}\sum_{\sigma_{d}\in OP_{d}} \frac{1}{(2\pi\sqrt{-1})^{l(\sigma_{d})+1}
\prod_{j=0}^{l(\sigma_{d})}d_{j}} \oint_{C_{0}} \frac{dz_{0}}{(z_{0})^N}\cdots
\oint_{C_{0}} \frac{d z_{l(\sigma_{d})} }{(z_{l(\sigma_{d})})^N}
(z_{0})^{N-2-n}\times\no\\
&&\times\prod_{j=1}^{l(\sigma_{d})-1}\frac{1}
{\biggl( \frac{z_{j}-z_{j-1}}{d_{j}}+\frac{z_{j}-z_{j+1}}{d_{j+1}}\biggr)kz_{j}}\prod_{j=1}^{l(\sigma_{d})}\frac{e(k,d_j;z_{j-1},z_{j})}
{t(N,d_j;z_{j-1},z_{j})}(z_{l(\sigma_{d})})^{n-1+(N-k)d}=\no\\
&=&\frac{1}{k}\frac{1}{(2\pi\sqrt{-1})^{d+1}}\oint_{C_{0}}dz_{0}\oint_{E_{1}}dz_{1}\cdots\oint_{E_{d-1}}dz_{d-1}\oint_{C_0}dz_{d}
\frac{(z_{0})^{N-2-n}(z_{d})^{n-1+(N-k)d}}{(z_{0})^N(z_{d})^N\prod_{i=1}^{d-1}((z_{i})^N(2z_{i}-z_{i-1}-z_{i+1}))}\times\no\\
&&\times \frac{\prod_{j=1}^{d}e(k,1;z_{j-1},z_{j})}{\prod_{i=1}^{d-1}(kz_{i})} ,
\label{simple}
\end{eqnarray}
where $\frac{1}{2\pi\sqrt{-1}}\oint_{E_{j}},\;(i=1,\cdots,d-1)$ represents $\frac{1}{2\pi\sqrt{-1}}\oint_{C_{(0,\frac{z_{j-1}+z_{j+1}}{2})}}$.
\end{lem}
{\it proof of Lemma 1)} We first pay attention to the fact that $\oint_{E_{j}}dz_j$ decomposed into $\oint_{C_{0}}dz_j+\oint_{C_{\frac{z_{j-1}+z_{j+1}}{2}}}dz_j$ for $j=1,2,\cdots,d-1$. Therefore, the r.h.s. of
(\ref{simple}) can be rewritten as follows:
\begin{eqnarray}
&&\frac{1}{k}
\frac{1}{(2\pi\sqrt{-1})^{d+1}}\sum_{n=0}^{d-1}\sum_{1\leq j_1<j_2<\cdots<j_n\leq d-1}\oint_{C_{0}}dz_{0}\cdots\oint_{C_{0}}dz_{j_1-1}
\oint_{C_{\frac{z_{j_1-1}+z_{j_1+1}}{2}}}dz_{j_1}\oint_{C_{0}}dz_{j_1+1}\cdots\oint_{C_{0}}dz_{j_2-1}\times\no\\
&&\times
\oint_{C_{\frac{z_{j_2-1}+z_{j_2+1}}{2}}}dz_{j_2}\oint_{C_{0}}dz_{j_2+1}\cdots\cdots \oint_{C_{0}}dz_{j_n-1}
\oint_{C_{\frac{z_{j_n-1}+z_{j_n+1}}{2}}}dz_{j_n}\oint_{C_{0}}dz_{j_n+1}\cdots\oint_{C_{0}}dz_d\times\no\\
&&\times
\frac{(z_{0})^{N-2-n}(z_{d})^{n-1+(N-k)d}}{(z_{0})^N(z_{d})^N\prod_{i=1}^{d-1}((z_{i})^N(2z_{i}-z_{i-1}-z_{i+1}))}
\cdot\frac{\prod_{j=1}^{d}e(k,1;z_{j-1},z_{j})}{\prod_{i=1}^{d-1}(kz_{i})}.
\label{dec0}
\end{eqnarray}
Then we change integration variables of the summand that corresponds to $1\leq j_1<j_2<\cdots<j_n\leq d-1$ 
as follows:
\begin{eqnarray}
u_{i}=z_{i}\;&&\mbox{if}\;i\notin\{j_1,j_2,\cdots,j_n\},\no\\
u_{i}=2z_{i}-z_{i-1}-z_{i+1}\;&&\mbox{if}\;i\in\{j_1,j_2,\cdots,j_n\}.
\label{ch1}
\end{eqnarray}
Let $\{i_{1},i_{2},\cdots,i_{l-1}\}$ be $\{1,2,\cdots,d-1\}-\{j_1,j_2,\cdots,j_n\}$ where 
\begin{eqnarray}
0=:i_0<i_1<i_2<\cdots<i_{l-1}<i_l:=d,\;\;l=d-n.
\end{eqnarray}
Inversion of (\ref{ch1}) results in,
\begin{eqnarray}
z_{j}(u_{*})=u_{j}\;&&\mbox{if}\;j\in \{i_{0},i_{1},i_{2},\cdots,i_{l-1},i_{l}\},\no\\
z_{j}(u_{*})=\frac{(i_{m}-j)u_{i_{m-1}}+(j-i_{m-1})u_{i_m}+\sum_{h=i_{m-1}+1}^{i_{m}-1}C_{h}^{j}u_h}{i_m-i_{m-1}}
&&\mbox{if}\;i_{m-1}+1\leq j\leq i_{m}-1,
\label{invch}
\end{eqnarray}
where $C_{h}^{j}$ is some positive integer. The Jacobian of this coordinate change is given by,
\begin{eqnarray}
\frac{1}{\prod_{m=1}^{l}(i_{m}-i_{m-1})}.
\label{jac}
\end{eqnarray}
In this way, the term corresponding to $1\leq j_1<j_2<\cdots<j_n\leq d-1$ in (\ref{dec0}) can be rewritten as follows:
\begin{eqnarray}
&&\frac{1}{k}\frac{1}{(2\pi\sqrt{-1})^{d+1}}
\frac{1}{\prod_{m=1}^{l}(i_{m}-i_{m-1})}\oint_{C_{0}}du_{0}\oint_{C_{0}}du_{1}\cdots\oint_{C_{0}}du_{d}
\frac{(z_{0}(u_*))^{N-2-n}(z_{d}(u_*))^{n-1+(N-k)d}}{(z_{0}(u_{*}))^N(z_{d}(u_{*}))^N}\times\no\\
&&\times\frac{1}{\prod_{m=1}^{l-1}\bigl((z_{i_m}(u_{*}))^N(2z_{i_m}(u_{*})-z_{i_m-1}(u_*)-z_{i_m+1}(u_*))\bigr)}
\cdot\frac{1}{\prod_{m=1}^{l}\prod_{j=i_{m-1}+1}^{i_m-1}u_j\cdot(z_{j}(u_{*}))^N}\times\no\\
&&\times\prod_{m=1}^{l-1}\frac{1}{kz_{i_m}(u_{*})}\prod_{m=1}^{l}
\frac{\prod_{j=i_{m-1}+1}^{i_{m}}e(k,1;z_{j-1}(u_{*}),z_j(u_{*}))}{\prod_{j=i_{m-1}+1}^{i_m-1}kz_{j}(u_*)}.
\label{urep}
\end{eqnarray}
Looking at (\ref{urep}), we observe that the integrand has only a simple pole at $u_{j_{h}}=0$ $(h=1,2,\cdots,n)$. 
Therefore, we can take the residue of $u_{j_h}$ before $u_{i_m}$ $(m=0,1,\cdots,l)$. 
 After this operation, (\ref{invch})
reduces to,  
\begin{eqnarray}
z_{j}(u_{*})=u_{j}\;&&\mbox{if}\;j\in \{i_{0},i_{1},i_{2},\cdots,i_{l-1},i_{l}\},\no\\
z_{j}(u_{*})=\frac{(i_{m}-j)u_{i_{m-1}}+(j-i_{m-1})u_{i_m}}{i_m-i_{m-1}}
&&\mbox{if}\;i_{m-1}+1\leq j\leq i_{m}-1.
\label{redinvch}
\end{eqnarray}
With (\ref{redinvch}) and some algebra, we can easily derive,
\begin{eqnarray} 
2z_{i_m}(u_*)-z_{i_m-1}(u_*)-z_{i_m+1}(u_*)&=&\frac{u_{i_m}-u_{i_{m-1}}}{i_m-i_{m-1}}+\frac{u_{i_m}-u_{i_{m+1}}}{i_{m+1}-i_{m}},\no\\
\prod_{m=1}^{l}\prod_{j=i_{m-1}+1}^{i_m-1}(z_{j}(u_{*}))^N&=&\prod_{m=1}^{l}t(N,i_m-i_{m-1};u_{i_{m-1}},u_{i_m}),\no\\
\frac{\prod_{j=i_{m-1}+1}^{i_{m}}e(k,1;z_{j-1}(u_{*}),z_j(u_{*}))}{\prod_{j=i_{m-1}+1}^{i_m-1}kz_{j}(u_*)}
&=&e(k,i_m-i_{m-1};u_{i_{m-1}},u_{i_m}).
\label{reduce}
\end{eqnarray}
And (\ref{urep}) equals, 
\begin{eqnarray}
&&\frac{1}{k}\frac{1}{(2\pi\sqrt{-1})^{l+1}}\frac{1}{\prod_{m=1}^{l}(i_{m}-i_{m-1})}
\oint_{C_{0}}du_{0}\oint_{C_{0}}du_{i_1}\oint_{C_{0}}du_{i_2}\cdots\oint_{C_{0}}du_{i_l}
\frac{u_0^{N-2-n}u_{d}^{n-1+(N-k)d}}{(u_{0})^N(u_{d})^N}\times\no\\
&&\times\frac{1}{\prod_{m=1}^{l-1}\bigl((u_{i_m})^N(\frac{u_{i_m}-u_{i_{m-1}}}{i_m-i_{m-1}}+\frac{u_{i_m}-u_{i_{m+1}}}{i_{m+1}-i_{m}})\bigr)
ku_{i_m}}
\cdot\prod_{m=1}^{l}\frac{e(k,i_m-i_{m-1};u_{i_{m-1}},u_{i_m})}{t(N,i_m-i_{m-1};u_{i_{m-1}},u_{i_m})}.
\label{urep2}
\end{eqnarray}
By setting $d_{m}=i_{m}-i_{m-1}$ and $z_{m}=u_{i_m}$, (\ref{urep2}) turns out be the summand of the l.h.s. of (\ref{simple}) corresponding to 
$\sigma_{d}=(d_1,d_2,\cdots,d_l)\;(l=l(\sigma_{d}))$. $\Box$

Next, we note the following elementary identity:
\begin{eqnarray}
&&\prod_{j=0}^{l(\sigma_{d})}\frac{1}{(z_j)^N}
\frac{1}{\prod_{i=1}^{d-1}(2z_{i}-z_{i-1}-z_{i+1})}\frac{\prod_{j=1}^{d}e(k,1;z_{j-1},z_{j})}{\prod_{i=1}^{d-1}(kz_{i})}
=\no\\
&&=\prod_{j=0}^{l(\sigma_{d})}\frac{1}{(z_j)^{N+1}}
\frac{1}{\prod_{i=1}^{d-1}(2z_{i}-z_{i-1}-z_{i+1})}\frac{\prod_{j=1}^{d}e(k,1;z_{j-1},z_{j})}{\prod_{i=1}^{d-1}(kz_{i})}(z_{0}z_{1}\cdots z_{d}).
\label{reduction}
\end{eqnarray}
(\ref{reduction}) tells us that the recursive formula (\ref{mr}) for arbitrary $d$ can be derived by sufficiently decomposing $z_0 z_1\cdots z_d$.
Let $r_{i}(z_*)$ be $2z_{i}-z_{i-1}-z_{i+1}$ $(i=1,\cdots,d-1)$. 
We introduce here the following decomposition of $\prod_{j=1}^{d-1}z_{j}$:  
\begin{equation}
\prod_{j=1}^{d-1}z_{j} =\sum_{k=0}^{d-1}\sum_{1\leq i_{1}<i_{2}<\cdots<i_{k}\leq d-1}f_{(i_{1},\cdots,i_{k})}(z_{0},z_{d},z_{i_{1}},\cdots,z_{i_{k}})\prod_{j=1}^{k}r_{i_{j}}(z_*).
\label{dec1}
\end{equation}
where $f_{(i_{1},\cdots,i_{k})}(z_{0},z_{d},z_{i_{1}},\cdots,z_{i_{k}})$ is a homogeneous polynomial in $z_{0},z_{d},z_{i_{1}},\cdots,z_{i_{k}}$
of degree $d-1-k$. 
\begin{lem}
\begin{eqnarray}
&&f_{(i_{1},\cdots,i_{k})}(z_{0},z_{d},z_{i_{1}},\cdots,z_{i_{k}})=\no\\
&&=\biggl(\prod_{j=0}^{k}(i_{j+1}-i_{j})\biggr)\cdot\frac{1}{(2\pi\sqrt{-1})^{d+1}}
\oint_{C_{z_{0}}}\frac{du_{0}}{u_0-z_0}\oint_{C_{z_{d}}}\frac{du_{d}}{u_d-z_d}
\oint_{D_{1}}du_{1}\cdots\oint_{D_{d-1}}du_{d-1}\prod_{j=1}^{d-1}\frac{u_{j}}{r_{j}(u_*)}\prod_{j=1}^{k}\frac{1}{u_{i_j}-z_{i_j}}. \no\\
\label{res1}
\end{eqnarray}
\end{lem}
In (\ref{res1}), the r.h.s does \it not \normalfont depend on order of integration, because residue integral in (\ref{res1}) 
takes all possible residues of each variable. 
\\ 
{\it proof of lemma 2)}  We first show that the decomposition in (\ref{dec1}) does exist. 
As a first step, we express $z_i$ $(i=1,2,\cdots,d-1)$ as a linear combination of $z_{0}$, $z_{d}$ and $
r_i$ $(i=1,2,\cdots,d-1)$:
\begin{eqnarray}
z_{i}=\frac{(d-i)z_0+iz_d}{d}+\sum_{j=1}^{d-1}C_{i}^{j}r_j,
\label{lin0}
\end{eqnarray}
where $C_i^j$ is some positive rational number. Insertion of the above expression into $z_{1}z_{2}\cdots z_{d-1}$ results 
in the following expression:
\begin{eqnarray}
z_{1}z_{2}\cdots z_{d-1}=\sum_{k=0}^{d-1}\sum_{1\leq i_{1}<\cdots<i_{k}\leq d-1}
\sum_{m_{i_1},\cdots,m_{i_k}\geq 1}g_{(i_1,\cdots,i_k)}^{(0),(m_{i_1},\cdots,m_{i_k})}
(z_{0},z_{d},z_{i_1},\cdots,z_{i_k})\prod_{j=1}^{k}(r_{i_j})^{m_{i_j}},
\label{step1}
\end{eqnarray}
where $g_{(i_1,\cdots,i_k)}^{(0),(m_{i_1},\cdots,m_{i_k})}
(z_{0},z_{d},z_{i_1},\cdots,z_{i_k})$ is a homogeneous polynomial in $z_0$, $z_d$ and $z_{i_j}$ $(j=1,2,\cdots,k)$ 
of degree $d-1-\sum_{j=1}^{k}m_{i_j}$(actually, it depends only on $z_{0}$ and $z_{d}$ at this step). At this stage, 
we focus on terms of the following type:
\begin{eqnarray}
g_{(i_1)}^{(0),(m_{i_1})}(z_{0},z_{d},z_{i_1})(r_{i_1})^{m_{i_1}}.
\label{pay1}
\end{eqnarray}
Then we express $z_{j}$ $(j\neq i_1)$ as a linear combination of $z_{0}$, $z_{d}$, $z_{i_1}$ and $r_{k}$ $(k\neq i_1)$. Inserting 
this expression into $r_{i_1}=2z_{i_1}-z_{i_1-1}-z_{i_1+1}$, we can express $r_{i_1}$ as a linear combination of these variables.
Let $l_{i_1}(z_{0},z_{d},z_{i_1},r_{k}\;(k\neq i_1))$ be the resulting expression of $r_{i_1}$. Then we rewrite the terms given 
in (\ref{pay1}) in the following form:
\begin{eqnarray}
g_{(i_1)}^{(0),(m_{i_1})}(z_{0},z_{d},z_{i_1})r_{i_{1}}(l_{i_1}(z_{0},z_{d},z_{i_1},r_{k}\;(k\neq i_1)))^{m_{i_1}-1}.
\label{rewrite1}
\end{eqnarray}
After this operation, we obtain a new expression for $z_1z_2\cdots z_{d-1}$:
\begin{eqnarray}
z_{1}z_{2}\cdots z_{d-1}=\sum_{k=0}^{d-1}\sum_{1\leq i_{1}<\cdots<i_{k}\leq d-1}
\sum_{m_{i_1},\cdots,m_{i_k}\geq 1}g_{(i_1,\cdots,i_k)}^{(1),(m_{i_1},\cdots,m_{i_k})}
(z_{0},z_{d},z_{i_1},\cdots,z_{i_k})\prod_{j=1}^{k}(r_{i_j})^{m_{i_j}}.
\label{step2}
\end{eqnarray}
In the above expression, terms of type:
\begin{eqnarray} 
g_{(i_1)}^{(1),(m_{i_1})}(z_{0},z_{d},z_{i_1})(r_{i_1})^{m_{i_1}},\;\;\;(m_{i_1}\geq 2)
\label{dont2}
\end{eqnarray}
do not appear. At this stage, we look at terms of the following type:
\begin{eqnarray}
g_{(i_1,i_2)}^{(1),(m_{i_1},m_{i_2})}(z_{0},z_{d},z_{i_1},z_{i_2})(r_{i_1})^{m_{i_1}}(r_{i_2})^{m_{i_2}}.
\label{pay2}
\end{eqnarray}
We then express $r_{i_1}$ and $r_{i_2}$ as linear combinations of $z_{0}$, $z_{d}$, $z_{i_1}$, $z_{i_2}$ and 
$r_{k}$ $(k\neq i_1, i_2)$ in the same way as the previous step. Let 
$l_{i_1}(z_{0},z_{d},z_{i_1},z_{i_2},r_{k}\;(k\neq i_1,i_2))$ and 
$l_{i_2}(z_{0},z_{d},z_{i_1},z_{i_2},r_{k}\;(k\neq i_1,i_2))$ be the resulting expressions. Next, we rewrite the terms 
given in (\ref{pay2}) in the form:
\begin{eqnarray}
g_{(i_1,i_2)}^{(1),(m_{i_1},m_{i_2})}(z_{0},z_{d},z_{i_1},z_{i_2})r_{i_{1}}r_{i_2}
(l_{i_1}(z_{0},z_{d},z_{i_1},z_{i_2},r_{k}\;(k\neq i_1,i_2)))^{m_{i_1}-1}
(l_{i_2}(z_{0},z_{d},z_{i_1},z_{i_2},r_{k}\;(k\neq i_1,i_2)))^{m_{i_2}-1}.\no\\
\label{rewrite2}
\end{eqnarray}
After this operation, we again obtain new expression of  $z_1z_2\cdots z_{d-1}$:
\begin{eqnarray}
z_{1}z_{2}\cdots z_{d-1}=\sum_{k=0}^{d-1}\sum_{1\leq i_{1}<\cdots<i_{k}\leq d-1}
\sum_{m_{i_1},\cdots,m_{i_k}\geq 1}g_{(i_1,\cdots,i_k)}^{(2),(m_{i_1},\cdots,m_{i_k})}
(z_{0},z_{d},z_{i_1},\cdots,z_{i_k})\prod_{j=1}^{k}(r_{i_k})^{m_{i_k}}.
\label{step3}
\end{eqnarray}
In the above expression, the terms of the following types:
\begin{eqnarray}
&&g_{(i_1)}^{(2),(m_{i_1})}(z_{0},z_{d},z_{i_1})(r_{i_1})^{m_{i_1}},\;\;\;(m_{i_1}\geq 2)\no\\
&&g_{(i_1,i_2)}^{(2),(m_{i_1},m_{i_2})}(z_{0},z_{d},z_{i_1},z_{i_2})(r_{i_1})^{m_{i_1}}(r_{i_2})^{m_{i_2}},
\;\;\;(m_{i_1}\geq 2\;\;\mbox{or}\;\;m_{i_2}\geq 2),
\label{dont3}
\end{eqnarray}
do not appear. In general, we can inductively construct a new expression of $z_1z_2\cdots z_{d-1}$:
\begin{eqnarray}
z_{1}z_{2}\cdots z_{d-1}=\sum_{k=0}^{d-1}\sum_{1\leq i_{1}<\cdots<i_{k}\leq d-1}
\sum_{m_{i_1},\cdots,m_{i_k}\geq 1}g_{(i_1,\cdots,i_k)}^{(h),(m_{i_1},\cdots,m_{i_k})}
(z_{0},z_{d},z_{i_1},\cdots,z_{i_k})\prod_{j=1}^{k}(r_{i_j})^{m_{i_j}}.
\label{steph}
\end{eqnarray}
by rewriting the terms of the type:
\begin{eqnarray}
g_{(i_1,\cdots,i_h)}^{(h-1),(m_{i_1},\cdots,m_{i_h})}
(z_{0},z_{d},z_{i_1},\cdots,z_{i_h})\prod_{j=1}^{h}(r_{i_j})^{m_{i_j}},
\;\;(m_{i_1}\geq2\:\;\mbox{or}\;\;\cdots\;\;\mbox{or}\;\;m_{i_h}\geq2), 
\label{rewriteh}
\end{eqnarray}
in the same way as the first two steps. Finally, the expression:
\begin{eqnarray}
z_{1}z_{2}\cdots z_{d-1}=\sum_{k=0}^{d-1}\sum_{1\leq i_{1}<\cdots<i_{k}\leq d-1}
\sum_{m_{i_1},\cdots,m_{i_k}\geq 1}g_{(i_1,\cdots,i_k)}^{(d-1),(m_{i_1},\cdots,m_{i_k})}
(z_{0},z_{d},z_{i_1},\cdots,z_{i_k})\prod_{j=1}^{k}(r_{i_j})^{m_{i_j}}.
\label{stepf}
\end{eqnarray}
is nothing but the desired decomposition.

We have shown that the decomposition (\ref{dec1}) does exist. Therefore, we can insert,
\begin{eqnarray}
\prod_{j=1}^{d-1}u_{j} =\sum_{m=0}^{d-1}\sum_{1\leq h_{1}<h_{2}<\cdots<h_{m}\leq d-1}f_{(h_{1},\cdots,h_{m})}(u_{0},u_{d},u_{h_{1}},\cdots,u_{h_{m}})\prod_{j=1}^{m}r_{h_{j}}(u_*),
\label{dec2}
\end{eqnarray}
into the r.h.s. of (\ref{res1}). It then becomes, 
\begin{eqnarray}
&&\sum_{m=0}^{d-1}\sum_{1\leq h_{1}<h_{2}<\cdots<h_{m}\leq d-1}\biggl(\prod_{j=0}^{k}(i_{j+1}-i_{j})\biggr)\cdot\frac{1}{(2\pi\sqrt{-1})^{d+1}}
\oint_{C_{z_{0}}}\frac{du_{0}}{u_0-z_0}\oint_{C_{z_{d}}}\frac{du_{d}}{u_d-z_d}
\oint_{D_{1}}du_{1}\cdots\oint_{D_{d-1}}du_{d-1}\times\no\\
&&\times f_{(h_{1},\cdots,h_{m})}(u_{0},u_{d},u_{h_{1}},\cdots,u_{h_{m}})
\left(\prod_{j\in\{1,2,\cdots,d-1\}-\{h_{1},h_{2},\cdots,h_{m}\}}\frac{1}{r_{j}(u_*)}\right)\prod_{j=1}^{k}\frac{1}{u_{i_j}-z_{i_j}}.
\label{ins1}
\end{eqnarray}
At this stage, we use the fact that the above expression does not depend on order of integration.
If $(\{1,2,\cdots,d-1\}-\{h_{1},h_{2},\cdots,h_{m}\})\cap\{i_1,i_2,\cdots,i_k\}\neq\emptyset$, the summand corresponding 
to $1\leq h_{1}<h_{2}<\cdots<h_{m}\leq d-1$ vanishes because for $j\in (\{1,2,\cdots,d-1\}-\{h_{1},h_{2},\cdots,h_{m}\})\cap\{i_1,i_2,\cdots,i_k\}$, 
\begin{eqnarray}
\oint_{D_{j}}du_j\frac{1}{(2u_j-u_{j-1}-u_{j+1})(u_j-z_j)}=0.
\label{van1}
\end{eqnarray} 
If $(\{1,2,\cdots,d-1\}-\{h_{1},h_{2},\cdots,h_{m}\})\cup\{i_1,i_2,\cdots,i_k\}\neq \{1,2,\cdots,d-1\}$, it also vanishes 
because the integrand has no poles of the variable $u_{j}$ 
$(j\notin (\{1,2,\cdots,d-1\}-\{h_{1},h_{2},\cdots,h_{m}\})\cup\{i_1,i_2,\cdots,i_k\})$. In this way, only the summand 
that satisfies $\{h_1,h_2,\cdots,h_m\}=\{i_1,i_2,\cdots,i_k\}$ survives. Hence (\ref{ins1}) becomes,
\begin{eqnarray}
&&\biggl(\prod_{j=0}^{k}(i_{j+1}-i_{j})\biggr)\cdot\frac{1}{(2\pi\sqrt{-1})^{d+1}}
\oint_{C_{z_{0}}}\frac{du_{0}}{u_0-z_0}\oint_{C_{z_{d}}}\frac{du_{d}}{u_d-z_d}
\oint_{D_{1}}du_{1}\cdots\oint_{D_{d-1}}du_{d-1}\times\no\\
&&\times f_{(i_{1},\cdots,i_{k})}(u_{0},u_{d},u_{i_{1}},\cdots,u_{i_{k}})
\left(\prod_{j\in\{1,2,\cdots,d-1\}-\{i_{1},i_{2},\cdots,i_{k}\}}\frac{1}{r_{j}(u_*)}\right)\prod_{j=1}^{k}\frac{1}{u_{i_j}-z_{i_j}}.
\label{ins2}
\end{eqnarray}
We then perform the following coordinate change of integration variables:
\begin{eqnarray}
w_{j}&=&u_{j}\;\;\mbox{if}\; j\in\{i_1,i_2,\cdots,i_k\}\cup\{0,d\},\no\\
w_{j}&=&r_{j}(u_{*})=2u_j-u_{j-1}-u_{j+1}\;\;\mbox{if}\; j\in\{1,2,\cdots,d-1\}-\{i_1,i_2,\cdots,i_k\}.
\label{chpr} 
\end{eqnarray}
Since the Jacobian of the above coordinate change is given by $\prod_{j=0}^{k}\frac{1}{(i_{j+1}-i_{j})}$, (\ref{ins2}) becomes, 
 \begin{eqnarray}
&&\frac{1}{(2\pi\sqrt{-1})^{d+1}}
\oint_{C_{z_{0}}}\frac{dw_{0}}{w_0-z_0}\oint_{C_{z_{d}}}\frac{dw_{d}}{w_d-z_d}
\prod_{j=1}^{k}\oint_{C_{z_{i_j}}}\frac{dw_{i_j}}{w_{i_j}-z_{i_j}}\left(
\prod_{j\in\{1,2,\cdots,d-1\}-\{i_1,i_2,\cdots,i_k\}}\oint_{C_{0}}\frac{dw_{j}}{w_j}\right)\times\no\\
&&\times f_{(i_{1},\cdots,i_{k})}(w_{0},w_{d},w_{i_{1}},\cdots,w_{i_{k}})=
f_{(i_{1},\cdots,i_{k})}(z_{0},z_{d},z_{i_{1}},\cdots,z_{i_{k}}).
\label{ins3}
\end{eqnarray} 
$\Box$
\\
{\it proof of Theorem 1)}
As our first step, we write down the explicit form of the recursive formula (\ref{mr}) used in the definition of 
$\tilde{L}_{n}^{N,k,d}$. Since $\frac{u_j}{u_j-z_j}=1+\frac{z_j}{u_j-z_j}$, we can rewrite $Poly_d$ in (\ref{trial2}) 
as follows:
\begin{eqnarray}
&&Poly_{d}(z_0,z_d,z_1,z_2,\cdots,z_{d-1})=\no\\
&&=\frac{d}{(2\pi\sqrt{-1})^{d-1}}\oint_{D_1}du_1\cdots\oint_{D_{d-1}}du_{d-1}
\prod_{j=1}^{d-1}\left(\frac{u_j}{(2u_j-u_{j-1}-u_{j+1})}(1+\frac{z_j}{u_j-z_j})\right)=\no\\
&&=\sum_{l=1}^{d}\sum_{1\leq i_1<\cdots <i_{l-1}\leq d-1}\frac{d}{\prod_{j=1}^{l}(i_{j}-i_{j-1})}
\left(\prod_{j=1}^{l-1}z_{i_j}\right)f_{(i_1,\cdots,i_{l-1})}(z_{0},z_{d},z_{i_1},\cdots,z_{i_{l-1}}),
\label{recr1}
\end{eqnarray}
where we formally set $i_{0}$ (resp. $i_l$) to $0$ (resp. $d$). In deriving (\ref{recr1}), we used Lemma 2.\\
Since $f_{(i_1,\cdots,i_{l-1})}(z_{0},z_{d},z_{i_1},\cdots,z_{i_{l-1}})$ is a homogeneous polynomial of 
degree $d-l$, it can be expanded as follows:
\begin{eqnarray}
f_{(i_1,\cdots,i_{l-1})}(z_{0},z_{d},z_{i_1},\cdots,z_{i_{l-1}})
=\sum_{m_j\geq 0,\;\sum_{j=0}^{l}m_j=d-l}C_{(i_1,\cdots,i_{l-1})}^{(m_0,m_1,\cdots,m_l)}\prod_{j=0}^{l}(z_{i_j})^{m_j},
\label{rex}
\end{eqnarray}
where $C_{(i_1,\cdots,i_{l-1})}^{(m_0,m_1,\cdots,m_l)}$ is some rational number. With these notations, $Poly_d$ is explicitly 
given by,
\begin{eqnarray}
&&Poly_{d}(z_0,z_d,z_1,z_2,\cdots,z_{d-1})=\no\\
&&=\sum_{l=1}^{d}\sum_{1\leq i_1<\cdots <i_{l-1}\leq d-1}\frac{d}{\prod_{j=1}^{l}(i_{j}-i_{j-1})}
\sum_{m_j\geq 0,\;\sum_{j=0}^{l}m_j=d-l}C_{(i_1,\cdots,i_{l-1})}^{(m_0,m_1,\cdots,m_l)}(z_{0})^{m_{0}}(z_{d})^{m_l}
\prod_{j=1}^{l-1}(z_{i_j})^{m_j+1}.\no\\
\label{expoly}
\end{eqnarray} 
Using the definition of ${\bf Q}$-linear map $\phi$ in Definition 1, we obtain an explicit form of the recursive formula (\ref{mr}):
\begin{eqnarray}
&&\frac{\tilde{L}_{n}^{N,k,d}}{d}=\sum_{l=1}^{d}\sum_{1\leq i_1<\cdots <i_{l-1}\leq d-1}
\sum_{m_j\geq 0,\;\sum_{j=0}^{l}m_j=d-l}C_{(i_1,\cdots,i_{l-1})}^{(m_0,m_1,\cdots,m_l)}\prod_{j=1}^{l}
\left(\frac{\tilde{L}^{N+1,k,i_j-i_{j-1}}_{n+i_{j-1}(N-k+1)+l-d-j+1+\sum_{h=j}^{l}m_h}}{(i_j-i_{j-1})}\right).\no\\
\label{rexf}
\end{eqnarray}
On the other hand, let $T_n^{N,k,d}$ be the r.h.s of (\ref{simple}), i.e.,
\begin{eqnarray}
T_n^{N,k,d}&:=&\frac{1}{k}\frac{1}{(2\pi\sqrt{-1})^{d+1}}\oint_{C_{0}}dz_{0}\oint_{E_{1}}dz_{1}\cdots\oint_{E_{d-1}}dz_{d-1}
\oint_{C_{0}}dz_{d}
\frac{(z_{0})^{N-2-n}(z_{d})^{n-1+(N-k)d}}{(z_{0})^N(z_{d})^N\prod_{i=1}^{d-1}((z_{i})^N(2z_{i}-z_{i-1}-z_{i+1}))}\times\no\\
&&\times \frac{\prod_{j=1}^{d}e(k,1;z_{j-1},z_{j})}{\prod_{i=1}^{d-1}(kz_{i})}.
\label{tdef}
\end{eqnarray}
To prove the assertion of Theorem 1, it suffices to show that $T_n^{N,k,d}$ satisfies the same initial conditions and recursive 
formulas as the those of $\frac{\tilde{L}_{n}^{N,k,d}}{d}$. By looking back at (\ref{reduction}) and (\ref{dec1}), we can deduce,
\begin{eqnarray}
&&T_n^{N,k,d}=\no\\
&=&\sum_{l=1}^{d}\sum_{1\leq i_1<\cdots <i_{l-1}\leq d-1}
\sum_{m_j\geq 0,\;\sum_{j=0}^{l}m_j=d-l}\frac{1}{k}\frac{1}{(2\pi\sqrt{-1})^{d+1}}\oint_{C_{0}}dz_{0}\oint_{E_{1}}dz_{1}\cdots\oint_{E_{d-1}}dz_{d-1}\oint_{C_{0}}dz_{d}
\times\no\\
&&\times C_{(i_1,\cdots,i_{l-1})}^{(m_0,m_1,\cdots,m_l)}\frac{(z_{0})^{N-1-n+m_0}(z_{d})^{n+(N-k)d+m_l}\prod_{j=1}^{l-1}((z_{i_j})^{m_j}r_{i_j}(z_*))}{(z_{0})^{N+1}(z_{d})^{N+1}\prod_{i=1}^{d-1}((z_{i})^{N+1}r_{i}(z_{*}))}\frac{\prod_{j=1}^{d}e(k,1;z_{j-1},z_{j})}{\prod_{i=1}^{d-1}(kz_{i})}=\no\\
&=&\sum_{l=1}^{d}\sum_{1\leq i_1<\cdots <i_{l-1}\leq d-1}
\sum_{m_j\geq 0,\;\sum_{j=0}^{l}m_j=d-l}C_{(i_1,\cdots,i_{l-1})}^{(m_0,m_1,\cdots,m_l)}\times\\
&&\times\frac{1}{k}\frac{1}{(2\pi\sqrt{-1})^{i_1+1}}\oint_{C_{0}}dz_{0}\cdots\oint_{E_{i_1}}dz_{i_1}
\frac{(z_0)^{N-1-n+m_0}}{(z_{0})^{N+1}(z_{i_1})^{N+1}\prod_{j=1}^{i_1-1}(r_{j}(z_{*})(z_j)^{N+1})}
\frac{\prod_{j=1}^{i_1}e(k,1;z_{j-1},z_{j})}{\prod_{j=1}^{i_1-1}(kz_{j})}\times\no\\
&&\times\frac{1}{k}\frac{1}{(2\pi\sqrt{-1})^{i_2-i_1}}\oint_{E_{i_1+1}}dz_{i_1+1}\cdots\oint_{E_{i_2}}dz_{i_2}
\frac{(z_{i_1})^{m_1-1}}{(z_{i_2})^{N+1}\prod_{j=i_1+1}^{i_2-1}(r_{j}(z_{*})(z_j)^{N+1})}
\frac{\prod_{j=i_1+1}^{i_2}e(k,1;z_{j-1},z_{j})}{\prod_{j=i_1}^{i_2-1}(kz_{j})}\times\no\\
&&\times\frac{1}{k}\frac{1}{(2\pi\sqrt{-1})^{i_3-i_2}}\oint_{E_{i_2+1}}dz_{i_1+1}\cdots\oint_{E_{i_3}}dz_{i_3}
\frac{(z_{i_2})^{m_2-1}}{(z_{i_3})^{N+1}\prod_{j=i_2+1}^{i_3-1}(r_{j}(z_{*})(z_j)^{N+1})}
\frac{\prod_{j=i_2+1}^{i_3}e(k,1;z_{j-1},z_{j})}{\prod_{j=i_2}^{i_3-1}(kz_{j})}\times\no\\
&&\times \hspace{6.8cm}\cdots\cdots\cdots\hspace{6.8cm}\times\no\\
&&\times\frac{1}{k}\frac{1}{(2\pi\sqrt{-1})^{i_l-i_{l-1}}}\oint_{E_{i_{l-1}+1}}dz_{i_{l-1}+1}\cdots\oint_{C_0}dz_{d}
\frac{(z_{i_{l-1}})^{m_{l-1}-1}(z_{d})^{n+(N-k)d+m_{l}}　}{(z_{i_{l-1}})^{N+1}(z_{d})^{N+1}\prod_{j=i_{l-1}+1}^{i_l-1}(r_{j}(z_{*})(z_j)^{N+1})}
\times\no\\
&&\times\frac{\prod_{j=i_{l-1}+1}^{i_l}e(k,1;z_{j-1},z_{j})}{\prod_{j=i_{l-1}}^{i_l-1}(kz_{j})}=\no\\
&=&\sum_{l=1}^{d}\sum_{1\leq i_1<\cdots <i_{l-1}\leq d-1}
\sum_{m_j\geq 0,\;\sum_{j=0}^{l}m_j=d-l}C_{(i_1,\cdots,i_{l-1})}^{(m_0,m_1,\cdots,m_l)}\prod_{j=1}^{l}
T^{N+1,k,i_j-i_{j-1}}_{n+i_{j-1}(N-k+1)+l-d-j+1+\sum_{h=j}^{l}m_h}.
\label{prrec}
\end{eqnarray}
Therefore,  $T_n^{N,k,d}$ indeed satisfies the same recursive formulas as $\frac{\tilde{L}_{n}^{N,k,d}}{d}$.
We can easily confirm that the initial conditions are the same by direct computation. $\Box$

As the final remark in this section, we go back to the formula (\ref{gw1}).
 The result of the computation of (\ref{gw1}) by the localization theorem coincided with the formula in the r.h.s. of (\ref{int}), and 
we concluded in \cite{vs} that the virtual structure constants can be interpreted as intersection numbers of the moduli space 
of polynomial maps $\widetilde{Mp}_{0,2}(N,d)$.  But by combining the r.h.s. of (\ref{simple}) with the relation (\ref{rel}), 
we can obtain an interesting formula:
\begin{eqnarray}
\frac{w({\cal O}_{h^{N-2-n}}{\cal O}_{h^{n-1+(N-k)d}})_{0,d}}{k}=\frac{\tilde{L}_{n}^{N,k,d}}{d} &=
\displaystyle{\frac{1}{k}}&\int_{ \widetilde{Mp}_{0,2}(N,d) }(h_{0})^{N-2-n}(h_{d})^{n-1+(N-k)d}\frac{\prod_{j=1}^{d}e(k,1;h_{j-1},h_{j})}{\prod_{i=1}^{d-1}(kh_{i})}, \no\\
\label{concl}
\end{eqnarray} 
where we apply normalization:
\begin{eqnarray}
\int_{ \widetilde{Mp}_{0,2}(N,d) }(h_{0})^{N-1}(h_{d})^{N-1}\prod_{j=1}^{d-1}(h_{j})^N=\frac{1}{d}.
\label{normal}
\end{eqnarray}
This formula gives an alternate expression of the virtual structure constant $\frac{\tilde{L}_{n}^{N,k,d}}{d}$ as an intersection number of 
$\widetilde{Mp}_{0,2}(N,d)$.

\newpage
\section{Generalizations to Toric Manifolds with Two K\"ahler Forms}
\subsection{$K_{F_{0}}$}
\subsubsection{ Construction of the Moduli Space $\widetilde{Mp}_{0,2}(F_{0},(d_{a},d_{b}))$}
The Hirzebruch surface $F_0$ is nothing but a product manifold of two ${\bf P}^1$'s. Therefore, 
it is given by, 
\begin{eqnarray}
F_{0}:&=&\{({\bf a},{\bf b})\;|\;{\bf a},{\bf b}\in{\bf C}^{2},{\bf a}, {\bf b}\neq {\bf 0}\}/({\bf C}^{\times})^2,
\end{eqnarray}
where the two ${\bf C}^{\times}$ actions act on ${\bf a}$ and ${\bf b}$ respectively:
\begin{eqnarray}
({\bf a}, {\bf b})\rightarrow (\mu{\bf a}, {\bf b}),\;\;({\bf a}, {\bf b})\rightarrow ({\bf a}, \nu{\bf b}).
\label{f0two}
\end{eqnarray}
Let $\pi_1$ (resp. $\pi_2$) be projection from $F_0$ to the first (resp. the second) ${\bf P}^1$. 
We denote $\pi_1^{*}{\cal O}_{{\bf P}^1}(1)$ (resp. $\pi_2^{*}{\cal O}_{{\bf P}^1}(1)$) by 
${\cal O}_{F_0}(a)$ (resp. ${\cal O}_{F_0}(b)$). Classical cohomology ring of $F_{0}$ is generated 
by two K\"ahler forms $z:=c_1({\cal O}_{F_0}(a))$ and $w:=c_1({\cal O}_{F_0}(b))$. They obey the two 
relations:
\begin{eqnarray}
z^2=0,\;\;w^2=0.
\label{crelf0}
\end{eqnarray}
Integration of $\alpha\in H^{*}(F_{0},{\bf C})$ over $F_0$ is realized as residue integral in $z$ 
and $w$:
\begin{eqnarray}
\int_{F_0}\alpha=\frac{1}{(2\pi\sqrt{-1})^2}\oint_{C_{0}}\frac{dz}{z^2}\oint_{C_{0}}\frac{dw}{w^2}\alpha,
\label{intf0}
\end{eqnarray}
where $\alpha$ in the r.h.s. should be regarded as a polynomial in $z$ and $w$.
Let us consider a polynomial map from $CP^1$ to $F_0$. Since $F_0$ has two K\"ahler forms, it 
is classified by bi-degree $\dd=(d_a, d_b)$. A polynomial map from $CP^1$ to $F_0$ of 
bi-degree $(d_a, d_b)$ is explicitly given as follows:
\begin{eqnarray}
&&p:{\bf C}^{2}\rightarrow {\bf C}^{2}\times{\bf C}^2\no\\
&&p(s,t)= ( \sum_{j=0}^{d_{a}}{\bf a}_{j}s^{d_{a}-j}t^{j} , \sum_{j=0}^{d_{b}}{\bf b}_{j}s^{d_{b}-j}t^{j}),\no\\
&& {\bf a}_{j},{\bf b}_{j}\in{\bf C}^{2}, {\bf a}_{0},{\bf b}_{0},{\bf a}_{d_{a}},{\bf b}_{d_{b}}\neq {\bf 0} .
\label{polyp2}
\end{eqnarray}
The conditions ${\bf a}_{0},{\bf b}_{0},{\bf a}_{d_{a}},{\bf b}_{d_{b}}\neq {\bf 0}$ come from 
requirement that it has a well-defined image in $[(1,0)], [(0,1)]\in CP^1$. 
The moduli space of polynomial maps from $CP^1$ to $F_0$ of bi-degree $(d_a,d_b)$ with two marked points, which we
denote by $Mp_{0,2}(F_{0},\dd)$, is defined as follows:
\begin{eqnarray}
Mp_{0,2}(F_{0},(d_{a},d_{b})) :=\{({\bf a}_{0},\cdots,{\bf a}_{d_{a}},{\bf b}_{0},\cdots,{\bf b}_{d_{b}})\:|\; 
{\bf a}_{j},{\bf b}_{j}\in{\bf C}^{2}, {\bf a}_{0},{\bf b}_{0},{\bf a}_{d_{a}},{\bf b}_{d_{b}}\neq {\bf 0} \}/({\bf C}^{\times})^3,
\label{mpdeff0}
\end{eqnarray} 
In (\ref{mpdeff0}), the three $C^{\times}$ actions are given by, 
\begin{eqnarray}
&& ({\bf a}_{0},\cdots,{\bf a}_{d_{a}},{\bf b}_{0},\cdots,{\bf b}_{d_{b}})  \rightarrow 
(\mu{\bf a}_{0},\cdots,\mu{\bf a}_{d_{a}},{\bf b}_{0},\cdots,{\bf b}_{d_{b}}),\no\\
&&({\bf a}_{0},\cdots,{\bf a}_{d_{a}},{\bf b}_{0},\cdots,{\bf b}_{d_{b}}) \rightarrow 
({\bf a}_{0},\cdots,{\bf a}_{d_{a}},\nu{\bf b}_{0},\cdots,\nu{\bf b}_{d_{b}}),\no\\
&&({\bf a}_{0},\cdots,{\bf a}_{d_{a}},{\bf b}_{0},\cdots,{\bf b}_{d_{b}}) \rightarrow 
({\bf a}_{0},\lambda{\bf a}_{1},{\lambda}^2{\bf a}_{2},\cdots,{\lambda}^{d_{a}}{\bf a}_{d_{a}},
{\bf b}_{0},{\lambda}{\bf b}_{1},{\lambda}^2{\bf b}_{2}\cdots,{\lambda}^{d_{b}}{\bf b}_{d_{b}}).
\end{eqnarray}
The first two actions are induced from the two ${\bf C}^{\times}$ actions in (\ref{f0two}), and 
the third one comes from automorphism group of $CP^1$ fixing two marked points.
We denote the toric compactification of $Mp_{0,2}(F_{0},\dd)$ by $\widetilde{Mp}_{0,2}(F_{0},\dd)$.
In order to compactify $Mp_{0,2}(F_{0},\dd)$, we add the boundary divisor $E_{(i_{a},i_{b})},\;((i_a,i_b)\neq (0,0), (d_{a},d_{b}),\;\; 0\leq i_a\leq d_a, 0\leq i_b\leq d_b)$ that correspond to 
chains of two polynomial maps: 
\begin{eqnarray}
( \sum_{j=0}^{i_{a}}{\bf a}_{j}s_{1}^{i_{a}-j}t_{1}^{j} , \sum_{j=0}^{i_{b}}{\bf b}_{j}s_{1}^{i_{b}-j}t_{1}^{j})
\cup ( \sum_{j=0}^{d_{a}-i_{a}}{\bf a}_{j}s_{2}^{d_{a}-i_{a}-j}t_{2}^{j} , \sum_{j=0}^{d_{b}-i_{b}}{\bf b}_{j}s_{2}^{d_{b}-i_{b}-j}t_{2}^{j}).
\label{chain1}
\end{eqnarray} 
Now, we present an explicit construction of $\widetilde{Mp}_{0,2}(F_{0},\dd)$. To this end, we introduce a partial ordering of 
bi-degree $(d_a, d_b)$ of $F_{0}$:
\begin{eqnarray}
(i_a,i_b)>(j_a,j_b)\stackrel{def.}{\Longleftrightarrow} i_{a}\geq j_{a},\;\; i_{b}\geq j_{b}\;\; \mbox{and}\;\; (i_a,i_b)\neq(j_a,j_b).
\label{order}
\end{eqnarray} 
As in the case of $CP^{N-1}$, $\widetilde{Mp}_{0,2}(F_{0},\dd)$ is given as a toric orbifold with boundary divisor coordinate $u_{(i_a,i_b)}$ 
that corresponds to $E_{(i_a,i_b)}$:
\begin{eqnarray}
\widetilde{Mp}_{0,2}(F_{0},\dd)&=&\{({\bf a}_{0},\cdots,{\bf a}_{d_{a}},{\bf b}_{0},\cdots,{\bf b}_{d_{b}},u_{(1,0)},u_{(2,0)},\cdots,
u_{(i_{a},i_{b})},\cdots,u_{(d_a-1,d_b)})\;|\no\\
&&{\bf a}_{i},{\bf b}_{j}\in{\bf C}^{2},\;u_{(i_a,i_b)}\in{\bf C},\;{\bf a}_{0},{\bf a}_{d_a},{\bf b}_{0},{\bf b}_{d_b}\neq {\bf 0},\;\no\\
&&({\bf a}_{i},\prod_{k=0}^{d_{a}}u_{(i,k)})\neq{\bf 0},\;({\bf b}_{j},\prod_{k=0}^{d_{b}}u_{(k,j)})\neq{\bf 0},\;(1\leq i\leq d_a-1,\;
1\leq j\leq d_b-1),\;\;\no\\
&&(u_{(i_a,i_b)},u_{(j_a,j_b)})\neq(0,0)\;\; \mbox{unless}\;\; (i_a,i_b)<(j_a,j_b)\; \mbox{or}\; (i_a,i_b)>(j_a,j_b) \}/ ({\bf C}^{\times})^{(d_a+1)(d_b+1)+1} .\no\\
\label{mpf0}
\end{eqnarray}
We have to explain the origin of the last two conditions in (\ref{mpf0}), which look a little bit complicated. In this construction, 
$u_{(i_{a},i_{b})}=0$ corresponds to the locus where polynomial maps are split into chains of two polynomial maps given in (\ref{chain1}). 
Therefore, if $u_{(i_{a},i_{b})}=0$, we need ${\bf a}_{i_{a}},{\bf b}_{i_{b}}\neq {\bf 0}$. This explains the meaning of 
the second condition. If $u_{(i_a,i_b)}=u_{(j_a,j_b)}=0$, this corresponds to the locus where polynomial maps split into 
chains of three polynomial maps. Therefore it is impossible unless $(i_a,i_b)<(j_a,j_b)$ or $(i_a,i_b)>(j_a,j_b)$. 
The $({\bf C}^{\times})^{ (d_a+1)(d_b+1)+1 }$ action is given by the $ ((d_a+1)(d_b+1)+1) \times ((d_a+2)(d_b+2)-3) $ weight matrix 
$W_{(d_a,d_b)}$. If $(d_{a},d_{b})=(d,0)$, $W_{(d,0)}$ is given by trivial generalization of $W_{d}$ in (\ref{toric1}):
\begin{eqnarray} 
&&W_{(d,0)}:=\no\\
&&\bordermatrix{
&{\bf a}_{0}&{\bf a}_{1}&{\bf a}_{2}&\cdots&{\bf a}_{d-3}&{\bf a}_{d-2}&{\bf a}_{d-1}&{\bf a}_{d}&{\bf b}_{0}&
u_{(1,0)}&u_{(2,0)}&u_{(3,0)}&\cdots&u_{(d-2,0)}&u_{(d-1,0)}\cr
                                  z_{0}&1&0&0&\cdots&0&0&0&0&0&-1&0&0&\cdots&0&0\cr 
                                  z_{1}&0&1&0&\cdots&0&0&0&0&0&2&-1&0&\cdots&0&0\cr  
                                  z_{2}&0&0&1&\ddots&0&0&0&0&0&-1&2&-1&\ddots&\vdots&0\cr                                  \vdots&\vdots&\vdots&\ddots&\ddots&\ddots&\vdots
                                  &\vdots&\vdots&\vdots&0&\ddots&\ddots&\ddots&0&\vdots\cr 
                       \vdots&\vdots&\vdots&\vdots&\ddots&1&0&0&0&0&\vdots&0&\ddots&\ddots&\ddots&0\cr  
                                 \vdots&\vdots&\vdots&0&\ddots&0&1&0&0&0&0&\vdots&\ddots&-1&2&-1\cr    
                                z_{d-1}&0&0&0&\cdots&0&0&1&0&0&0&0&\cdots&0&-1&2\cr
                                  z_{d}&0&0&0&\cdots&0&0&0&1&0&0&0&\cdots&0&0&-1\cr
                                  w_{0}&0&0&0&\cdots&0&0&0&0&1&0&0&\cdots&0&0&0\cr }.\no\\
\label{toric2}
\end{eqnarray}
$W_{(0,d)}$ is obtained in the same way with the roles of ${\bf a}$ and ${\bf b}$ interchanged. If $d_{a}, d_{b}\geq 1$, the construction of 
$W_{(d_a, d_b)}$ becomes non-trivial. As an example, we present $W_{(1,1)}$: 
\begin{eqnarray}
W_{(1,1)}=\bordermatrix{
&{\bf a}_{0}&{\bf a}_{1}&{\bf b}_{0}&{\bf b}_{1}&u_{(1,0)}&u_{(0,1)}\cr
        z_{0}&1&0&0&0&-1&0\cr
        z_{1}&0&1&0&0&0&-1\cr
        w_{0}&0&0&1&0&0&-1\cr
        w_{1}&0&0&0&1&-1&0\cr
        f_{(1,1)}&0&0&0&0&1&1\cr
}.
\label{mp11}
\end{eqnarray}
Let us see how the above weight matrix works in the definition of $\widetilde{Mp}_{0,2}(F_0,(1,1))$.
We can trivialize the last two entries of $({\bf a}_{0},{\bf a}_{1},{\bf b}_{0},{\bf b}_{1},u_{(1,0)},u_{(0,1)})$ by using 
two of the five ${\bf C}^{\times}$ actions as follows.
\begin{eqnarray} 
&&[({\bf a}_{0},{\bf a}_{1},{\bf b}_{0},{\bf b}_{1},u_{(1,0)},u_{(0,1)})]=
[({\bf a}_{0},\frac{u_{(0,1)}}{u_{(1,0)}}{\bf a}_{1},{\bf b}_{0},{\bf b}_{1},1,1)]=[({\bf a}_{0},{\bf a}_{1},\frac{u_{(0,1)}}{u_{(1,0)}}{\bf b}_{0},{\bf b}_{1},1,1)].
\label{bd11a}
\end{eqnarray}
The second representation corresponds to the polynomial map:
\begin{equation}
({\bf a}_{0}s+\frac{u_{(0,1)}}{u_{(1,0)}}{\bf a}_{1}t,{\bf b}_{0}s+{\bf b}_{1}t),
\label{bd11b}
\end{equation}
and the third representation corresponds to, 
\begin{equation}
({\bf a}_{0}s+{\bf a}_{1}t,\frac{u_{(0,1)}}{u_{(1,0)}}{\bf b}_{0}s+{\bf b}_{1}t).
\label{bd11c}
\end{equation}
If we set $u_{(0,1)}=0$, (\ref{bd11b}) (resp. (\ref{bd11c})) turns into,
\begin{equation}
({\bf a}_{0},{\bf b}_{0}s+{\bf b}_{1}t),\;\;(\mbox{resp.}\;\; ({\bf a}_{0}s+{\bf a}_{1}t,{\bf b}_{1})), 
\end{equation}
by projective equivalence. In this way, the locus given by $u_{(0,1)}=0$ corresponds to the boundary component described by the following 
chain of polynomial map:
\begin{equation}
({\bf a}_{0},{\bf b}_{0}s+{\bf b}_{1}t)\cup ({\bf a}_{0}s+{\bf a}_{1}t,{\bf b}_{1}).
\end{equation}
We can also see that the locus given by $u_{(1,0)}=0$ corresponds to the boundary component described by 
$({\bf a}_{0}s+{\bf a}_{1}t,{\bf b}_{0})\cup ({\bf a}_{0},{\bf b}_{0}s+{\bf b}_{1}t)$.
 
In general, $W_{(d_a,d_b)}$ consists of $(d_a+1)(d_b+1)+1$ rows labeled by 
$z_{i}\; (i=0,1,\cdots, d_{a}),\;w_{j}\;(j=0,1,\cdots, d_{b})$, $f_{(i,j)} (i=1,\cdots, d_a,\;j=1,\cdots, d_b)$ and 
$(d_a+2)(d_b+2)-3$ columns labeled by ${\bf a}_{j},\; {\bf b}_{j},\; u_{(i_a,i_b)}$. 
Elements of the matrix $W_{(d_a,d_b)}$ are described as follows.
\begin{description}
\item{column ${\bf a}_i\;\;(0\leq i\leq d_a)$:} $z_i$ element is $1$ and the other elements are $0$. 
\item{column ${\bf b}_i\;\;(0\leq i\leq d_b)$:} $w_i$ element is $1$ and the other elements are $0$. 
\item{column $u_{(i,j)}\;\;(1\leq i\leq d_a-1,\;1\leq j\leq d_b-1)$:} $f_{(i+1,j)}$ and $f_{(i,j+1)}$ elements are $1$,
$f_{(i,j)}$ and $f_{(i+1,j+1)}$ elements are $-1$ and the other elements are $0$. 
\item{column $u_{(i,0)}\;\;(1\leq i\leq d_a-1)$:} $z_{i-1}$ element is $-1$, $z_{i}$ element is $1$, $f_{(i,1)}$ element is $1$,
$f_{(i+1,1)}$ element is $-1$ and the other elements are $0$.
\item{column $u_{(i,d_b)}\;\;(1\leq i\leq d_a-1)$:} $z_{i}$ element is $1$, $z_{i+1}$ element is $-1$, $f_{(i+1,d_b)}$ element is $1$,
$f_{(i,d_b)}$ element is $-1$ and the other elements are $0$.
\item{column $u_{(0,j)}\;\;(1\leq j\leq d_b-1)$:} $w_{j-1}$ element is $-1$, $w_{j}$ element is $1$, $f_{(1,j)}$ element is $1$,
$f_{(1,j+1)}$ element is $-1$ and the other elements are $0$.
\item{column $u_{(d_a,j)}\;\;(1\leq j\leq d_b-1)$:} $w_{j}$ element is $1$, $w_{j+1}$ element is $-1$, $f_{(d_a,j+1)}$ element is $1$,
$f_{(d_a,j)}$ element is $-1$ and the other elements are $0$.
\item{column $u_{(0,d_b)}$:} $z_{1}$ and $w_{d_b-1}$ elements are $-1$, $f_{(1,d_b)}$ element is $1$
and the other elements are $0$.
\item{column $u_{(d_a,0)}$:} $z_{d_a-1}$ and $w_{1}$ elements are $-1$, $f_{(d_a,1)}$ element is $1$
and the other elements are $0$.
\end{description}
As an example, we write down $W_{(2,1)}$, $W_{(3,1)}$ and $W_{(2,2)}$  below :
\begin{eqnarray}
W_{(2,1)}=\bordermatrix{
&{\bf a}_{0}&{\bf a}_{1}&{\bf a}_{2}&{\bf b}_{0}&{\bf b}_{1}&u_{(0,1)}&u_{(1,1)}&u_{(1,0)}&u_{(2,0)}\cr
        z_{0}&1&0&0&0&0&0&0&-1&0\cr
        z_{1}&0&1&0&0&0&-1&1&1&-1\cr
        z_{2}&0&0&1&0&0&0&-1&0&0\cr
        w_{0}&0&0&0&1&0&-1&0&0&0\cr
        w_{1}&0&0&0&0&1&0&0&0&-1\cr
        f_{(1,1)}&0&0&0&0&0&1&-1&1&0\cr
        f_{(2,1)}&0&0&0&0&0&0&1&-1&1\cr},
\end{eqnarray}
\begin{eqnarray}
W_{(3,1)}=\bordermatrix{
&{\bf a}_{0}&{\bf a}_{1}&{\bf a}_{2}&{\bf a}_{3}&{\bf b}_{0}&{\bf b}_{1}&u_{(0,1)}&u_{(1,1)}&u_{(2,1)}&u_{(1,0)}&u_{(2,0)}&u_{(3,0)}\cr
        z_{0}&1&0&0&0&0&0& 0&0&0&-1&0&0\cr
        z_{1}&0&1&0&0&0&0& -1&1&0&1&-1&0\cr
        z_{2}&0&0&1&0&0&0& 0&-1&1&0&1&-1\cr
        z_{3}&0&0&0&1&0&0& 0&0&-1&0&0&0\cr
        w_{0}&0&0&0&0&1&0& -1&0&0&0&0&0\cr
        w_{1}&0&0&0&0&0&1& 0&0&0&0&0&-1\cr
        f_{(1,1)}&0&0&0&0&0&0& 1&-1&0&1&0&0\cr
        f_{(2,1)}&0&0&0&0&0&0& 0&1&-1&-1&1&0\cr
        f_{(3,1)}&0&0&0&0&0&0& 0&0&1&0&-1&1\cr},
\end{eqnarray}
\begin{eqnarray}
W_{(2,2)}=\bordermatrix{
&{\bf a}_{0}&{\bf a}_{1}&{\bf a}_{2}&{\bf b}_{0}&{\bf b}_{1}&{\bf b}_{2}&u_{(1,0)}&u_{(2,0)}&u_{(0,1)}&u_{(1,1)}&u_{(2,1)}&u_{(0,2)}&u_{(1,2)}\cr
        z_{0}&1&0&0&0&0&0& -1&0&0&0&0&0&0\cr
        z_{1}&0&1&0&0&0&0& 1&-1&0&0&0&-1&1\cr
        z_{2}&0&0&1&0&0&0& 0&0&0&0&0&0&-1\cr
        w_{0}&0&0&0&1&0&0& 0&0&-1&0&0&0&0\cr
        w_{1}&0&0&0&0&1&0& 0&-1&1&0&1&-1&0\cr
        w_{2}&0&0&0&0&0&1& 0&0&0&0&-1&0&0\cr
        f_{(1,1)}&0&0&0&0&0&0& 1&0&1&-1&0&0&0\cr
        f_{(2,1)}&0&0&0&0&0&0& -1&1&0&1&-1&0&0\cr  
        f_{(1,2)}&0&0&0&0&0&0& 0&0&-1&1&0&1&-1\cr
        f_{(2,2)}&0&0&0&0&0&0& 0&0&0&-1&1&0&1\cr
}.
\end{eqnarray}
To understand the rule for determining elements of these matrices, it is convenient to write degree diagrams presented in Fig. 1 and Fig. 2. 
The degree diagram of type $(d_a,d_b)$ consists of vertices $(i,j)\;(0\leq i\leq d_a,\; 0\leq j\leq d_b)$ ordered in a rectangular shape 
with arrows from $(i,j)$ to $(i-1,j)$ and to $(i,j-1)$. The symbol $f_{(i,j)}$ is located at the center of the block surrounded by 
the vertices $(i-1,j-1)$, $(i,j-1)$, $(i-1,j)$ and $(i,j)$. The vertex $(i,j)$ corresponds to the coordinate $u_{(i,j)}$ if 
$(i,j)\neq (0,0), (d_a,d_b)$. The complicated rule of the description of $z_{*}, w_{*}$ element of the column $u_{(i,j)}$ arises from 
whether the vertex $(i,j)$ is located in the interior, or on the edge, or on the apex of the big rectangle whose four corner vertices are
given by, 
\begin{eqnarray}
(0,0),\;(d_a,0),\;(0,d_b),\;(d_a,d_b).
\label{4vertices}
\end{eqnarray}
 We can give 
graphical explanation of the description of $f_{(*,*)}$ element of the column $u_{(i,j)}$ with the diagram. If the vertex $(i,j)$ is 
located at the upper-left or lower-right corner of one of the blocks with $f_{(k,l)}$ at its center, the $f_{(k.l)}$ element of the column $u_{(i,j)}$ is $1$.
If $(i,j)$ is located at upper-right or lower-left corner of one of the $f_{(k,l)}$ blocks, the $f_{(k.l)}$ element of the column $u_{(i,j)}$ is $-1$. 
Otherwise, the $f_{(k.l)}$ element of the column $u_{(i,j)}$ is $0$. 
\begin{figure}[h]
      \epsfxsize=4cm
     \centerline{\epsfbox{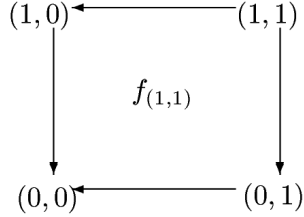}}
    \caption{\bf Degree Diagram of Type $(1,1)$}
\label{11}
\end{figure}
\begin{figure}[h]
      \epsfxsize=16cm
     \centerline{\epsfbox{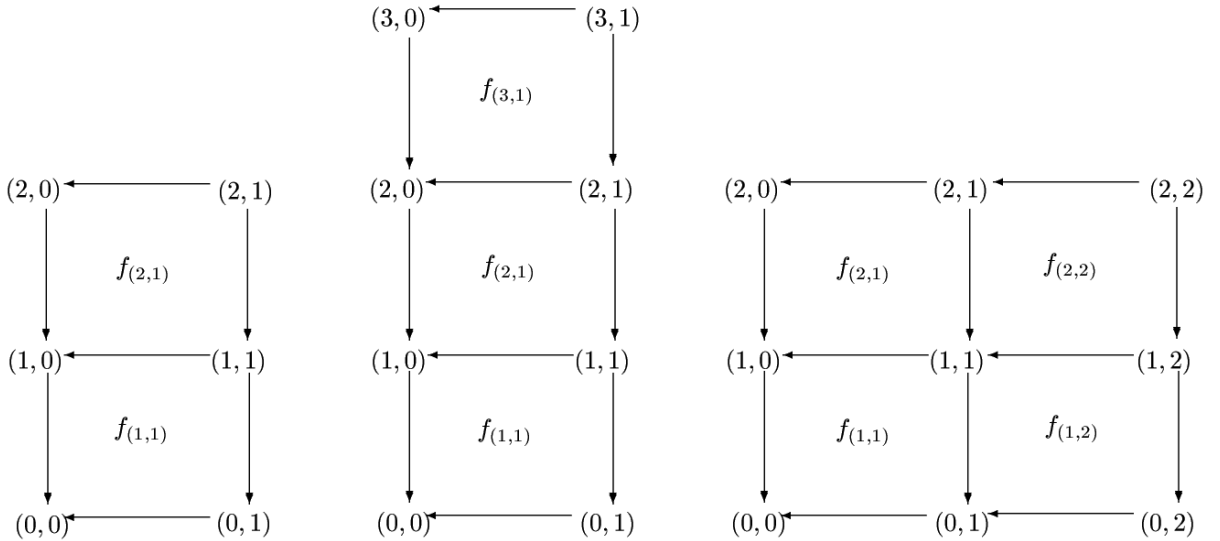}}
    \caption{\bf Degree Diagrams of Type $(2,1)$, $(3,1)$ and $(2,2)$}
\label{diags}
\end{figure}
With this setup, let us explain how the locus $u_{(i_a,i_b)}=0$ describes chains of two polynomial maps:
\begin{eqnarray}
( \sum_{j=0}^{i_{a}}{\bf a}_{j}s_{1}^{i_{a}-j}t_{1}^{j} , \sum_{j=0}^{i_{b}}{\bf b}_{j}s_{1}^{i_{b}-j}t_{1}^{j})
\cup ( \sum_{j=0}^{d_{a}-i_{a}}{\bf a}_{j}s_{2}^{d_{a}-i_{a}-j}t_{2}^{j} , \sum_{j=0}^{d_{b}-i_{b}}{\bf b}_{j}s_{2}^{d_{b}-i_{b}-j}t_{2}^{j}).
\label{chain2}
\end{eqnarray}
From the conditions:
$$ ({\bf a}_{i_a},\prod_{k=0}^{d_b}u_{(i_a,k)})\neq {\bf 0},  ({\bf b}_{i_b},\prod_{k=0}^{d_a}u_{(k,i_b)})\neq {\bf 0},$$
we can see that $u_{(i_a,i_b)}=0$ implies ${\bf a}_{i_a},{\bf b}_{i_b}\neq {\bf 0}$. The last condition in (\ref{mpf0}) tells 
us that $u_{(i_a,i_b)}=0$ also implies $u_{(k,l)}\neq 0$ if $(k,l)$ is no bigger or no smaller than $(i_a, i_b)$. Therefore, 
we can trivialize these coordinates by using the torus action $f_{(*,*)}$ whose block is the upper-left or lower-right of the 
vertex $(i_a,i_b)$. After this operation, we can define new coordinates $\tilde{u}_{(i,j)}\;((i,j)<(i_a,i_b)\;\mbox{or}\;
(i,j)>(i_a,i_b))$ as follows:
\begin{eqnarray}
\tilde{u}_{(i,j)}=\left\{\begin{array}{cc}u_{(i,j)} & (i\neq i_a\; \mbox{and}\;j\neq i_b), \\
                                          \prod_{k=i_a}^{d_a}u_{(k,j)}& (i= i_a\; \mbox{and}\;j<i_b), \\
                                          \prod_{k=0}^{i_a}u_{(k,j)}& (i= i_a\; \mbox{and}\;j>i_b), \\
                                          \prod_{l=i_b}^{d_b}u_{(i,l)}& (i<i_a\; \mbox{and}\;j=i_b), \\
                                          \prod_{l=0}^{i_b}u_{(i,l)}& (i>i_a\; \mbox{and}\;j=i_b). \\
\end{array}
\right. 
\end{eqnarray}
If we write down the corresponding weight matrix with columns labeled by ${\bf a}_i$, ${\bf b}_j$ and $\tilde{u}_{(i,j)}$ and with 
rows labeled by $z_i$, $w_j$ and $f_{(k.l)}\;((k,l)\leq (i_a,i_b)\;\mbox{or}\; (k,l)>(i_a,i_b))$, we observe that the 
locus $u_{(i,j)}=0$ describes the chains of two polynomial maps in (\ref{chain2}). Let us take the case when $(d_a, d_b)=(2,2)$ 
for example. If $u_{(2,0)}=0$, we introduce the new coordinates $\tilde{u}_{(1,0)}=u_{(1,0)}u_{(1,1)}u_{(1,2)} $ and 
$\tilde{u}_{(2,1)}=u_{(0,1)}u_{(1,1)}u_{(2,1)}$. Then the weight matrix associated with the locus is given as follows:
 \begin{eqnarray}
\bordermatrix{
&{\bf a}_{0}&{\bf a}_{1}&{\bf a}_{2}&{\bf b}_{0}&{\bf b}_{1}&{\bf b}_{2}&\tilde{u}_{(1,0)}&\tilde{u}_{(2,1)}\cr
        z_{0}&1&0&0&0&0&0& -1&0\cr
        z_{1}&0&1&0&0&0&0& 2&0\cr
        z_{2}&0&0&1&0&0&0& -1&0\cr
        w_{0}&0&0&0&1&0&0& 0&-1\cr
        w_{1}&0&0&0&0&1&0& 0&2\cr
        w_{2}&0&0&0&0&0&1& 0&-1\cr
       }, 
\end{eqnarray}
where the column of $\tilde{u}_{(1,0))}$ (resp. $\tilde{u}_{(2,1)}$) is obtained by adding up $ u_{(1,0)}$, $u_{(1,1)}$ and $u_{(1,2)}$ (resp. 
$ u_{(0,1)}$, $u_{(1,1)}$ and $u_{(2,1)}$) column vectors of $W_{(2,2)}$ and by eliminating unnecessary elements.    
From this matrix, we can easily see that the corresponding locus describes chains of two polynomial maps of 
degree $(2,0)$ and of degree $(0,2)$. If $u_{(0,1)}=0$, the new coordinates are given as follows:
\begin{eqnarray}
\tilde{u}_{(0,2)}=u_{(0,2)},\;\tilde{u}_{(1,2)}=u_{(1,2)}, \;\tilde{u}_{(1,1)}=u_{(1,0)}u_{(1,1)},\;\tilde{u}_{(2,1)}=u_{(2,0)}u_{(2,1)},
\end{eqnarray}
and the corresponding weight matrix becomes,
\begin{eqnarray}
\bordermatrix{
&{\bf a}_{0}&{\bf a}_{1}&{\bf a}_{2}&{\bf b}_{0}&{\bf b}_{1}&{\bf b}_{2}& \tilde{u}_{(0,2)}&\tilde{u}_{(1,2)}&\tilde{u}_{(1,1)}&\tilde{u}_{(2,1)}\cr
        z_{0}&1&0&0&0&0&0& 0&0&-1&0\cr
        z_{1}&0&1&0&0&0&0& -1&1&1&-1\cr
        z_{2}&0&0&1&0&0&0& 0&-1&0&0\cr
        w_{0}&0&0&0&1&0&0& 0&0&0&0\cr
        w_{1}&0&0&0&0&1&0& -1&0&0&0\cr
        w_{2}&0&0&0&0&0&1& 0&0&0&-1\cr
    f_{(1,2)}&0&0&0&0&0&0& 1&-1&1&0\cr
    f_{(2,2)}&0&0&0&0&0&0& 0&1&-1&1\cr}.
\end{eqnarray}
This matrix includes a copy of $W_{(2,1)}$. Hence it describes chains of two polynomial maps of 
degree $(0,1)$ and $(2,1)$. In this way, we can observe that the locus $u_{(i_a, i_b)}=0$ corresponds to 
chains of two polynomial maps in (\ref{chain2}). In the same way as the $CP^{N-1}$ case, we can consider multi-zero locus: 
\begin{eqnarray}
&&u_{(d_{a,1},d_{b,1})}=u_{(d_{a,2},d_{b,2})}=\cdots=u_{(d_{a,l-1},d_{b,l-1})}=0,\no\\
&&((0,0)=(d_{a,0},d_{b,0})<(d_{a,1},d_{b,1})<(d_{a,2},d_{a,2})<\cdots
<(d_{a,l-1},d_{a,l-1})<(d_{a,l},d_{a,l})=(d_a,d_b)). 
\label{seq1}
\end{eqnarray}
This locus corresponds to chains of polynomial maps:
\begin{eqnarray}
&&\cup_{j=1}^{l}\bigl( \sum_{m_{a,j}=0}^{d_{a,j}-d_{a,j-1}} {\bf a}_{d_{a,j-1}+m_{a,j}}(s_{j})^{m_{a,j}}(t_{j})^{d_{a,j}-d_{a,j-1}-m_{a,j}}
, \sum_{m_{b,j}=0}^{d_{b,j}-d_{b,j-1}} {\bf b}_{d_{b,j-1}+m_{b,j}}(s_{j})^{m_{b,j}}(t_{j})^{d_{b,j}-d_{b,j-1}-m_{b,j}} \bigr),\no\\
&&\;\;\bigl( {\bf a}_{d_{a,j}} , {\bf b}_{d_{b,j}} \neq {\bf 0},\;\;j=0,1,\cdots,l\bigr).
\label{chain3}
\end{eqnarray}

\subsubsection{Localization Computation}
We have constructed the moduli space of polynomial maps of degree $\dd$ with two marked points,
$\widetilde{Mp}_{0,2}(F_{0},(\dd))$.
Next, we define and compute an analogue of the genus $0$ local Gromov-Witten invariant of $K_{F_{0}}$ defined by,
\begin{eqnarray}
\langle{\cal O}_{\alpha}{\cal O}_{\beta}\rangle_{0,\dd}:=
\int_{ \overline{M}_{0,2}(F_{0},\dd) }ev_{1}^{*}(\alpha)\wedge ev_{2}^{*}(\beta)\wedge c_{top}(R^{1}\pi_{*}ev_{3}^{*}
({\cal O}_{F_{0}}(-2a-2b))),
\label{lgw1} 
\end{eqnarray} 
by changing the moduli space of stable maps $ \overline{M}_{0,2}(F_{0},\dd) $ into $\widetilde{Mp}_{0,2}(F_{0},\dd)$.
In (\ref{lgw1}), $$ev_{i}: \overline{M}_{0,n}(F_{0},\dd) \rightarrow F_{0}$$ is the evaluation map at the $i$-th marked point 
of stable curves, and $\pi$ is the forgetful map that forgets the third marked point of $\overline{M}_{0,3}(F_{0},(d_a,d_b))$. To 
construct an analogue of $\langle{\cal O}_{\alpha}{\cal O}_{\beta}\rangle_{0,(d_a,d_b)}$, which should be given as an intersection 
number on $ \widetilde{Mp}_{0,2}(F_{0},\dd) $, we have to define cohomology classes which correspond to $ev_{1}^{*}(\alpha)$, 
$ev_{2}^{*}(\beta)$ and $c_{top}(R^{1}\pi_{*}ev_{3}^{*}({\cal O}_{F_{0}}(-2a-2b)))$ respectively.   
For the first two classes, our task is easily accomplished because we have evaluation maps $ev_{1}$ and $ev_{2}$ defined 
on  $\widetilde{Mp}_{0,2}(F_{0},\dd)$:
\begin{eqnarray}
&& ev_{1}( [({\bf a}_{0},\cdots,{\bf a}_{d_a},{\bf b}_{0},\cdots,{\bf b}_{d_b},u_{(1,0)},\cdots,u_{(d_a-1,d_b)})] )
=[({\bf a}_{0},{\bf b}_{0})]\in F_{0},\no\\
&& ev_{2}([({\bf a}_{0},\cdots,{\bf a}_{d_a},{\bf b}_{0},\cdots,{\bf b}_{d_b},u_{(1,0)},\cdots,u_{(d_a-1,d_b)})])
=[({\bf a}_{d_a},{\bf b}_{d_b})]\in F_{0},
\label{ev}
\end{eqnarray}
where $[*]$ represents equivalence class of torus actions. 
Let us turn to an analogue of $c_{top}(R^{1}\pi_{*}ev_{3}^{*}({\cal O}_{F_{0}}(-2a-2b)))$. If we look back 
at the discussion in Subsection 2.3, we can define a rank $2d_a+2d_b-1$ orbi-bundle ${\cal E}_{\dd}$ on $Mp_{0,2}(F_0,\dd)$ by using
Kodaira-Serre duality,
\begin{eqnarray}
H^{1}(CP^1,\varphi^{*}{\cal O}_{F_{0}}(-2a-2b))\simeq ( H^{0}(CP^1,\varphi^{*}{\cal O}_{F_{0}}(2a+2b)\otimes K_{CP^1}))^{\vee},
\label{ks2}
\end{eqnarray} 
where $\varphi$ is a polynomial map:$(\sum_{j=0}^{d_a}{\bf a}_{j}s^{d_a-j}t^{j},\sum_{j=0}^{d_b}{\bf b}_js^{d_b-j}t^{j})$.
We can extend this orbi-bundle to the whole $ \widetilde{Mp}_{0,2}(F_{0},\dd) $ by generalizing the exact sequence (\ref{exseq2}).
In this way, we can define an analogue 
of $\langle{\cal O}_{\alpha}{\cal O}_{\beta}\rangle_{0,\dd}$ as an intersection 
number of $ \widetilde{Mp}_{0,2}(F_{0},\dd) $:
\begin{eqnarray}
w({\cal O}_{\alpha}{\cal O}_{\beta})_{0,\dd}:=
\int_{ \widetilde{Mp}_{0,2}(F_{0},\dd) }ev_{1}^{*}(\alpha)\wedge ev_{2}^{*}(\beta)\wedge c_{top}({\cal E}_{\dd}).
\label{lgw2} 
\end{eqnarray}
From now on, we compute the intersection number $w({\cal O}_{\alpha}{\cal O}_{\beta})_{0,\dd}$ by using the localization 
computation. To apply this technique, we introduce a torus action flow to $ \widetilde{Mp}_{0,2}(F_{0},\dd) $ as follows: 
\begin{eqnarray}
[( e^{\lambda_{0}t}{\bf a}_{0},  e^{\lambda_{1}t}{\bf a}_{1},\cdots,  e^{\lambda_{d_a}t}{\bf a}_{d_{a}},
 e^{\mu_{0}t}{\bf b}_{0},  e^{\mu_{1}t}{\bf b}_{1},\cdots,  e^{\mu_{d_b}t}{\bf b}_{d_{b}}, u_{(1,0)} ,\cdots,u_{(d_a-1,d_b)})],
\label{flow2}
\end{eqnarray}  
where $\lambda_{i}$ and $\mu_{j}$ are characters of the torus action. We take these characters as generic as possible. 
We then have to determine the fixed point set of $ \widetilde{Mp}_{0,2}(F_{0},\dd) $ under the flow. Let us consider 
the case when all the $u_{(i,j)}$'s are non-zero. In this case, we can set these $u_{(i,j)}$'s to $1$ by using the $C^{\times}$
actions in the definition of  $ \widetilde{Mp}_{0,2}(F_{0},\dd) $ and represent a point in this locus as a single 
polynomial map:
\begin{eqnarray}
( \sum_{i=0}^{d_{a}}{\bf a}_{i}s^{d_{a}-i}t^{i} , \sum_{j=0}^{d_{b}}{\bf b}_{j}s^{d_{b}-j}t^{j}).
\end{eqnarray}  　
Looking back at (\ref{flow2}), we can see that fixed points do exist when $\dd=(d_a,0)\;\mbox{or}\;(0,d_b)$. In this case, fixed points are given by 
polynomial maps:
\begin{eqnarray}
({\bf a}_{0}s^{d_a}+{\bf a}_{d_a}t^{d_a},{\bf b}_{0}),\;\; ({\bf a}_{0},{\bf b}_{0}s^{d_b}+{\bf b}_{d_b}t^{d_b}),
\end{eqnarray} 
because the torus action flow given by  (\ref{flow2}) is canceled by the three remaining $C^{\times}$ actions used in the definition of the moduli space.
But if $d_a,d_b>0$, we can conclude that there are no fixed points in this locus. Naively, we might say that the map:
\begin{eqnarray}
({\bf a}_{0}s^{d_a}+{\bf a}_{d_a}t^{d_a},{\bf b}_{0}s^{d_b}+{\bf b}_{d_b}t^{d_b}),
\label{non}
\end{eqnarray} 
is a candidate; however, four independent characters $\lambda_{0}$, $\lambda_{d_a}$, $\mu_{0}$ and $\mu_{d_b}$ act on it. These cannot be 
canceled by the remaining three $C^{\times}$ actions. Therefore, the points represented by (\ref{non}) do "move" under 
the flow (\ref{flow2}).

Next, we consider the locus where we can pick up the sequence of bi-degrees (\ref{seq1}) and represent a point by the chain of polynomial maps (\ref{chain3}). 
From the previous discussion, we conclude that there exist non-trivial fixed points if and only if 
\begin{eqnarray}
d_{a,j}-d_{a,j-1}=0,\;\mbox{or}\;d_{b,j}-d_{b,j-1}=0\; \mbox{for all}\;j=1,2,\cdots,l. 
\label{fc}
\end{eqnarray}
If the above condition is satisfied, fixed points can be represented by chains of polynomial maps whose $j$-th component is given by, 
\begin{eqnarray}
({\bf a}_{d_{a,j-1}},{\bf b}_{d_{b,j-1}}(s_{j})^{d_{b,j}-d_{b,j-1}}+{\bf b}_{d_{b,j}}(t_{j})^{d_{b,j}-d_{b,j-1}})\;\mbox{or}\;
({\bf a}_{d_{a,j-1}}(s_{j})^{d_{a,j}-d_{a,j-1}}+{\bf a}_{d_{a,j}}(t_{j})^{d_{a,j}-d_{a,j-1}},{\bf b}_{d_{b,j-1}}),
\label{chain5}
\end{eqnarray}  
respectively. In this way, we have seen that fixed points are classified by the sequence of bi-degrees satisfying (\ref{fc}).
We introduce here a set of ordered partitions of bi-degree ${\dd}$:  
\begin{equation}
OP_{\dd}=\{\sigma_{\dd}=(\dd_{1},\dd_{2},\cdots,\dd_{l(\sigma_{\dd})})\;\;|\;\;
\sum_{j=1}^{l(\sigma_{d})}\dd_{j}=\dd\;\;,\;\;\dd_{j}=(d_{a,j},0)\;\mbox{or}\;\dd_{j}=(0,d_{b,j})\},
\label{part2} 
\end{equation}
whose element is in one-to-one correspondence with a sequence of bi-degrees satisfying (\ref{fc}).
We also introduce the notation:
\begin{eqnarray}
|\dd_{j}|:=\left\{\begin{array}{cc} d_{a,j}&\mbox{if}\;\dd_j=(d_{a,j},0),\\
d_{b,j}&\mbox{if}\;\dd_j=(0,d_{b,j}).\end{array}\right.
\end{eqnarray}
Let $F_{\sigma_{\dd}}$ be a connected component of the fixed point set labeled by $\sigma_{d}\in OP_{\dd}$. By relabeling subscripts, 
it consists of chains of polynomial maps of length $l(\sigma_{\dd})$ whose $j$-th component is given by, 
\begin{eqnarray}
( {\bf a}_{j-1}(s_{j})^{|\dd_{j}|}+{\bf a}_{j}(t_{j})^{|\dd_{j}|} ,{\bf b}_{j-1})\; \mbox{or} \;({\bf a}_{j-1},{\bf b}_{j-1}(s_{j})^{|\dd_{j}|}+{\bf b}_{j}(t_{j})^{|\dd_{j}|}),
\label{ch6}
\end{eqnarray} 
respectively if $\dd_{j}=(d_{a,j},0)$ or $\dd_{j}=(0,d_{b,j})$. Therefore, it is set-theoretically given by a subset of,
\begin{eqnarray}
(F_{0})_{0} \times(F_{0})_{1} \times(F_{0})_{2} \times\cdots \times(F_{0})_{l(\sigma_{\dd})},\; ((F_{0})_j=\{[({\bf a}_j,{\bf b}_j)]\}),
\label{pro1}
\end{eqnarray} 
defined by the following conditions:
\begin{eqnarray}
&&{\bf b}_{j-1}={\bf b}_{j}\;\;\;\mbox{if}\; \dd_{j}=(d_{a,j},0),\no\\
&&{\bf a}_{j-1}={\bf a}_{j} \;\;\;\mbox{if}\; \dd_{j}=(0,d_{b,j}).
\label{id}
\end{eqnarray}
We have to note one subtlety here. Though $F_{\sigma_{\dd}}$ is set-theoretically bijective to the space given in (\ref{pro1}), it should be considered as
an orbifold on which an abelian group $\oplus_{j=1}^{l(\sigma_{\dd})}\bigl({\bf Z}/(|\dd_{j}|{\bf Z})\bigr)$ acts. This group action comes from the $C^{\times}$ actions 
in the definition of $\widetilde{Mp}_{0,0}(F_{0},\dd)$ that keep the chains of polynomial maps in this component fixed.

We now describe normal bundle of $F_{\sigma_{\dd}}$ in  $\widetilde{Mp}_{0,0}(F_{0},\dd)$. As was discussed in our previous paper \cite{vs}, it has two degrees of freedom:
\begin{itemize}
\item[(i)] Deformations of each component of the chain of polynomial maps in $\widetilde{Mp}_{0,0}(F_{0},\dd)$.
\item[(ii)] Resolutions of nodal singularities of the image curve in $F_{0}$.
\end{itemize}
These can be easily realized as sheaves of the orbifold $F_{\sigma_{\dd}}$ by a  straightforward generalization of the discussion in \cite{vs}  to this case.
Let us introduce the notation:
\begin{eqnarray}
{\cal O}_{F_{0}}(\frac{m}{\dd_{j}}):=\left\{\begin{array}{cc}  {\cal O}_{F_{0}}(\frac{m}{|\dd_{j}|}a) &\;\mbox{if}\;\dd_j=(d_{a,j},0),\\
  {\cal O}_{F_{0}}(\frac{m}{|\dd_{j}|}b) &\;\mbox{if}\;\dd_j=(0,d_{b,j}).\end{array}\right.
\end{eqnarray}  
With this notation, we can write down the normal bundle as follows:    
\begin{eqnarray}
&& \mathop{\oplus}_{j=1}^{l(\sigma_{\dd})}\biggl(\mathop{\oplus}_{i=1}^{|\dd_{j}|-1} \bigl( 
{\cal O}_{(F_{0})_{j-1}}(\frac{i}{\dd_{j}}) \otimes {\cal O}_{(F_{0})_{j}}(\frac{|\dd_j|-i}{\dd_{j}})
\bigr)^{\oplus 2} \biggr)\oplus \no\\
&&\oplus \mathop{\oplus}_{j=1}^{l(\sigma_{\dd})-1}\bigl( {\cal O}_{(F_{0})_{j-1}}(-\frac{1}{\dd_{j}}) \otimes {\cal O}_{(F_{0})_{j}}(\frac{1}{\dd_{j}})
\otimes {\cal O}_{(F_{0})_{j}}(\frac{1}{\dd_{j+1}}) \otimes {\cal O}_{(F_{0})_{j+1}}(-\frac{1}{\dd_{j+1}}) \bigr) ,
\label{no}
\end{eqnarray}
where the first line (resp. the second line) corresponds to the degree of freedom (i) (resp. (ii)).

We have described the fixed point set of the torus action flow and the normal bundle of its connected components. What remains is to describe 
is restriction of the orbi-bundle ${\cal E}_{\dd}$ to $F_{\sigma_{\dd}}$. This task can also be accomplished by the direct 
generalization of the discussion in \cite{vs}. The result turns out to be,
\begin{eqnarray}
&& \mathop{\oplus}_{j=1}^{l(\sigma_{\dd})}\biggl(\mathop{\oplus}_{i=1}^{2|\dd_{j}|-1} \bigl( 
{\cal O}_{(F_{0})_{j-1}}(\frac{i}{\dd_{j}}) \otimes {\cal O}_{(F_{0})_{j}}(\frac{-i}{\dd_{j}})
\otimes{\cal O}_{(F_{0})_{j-1}}(-2a-2b)\bigr)\biggr)\oplus \no\\
&&\oplus \mathop{\oplus}_{j=1}^{l(\sigma_{\dd})-1}\bigl( {\cal O}_{(F_{0})_{j}}(-2a-2b)\bigr) .
\label{cano}
\end{eqnarray}
The first line of (\ref{cano}) comes from $H^{1}(CP^1,\varphi_{j}^{*}{\cal O}_{F_{0}}(-2a-2b))$ where 
$\varphi_{j}:CP^1\rightarrow F_{0}$ is the j-th component map of chains of polynomial maps in (\ref{ch6}).
The second line comes from effects of nodal singularities of the image curve.

Now, we are ready to apply the localization theorem to $ w({\cal O}_{\alpha}{\cal O}_{\beta})_{0,\dd} $ given in (\ref{lgw2}).
In the same way as was used in \cite{vs}, we take the non-equivariant limit $\lambda_{j},\mu_{j}\rightarrow 0$, with which 
we still can obtain well-defined results. To describe the results of the localization computation, we introduce the 
notation:
\begin{eqnarray}
&& z_j:=c_{1}({\cal O}_{(F_{0})_{j}}(a)), \; w_j:=c_{1}({\cal O}_{(F_{0})_j}(b)).
\label{chern}
\end{eqnarray}   
Since we have expressions for the normal bundle and ${\cal E}_{\dd}|_{F_{\sigma_{\dd}}}$ as sheaves on $F_{\sigma_{\dd}}$, 
it is straightforward to write down the formula we need. As the first step, we define the following 
rational function to express contributions from the first lines of (\ref{no}) and (\ref{cano}) 
\begin{eqnarray}
G(\dd;z_{0},z_{1},w_{0},w_1):=\left\{\begin{array}{cc}
\frac{\prod_{j=1}^{2|\dd|-1}\bigl(\frac{-jz_0-(2|\dd|-j)z_1-2w_{0}}{|\dd|}\bigr)}
{\prod_{j=1}^{|\dd|-1}\bigl(\frac{jz_{0}+(|\dd|-j)z_{1}}{|\dd|}\bigr)^{2}},&\;\;\mbox{if}\; \dd=(d,0) ,\\
\frac{\prod_{j=1}^{2|\dd|-1}\bigl(\frac{-jw_0-(2|\dd|-j)w_1-2z_{0}}{|\dd|}\bigr)}
{\prod_{j=1}^{|\dd|-1}\bigl(\frac{jw_{0}+(|\dd|-j)w_{1}}{|\dd|}\bigr)^{2}},&\;\;\mbox{if}\; \dd=(0,d) .
\end{array}\right.
\label{G}
\end{eqnarray}
To express contributions from the second lines of (\ref{no}) and (\ref{cano}), we introduce another rational 
function:
\begin{eqnarray}
H(\dd_1;\dd_2,z_0,z_1,z_2,w_0,w_1,w_2):=\left\{\begin{array}{cc}
\frac{(-2z_1-2w_1)}{ \bigl(\frac{z_1-z_0}{|\dd_1|}+ 
\frac{z_1-z_2}{|\dd_2|}\bigr)},&\;\mbox{if}\;\dd_{1}=(d_{a,1},0)\;\mbox{and}\;\dd_{2}=(d_{a,2},0),\\
\frac{(-2z_1-2w_1)}{ \bigl(\frac{z_1-z_0}{|\dd_1|}+ 
\frac{w_1-w_2}{|\dd_2|}\bigr)},&\;\mbox{if}\;\dd_{1}=(d_{a,1},0)\;\mbox{and}\;\dd_{2}=(0,d_{b,2}), \\
\frac{(-2z_1-2w_1)}{ \bigl(\frac{w_1-w_0}{|\dd_1|}+ 
\frac{z_1-z_2}{|\dd_2|}\bigr)},&\;\mbox{if}\;\dd_{1}=(0,d_{b,1})\;\mbox{and}\;\dd_{2}=(d_{a,2},0),\\
\frac{(-2z_1-2w_1)}{ \bigl(\frac{w_1-w_0}{|\dd_1|}+ 
\frac{w_1-w_2}{|\dd_2|}\bigr)},&\;\mbox{if}\;\dd_{1}=(0,d_{b,1})\;\mbox{and}\;\dd_{2}=(0,d_{b,2}).
\end{array}\right.
\label{Hf0}
\end{eqnarray}
With this setup, the contributions from ${\cal E}_{\dd}|_{F_{\sigma_{\dd}}}$ and from the normal bundle 
of $F_{\sigma_{\dd}}$ can be collected in the following integrand:
\begin{eqnarray}
K(\sigma_{\dd};z_{*},w_{*}):= \prod_{j=1}^{l(\sigma_{\dd})}G(\dd_j;z_{j-1},z_j,w_{j-1},w_j)
 \prod_{j=1}^{l(\sigma_{\dd})-1}H(\dd_j;\dd_{j+1},z_{j-1},z_j,z_{j+1},w_{j-1},w_j,w_{j+1}).
\label{K}
\end{eqnarray}
Next, we turn to the contributions from $ev_{1}^{*}(\alpha)$ and $ev_{2}^{*}(\beta)$. Since 
$\alpha,\beta\in H^{*}(F_{0},C)$, these can be written as $z^{s}w^{t}\; (s,t\in\{0,1\})$.
The definition of $ev_{i}$ in (\ref{ev}) ( we have to take care in relabeling subscripts ) directly 
leads us to, 
\begin{eqnarray}
ev_{1}^{*}(z^sw^t)=(z_{0})^s(w_{0})^t, \; ev_{2}^{*}(z^sw^t)=(z_{l(\sigma_{\dd})})^s(w_{l(\sigma_{\dd})})^t. 
\end{eqnarray} 
What remains to be done is to integrate out $ ev_{1}^{*}(\alpha)ev_{2}^{*}(\beta)K(\sigma_{\dd};z_{*},w_{*}) $
over $F_{\sigma_{\dd}}$. For this purpose, we note the following three facts: 
\begin{itemize}
\item[(i)] Integration of the cohomology element $\alpha\in H^{*}(F_{0},C)$ can be realized as the following residue 
integral in the variables $z$ and $w$:
\begin{eqnarray}
\int_{F_{0}}\alpha= \frac{1}{(2\pi\sqrt{-1})^2} \oint_{C_{0}} \frac{dz}{z^2} 
\oint_{C_{0}}\frac{dw}{w^2}\alpha.
\label{res0}
\end{eqnarray}
\item[(ii)] Looking back at (\ref{pro1}) and (\ref{id}), we make the identification:
\begin{eqnarray}
&&w_{j-1}= w_{j}\;\ ;\;\mbox{if}\; \dd_{j}=(d_{a,j},0),\no\\
&&z_{j-1}=z_{j} \;\;\;\mbox{if}\; \dd_{j}=(0,d_{b,j}).
\label{id2}
\end{eqnarray}
\item[(iii)]
$F_{\sigma_{\dd}}$ should be considered as an orbifold on which an abelian group 
$\oplus_{j=1}^{l(\sigma_{\dd})}\bigl({\bf Z}/(|\dd_{j}|{\bf Z})\bigr)$ acts.
\end{itemize}
Taking facts (i) and (ii) into account, we define the following operation on 
a rational function $f$ in $z_*$ and $w_*$:
\begin{eqnarray}
&& Res_{(F_{0})_j}(f):=\left\{\begin{array}{cc}
\bigl( \frac{1}{2\pi\sqrt{-1}}\oint_{C_{0}}\frac{dz_j}{(z_j)^2}f\bigr)|_{w_{j}=w_{j+1}}, 
&\;\mbox{if}\; \dd_{j+1}=(d_{a,j+1},0),\\
 \frac{1}{2\pi\sqrt{-1}}\oint_{C_{0}}\frac{dw_j}{(w_j)^2}(f|_{z_j=z_{j+1}}), 
&\;\mbox{if}\; \dd_{j+1}=(0,d_{b,j+1}),
\end{array}\right. (j=0,1,\cdots,l(\sigma_{\dd})-1),\no\\
&& Res_{(F_{0})_{l(\sigma_{\dd})}}(f):=\frac{1}{(2\pi\sqrt{-1})^2} \oint_{C_{0}} \frac{dz_{l(\sigma_{\dd})}}{(z_{l(\sigma_{\dd})})^2} 
\oint_{C_{0}}\frac{dw_{l(\sigma_{\dd})}}{(w_{l(\sigma_{\dd})})^2}f.
\label{resf}
\end{eqnarray}
With this definition and fact (iii) in mind, we conclude that the result of the integration is given by, 
\begin{eqnarray}
&& Amp(\sigma_{\dd};\alpha,\beta) =\biggl(\prod_{j=1}^{l(\sigma_{\dd})}\frac{1}{|\dd_{j}|}\biggr) Res_{(F_{0})_{l(\sigma_{\dd})}}( Res_{(F_{0})_{l(\sigma_{\dd})-1}}(
\cdots Res_{(F_{0})_{0}}( ev_{1}^{*}(\alpha)ev_{2}^{*}(\beta)K(\sigma_{\dd};z_{*},w_{*}) )\cdots)).\no\\
\label{amp}
\end{eqnarray}
Finally, the localization theorem tells us that,
\begin{eqnarray}
 w({\cal O}_{\alpha}{\cal O}_{\beta})_{0,\dd} =\sum_{\sigma_{\dd}\in OP_{\dd}}Amp(\sigma_{\dd};\alpha,\beta). 
\end{eqnarray}
\subsubsection{Numerical Results and the Mirror Computation}
In the previous sectionh, we obtained an explicit formula to compute $w({\cal O}_{\alpha}{\cal O}_{\beta})_{0,\dd}$. It is defined as an 
intersection number of $\widetilde{Mp}_{0,2}(F_{0},\dd)$ and has the same geometrical meaning as the local Gromov-Witten invariant 
$\langle{\cal O}_{\alpha}{\cal O}_{\beta}\rangle_{0,\dd}$ except that it is defined on the moduli space of polynomial maps 
instead of the moduli space of stable maps. These results then lead us naturally to the following questions. Is there any numerical difference between 
$w({\cal O}_{\alpha}{\cal O}_{\beta})_{0,\dd}$ and  $\langle{\cal O}_{\alpha}{\cal O}_{\beta}\rangle_{0,\dd}$? Can we compute 
$\langle{\cal O}_{\alpha}{\cal O}_{\beta}\rangle_{0,\dd}$ by using the data of $w({\cal O}_{\alpha}{\cal O}_{\beta})_{0,\dd}$?
In our previous paper \cite{vs}, we conjectured through explicit numerical computation that, in the $CP^{N-1}$ case, this new intersection 
number gives us the same information as the B-model used in the mirror computation. 
For example, $w({\cal O}_{1}{\cal O}_{h^{N-3+(N-k)d}})_{0,d}$ in (\ref{gw1}) reproduces the expansion coefficient of the mirror 
map in the $N\leq k$ case, regardless of whether the degree $k$ hypersurface in $CP^{N-1}$ is Calabi-Yau or of general type. Moreover, we 
can compute Gromov-Witten invariants of the hypersurface using the recipe of the standard mirror computation.
In the following, we demonstrate the mirror computation for $K_{F_{0}}$ by using the numerical data of $w({\cal O}_{\alpha}{\cal O}_{\beta})_{0,\dd}$
and argue that the same conjecture holds true in our current example. 

As the first step of mirror computation, we introduce the virtual classical intersection numbers used in our papers \cite{fj3}, \cite{fj4}:
\begin{eqnarray}
cl(z^3):=k,\;cl(z^2w):=-k,\;cl(zw^2)=k-\frac{1}{2},\;cl(w^2):=\frac{1}{2}-k,
\label{vcl}
\end{eqnarray}
where $k$ is a free parameter. If $z^sw^t$ is a monomial with $s+t\neq 3$, we set $cl(z^sw^t)=0$. Let $\eta_{\alpha\beta}$ and 
$C_{\alpha\beta\gamma}^{(0,0)}$ be symmetric tensors on the ${\bf C}$-vector space $H:=<1,z,w,z^2,zw,z^3>_{\bf C}$ defined by,
\begin{eqnarray}  
\eta_{\alpha\beta}:=cl(\alpha\beta),\; C_{\alpha\beta\gamma}^{(0,0)} :=cl(\alpha\beta\gamma).
\label{tens}
\end{eqnarray}
In (\ref{tens}), $\alpha$,$\beta$ and $\gamma$ take values in a basis of $H$ and should be considered as monomials in $z$ and $w$ 
in the r.h.s.. With these tensors, we can regard $H$ as the virtual classical intersection ring of $K_{F_{0}}$. 
As usual in the case of quantum cohomology ring, the relation $\eta_{\alpha\beta}=C_{1\alpha\beta}^{(0,0)}$ holds.   
For later use, 
we also define the symmetric tensor $\eta^{\alpha\beta}$ by the relation: $\eta_{\alpha\beta}\eta^{\beta\gamma}=\delta_{\alpha}^{\gamma}$. 
We present here $\eta_{\alpha\beta}$ and $\eta^{\alpha\beta}$ in matrix form:
\begin{eqnarray}  
(\eta_{\alpha\beta})=\bordermatrix{
&1&z& w& z^2&zw &z^3\cr
        1&0&0&0&0&0&k\cr
        z&0&0&0&k&-k&0\cr
        w&0&0&0&-k&k-\frac{1}{2}&0\cr
      z^2&0&k&-k&0&0&0\cr
       zw&0&-k&k-\frac{1}{2}&0&0&0\cr
      z^3&k&0&0&0&0&0\cr},
(\eta^{\alpha\beta})=\bordermatrix{
&1&z& w& z^2&zw &z^3\cr
        1&0&0&0&0&0&\frac{1}{k}\cr
        z&0&0&0&-\frac{2k-1}{k}&-2&0\cr
        w&0&0&0&-2&-2&0\cr
      z^2&0&-\frac{2k-1}{k}&-2&0&0&0\cr
       zw&0&-2&-2&0&0&0\cr
      z^3&\frac{1}{k}&0&0&0&0&0\cr}.
\label{etam}
\end{eqnarray}
Next, we give numerical results of the intersection number $ w({\cal O}_{\alpha}{\cal O}_{\beta})_{0,\dd} $
by using the generating function:
\begin{eqnarray}
w({\cal O}_{\alpha}{\cal O}_{\beta})_{0}:=C_{\alpha\beta z}^{(0,0)}x_{1}+C_{\alpha\beta w}^{(0,0)}x_{2}+
\sum_{\dd>(0,0)}w({\cal O}_{\alpha}{\cal O}_{\beta})_{0,\dd}e^{d_a x_1+d_b x_2}.
\label{gene1}    
\end{eqnarray}
Note that we add classical terms, defined through symmetric tensor in (\ref{tens}), to  $w({\cal O}_{\alpha}{\cal O}_{\beta})_{0}$.
In the following, we give numerical results for  $ w({\cal O}_{1}{\cal O}_{zw})_{0} $, $w({\cal O}_{z}{\cal O}_{z})_{0}$ and 
$w({\cal O}_{z}{\cal O}_{w})_{0}$ up to total degree $4$:
\begin{eqnarray}
w({\cal O}_{1}{\cal O}_{zw})_{0}&=&-kx_1+(k-\frac{1}{2})x_{2}-e^{x_{1}}-e^{x_{2}}
-\frac{3}{2}e^{2x_{1}}-6e^{x_{1}+x_{2}}-\frac{3}{2}e^{2x_{2}}- \no\\
&&-\frac {10}{3}e^{3x_{1}}- 30e^{2x_{1}+x_{2}}-30e^{x_{1}+2x_{2}} 
-\frac {10}{3}e^{3x_{2}}-\no \\
&&-\frac {35}{4}e^{4x_{1}}- 
140e^{3x_{1}+x_{2}} -315e^{2x_{1}+2x_{2}}-140e^{x_{1}+3x_{2}}-\frac {35}{4}e^{4x_{2}}-,\cdots\no\\ 
w({\cal O}_{z}{\cal O}_{z})_{0}&=&k{x_{1}}- k{x_{2}}
 - 2e^{x_{1}}- 5e^{2x_{1}} - 8e^{x_{1}+x_{2}} - \frac {44}{3}e^{3x_{1}}- 
76e^{2x_{1}+x_{2}}-32e^{x_{1}+2x_{2}}-\no\\
&&- \frac {93}{2} e^{4x_{1}}- 504e^{3x_{1}+x_{2}}
- 672e^{2x_{1}+2x_{2}}- 128e^{x_{1}+3x_{2}}-\cdots, \no\\
w({\cal O}_{z}{\cal O}_{w})_{0}&=&- k{x_{1}}+(k - \frac {1}{2}){x_{2}}
 - e^{x_{1}} - e^{x_{2}}
 - \frac {3}{2}e^{2x_{1}}- 
10e^{x_{1}+x_{2}} - \frac {3}{2}e^{2x_{2}}-\no\\
&&- \frac {10}{3} e^{3x_{1}}
 - 58e^{2x_{1}+x_{2}} - 58e^{x_{1}+2x_{2}}- 
\frac {10}{3}e^{3x_{2}}- \no\\
&&- \frac {35}{4}e^{4x_{1}}
 - 292e^{3x_{1}+x_{2}} - 749e^{2x_{1}+2x_{2}}- 
292e^{x_{1}+3x_{2}}- \frac {35}{4} e^{4x_{2}}-\cdots.
\label{nr1}
\end{eqnarray} 
We introduce here an auxiliary generating function:
\begin{eqnarray}
w({\cal O}_{1}{\cal O}_{z^2})_{0}=kx_1-kx_2.
\end{eqnarray}
Since $z^2=0$ in $H^{*}(F_{0},C)$, $w({\cal O}_{1}{\cal O}_{z^2})_{0,\dd}=0\;(\dd>(0,0))$ by definition. 
With these results, we can confirm that
\begin{eqnarray}
t_{1}(x_1,x_2)=\eta^{z\alpha}w({\cal O}_{1}{\cal O}_{\alpha})_{0}
&=&x_{1}+2e^{x_1}+2e^{x_2}+3e^{2x_1}+12e^{x_1+x_2}+3e^{2x_2}+\no\\
&&+\frac{20}{3}e^{3x_1}+60e^{2x_1+x_2}+60e^{x_1+2x_2}+\frac{20}{3}e^{3x_2}+\no\\
&&+\frac{35}{2}e^{4x_1}+280e^{3x_1+x_2}+
630e^{2x_1+2x_2}+280e^{x_1+3x_2}+\frac{35}{2}e^{4x_2}+\cdots
,\no\\
t_{2}(x_1,x_2)=\eta^{w\alpha}w({\cal O}_{1}{\cal O}_{\alpha})_{0}
&=&x_{2}+2e^{x_1}+2e^{x_2}+3e^{2x_1}+12e^{x_1+x_2}+3e^{2x_2}+\no\\
&&+\frac{20}{3}e^{3x_1}+60e^{2x_1+x_2}+60e^{x_1+2x_2}+\frac{20}{3}e^{3x_2}+\no\\
&&+\frac{35}{2}e^{4x_1}+280e^{3x_1+x_2}+
630e^{2x_1+2x_2}+280e^{x_1+3x_2}+\frac{35}{2}e^{4x_2}+\cdots
,\no\\
\label{mmap1}
\end{eqnarray}
coincide with the mirror map obtained from the standard Picard-Fuchs system 
used in mirror computation of $K_{F_{0}}$. Finally, we invert (\ref{mmap1}) and 
substitute $x_1=x_1(t_1,t_2)$ and $x_2=x_2(t_1,t_2)$ into $w({\cal O}_{\alpha}{\cal O}_{\beta})_{0}$
. We show here the result of this operation in the cases of $w({\cal O}_{z}{\cal O}_{z})_{0}$ and 
$w({\cal O}_{z}{\cal O}_{w})_{0}$:  
\begin{eqnarray}
w({\cal O}_{z}{\cal O}_{z})_{0}|_{x_1=x_1(t_1,t_2),x_2=x_2(t_1,t_2)}&=&
k{t_{1}} - k{t_{2}} - 2e^{t_{1}}
 -e^{2t_{1}} -4e^{t_{1}+t_{2}} - 
\frac {2}{3}e^{3t_{1}} - 24e^{2t_{1}+t_{2}} - 6e^{t_{1}+2t_{2}}
- \frac {1}{2}e^{4t_{1}}-\no\\
&&- 8e^{t_{1}+3t_{2}}- 72e^{3t_{1}+t_{2}}-130e^{2t_{1}+2t_{2}}-\cdots,\no\\
w({\cal O}_{z}{\cal O}_{w})_{0}|_{x_1=x_1(t_1,t_2),x_2=x_2(t_1,t_2)}&=&
- k{t_{1}} + (k-\frac{1}{2}){t_{2}}- 4e^{t_{1}+t_{2}}
- 12e^{2t_{1}+t_{2}} 
- 12e^{t_{1}+2t_{2}} - 24e^{t_{1}+3t_{2}}-\no \\
&&- 130e^{2t_{1}+2t_{2}} - 24e^{3t_{1}+t_{2}}-\cdots. 
\label{two}
\end{eqnarray}
These results indeed agree with the results of standard computation of local mirror symmetry \cite{ckyz}.
Therefore, they give us numerical evidence of Conjecture 1 in the case of $K_{F_{0}}$. 
\subsection{$F_{3}$}
\subsubsection{Notation and Polynomial Maps}
In this section, we treat Hirzebruch surface $F_{3}$, which is a more challenging example than $F_{0}$. 
In short, it is given as a projective bundle $\pi: {\bf P}( {\cal O}_{{\bf P}^{1}} \oplus {\cal O}_{{\bf P}^{1}}(-3))
\rightarrow {\bf P}^{1}$ and is a well-known example of non-nef complex manifold. Therefore, its quantum cohomology 
is difficult to analyze from the point of view of the mirror computation \cite{fj1}. Let $ {\cal O}_{F_{3}}(a) $ be 
$\pi^{*}{\cal O}_{{\bf P}^{1}}(1)$ and ${\cal O}_{F_{3}}(b)$ be dual line bundle of the universal bundle of 
${\bf P}( {\cal O}_{{\bf P}^{1}} \oplus {\cal O}_{{\bf P}^{1}}(-3))$.  We denote $c_1({\cal O}_{F_{3}}(a))$ 
(resp. $c_1({\cal O}_{F_{3}}(b))$ by $z$ (resp. $w$). $z$ and $w$ generate the cohomology ring $H^{*}(F_{3},C)$ and  
obey the relations:
\begin{eqnarray}      
z^2=0,\;w(w-3z)=0.
\label{relf3}
\end{eqnarray}
In this section, we identify $H^{*}(F_{3},C)$ with $<1,z,w,w^2>_{\bf C}$, i.e., we take $w^2$ as the representative 
of the basis of $H^{2,2}(F_{3},C)$. With these set-up's, integration of cohomology element $\alpha$ over 
$F_{3}$ can be realized by the residue integral:
\begin{eqnarray}
\frac{1}{(2\pi\sqrt{-1})^2}\oint_{C_{0}}\frac{dz}{z^2}\oint_{C_{(0,3z)}}\frac{dw}{w(w-3z)}\alpha.
\label{res3}
\end{eqnarray}
In (\ref{res3}), $\alpha$ should be considered as a polynomial in $z$ and $w$.

Like $F_{0}$, $F_{3}$ has the following toric construction:
\begin{eqnarray}    
F_{3}=\{({\bf a},{\bf b})\;|\; {\bf a}=(a_0,a_1)\in{\bf C}^2,\; {\bf b}=(b_0,b_1)\in{\bf C}^2 
\;,{\bf a}\neq{\bf 0},\;{\bf b}\neq{\bf 0}\;\;\}/({\bf C}^{\times})^2, 
\label{tof3}
\end{eqnarray}
where the two ${\bf C}^{\times}$ actions are given by, 
\begin{eqnarray}
(a_0,a_1,b_0,b_1)\rightarrow (\mu a_0,\mu a_1, {\mu}^{-3}b_0,b_1), \;
(a_0,a_1,b_0,b_1)\rightarrow (a_0,a_1, \nu b_0,\nu b_1). 
\end{eqnarray}
From now on, we denote by $[({\bf a},{\bf b})]$ the equivalence class of $({\bf a},{\bf b})$ 
under these two ${\bf C}^{\times}$ actions. It is well-known that the K\"ahler form $z$ (resp. $w$) is 
associated with the first (resp. second) ${\bf C}^{\times}$ action through the moment map construction. 

We then consider a polynomial map of $F_{3}$ of bi-degree $\dd=(d_a,d_b)$ where $d_a$ (resp. $d_b$) is the degree 
associated with the first (resp. second) ${\bf C}^{\times}$ action. It behaves in a more complicated 
way than in the $F_{0}$ case because the first ${\bf C}^{\times}$ action has a $\mu^{-3}$ factor. Since we consider the 
moduli space of polynomial maps with two marked points, we restrict our attention to polynomial maps such that the
images of $[(1,0)]$ and $[(0,1)]$ are well-defined. If $d_a>0$, a polynomial map of degree $\dd$ satisfying the above 
condition is given by, 
\begin{eqnarray}
( \sum_{j=0}^{\d_a}{\bf a}_js^jt^{d_a-j} , (0,\sum_{j=0}^{d_b}b_{1j}s^jt^{d_b-j})),\;({\bf a}_{0},{\bf a}_{d_a}\neq {\bf 0},\;
b_{10},b_{1 d_b}\neq 0). 
\label{polyf31}
\end{eqnarray}   
The first entry of ${\bf b}$ factor should be $0$ because of the $\mu^{-3}$ factor of the first ${\bf C}^{\times}$ action.
If $d_a=0$, the polynomial map we need is given as follows: 
\begin{eqnarray}
({\bf a}_{0}, \sum_{j=0}^{d_b}{\bf b}_{j}s^jt^{d_b-j})),\;({\bf a}_{0},{\bf b}_{0},{\bf b}_{d_b}\neq {\bf 0}). 
\label{polyf32}
\end{eqnarray}
In the same way as in the $F_{0}$ case, we define the moduli space of polynomial maps with two marked points $ Mp(F_{3},\dd) $:
\begin{eqnarray}    
Mp(F_{3},\dd)=\{({\bf a}_{0},{\bf a}_{1}\cdots,{\bf a}_{d_a},{\bf b}_{0},{\bf b}_1\cdots,{\bf b}_{d_b})\;|\;{\bf a}_{i},{\bf b}_j\in {\bf C}^{2},
{\bf a}_{0},{\bf a}_{d_a},{\bf b}_0,{\bf b}_{d_b}\neq {\bf 0}\}/({\bf C}^{\times})^{3},
\label{mpf3} 
\end{eqnarray}
where we have to set $b_{0j}=0\;(j=0,1,\cdots,d_b)$ if $d_a>0$. The three ${\bf C}^{\times}$ actions are given by, 
\begin{eqnarray}
&& ({\bf a}_{0},\cdots,{\bf a}_{d_{a}},{\bf b}_{0},\cdots,{\bf b}_{d_{b}})  \rightarrow 
(\mu{\bf a}_{0},\cdots,\mu{\bf a}_{d_{a}},{\bf b}_{0},\cdots,{\bf b}_{d_{b}}),\no\\
&&({\bf a}_{0},\cdots,{\bf a}_{d_{a}},{\bf b}_{0},\cdots,{\bf b}_{d_{b}}) \rightarrow 
({\bf a}_{0},\cdots,{\bf a}_{d_{a}},\nu{\bf b}_{0},\cdots,\nu{\bf b}_{d_{b}}),\no\\
&&({\bf a}_{0},\cdots,{\bf a}_{d_{a}},{\bf b}_{0},\cdots,{\bf b}_{d_{b}}) \rightarrow 
({\bf a}_{0},\lambda{\bf a}_{1},{\lambda}^2{\bf a}_{2},\cdots,{\lambda}^{d_{a}}{\bf a}_{d_{a}},
{\bf b}_{0},{\lambda}{\bf b}_{1},{\lambda}^2{\bf b}_{2}\cdots,{\lambda}^{d_{b}}{\bf b}_{d_{b}}).
\end{eqnarray}
The complex dimension of $Mp(F_{0},\dd)$ coincides with the expected dimension $1-d_a+2d_b$ if $d_a=0$. But if $d_a>0$, it becomes 
$2d_a+d_b$ and exceeds the expected dimension 
by $3d_a-d_b-1$. At this stage, we must note that the rational map $\varphi: CP^{1}\rightarrow F_{3}$ induced 
from the polynomial map given in (\ref{polyf31}) has non-trivial obstruction. Let $C\subset F_{3}$ be the image curve of $\varphi$. 
We first assume here that the vector-valued polynomial is "not" factorized into product of a homogeneous polynomial in $s$ and $t$ of 
positive degree $f(>0)$ and a vector-valued polynomial of positive degree:
\begin{eqnarray}
\sum_{j=0}^{\d_a}{\bf a}_js^jt^{d_a-j} =(\sum_{j=0}^{f}p_js^jt^{f-j})\cdot( \sum_{j=0}^{\d_a-f}{\bf a}^{\prime}_js^jt^{d_a-f-j}). 
\end{eqnarray}
Under this assumption, $C$ is identified with a section $\{[({\bf a},(0,1))]\;|\;{\bf a}\in {\bf C}^{2}\}$ and the normal bundle $N_{C}$ of 
$C$ in $F_{3}$ is identified with ${\cal O}_{F_{3}}(b-3a)$ through the Euler sequence:
\begin{eqnarray}
0\rightarrow {\bf C}^{2}\rightarrow {\cal O}_{F_{3}}(a)\oplus {\cal O}_{F_{3}}(a)\oplus {\cal O}_{F_{3}}(b-3a)\oplus {\cal O}_{F_{3}}(b)
\rightarrow T^{\prime}F_{3}\rightarrow 0.
\end{eqnarray}  　  
Since $\varphi^{*}{\cal O}_{F_{3}}(b-3a)={\cal O}_{CP^{1}}(d_b-3d_a)$, we have non-trivial obstruction $H^{1}(CP^1,\varphi^{*}N_C)=H^{1}(CP^1,{\cal O}(d_b-3d_a))$ 
of rank $3d_a-d_b-1$ if $d_b-3d_a<0$.
We can extend this obstruction of rank $3d_a-d_b-1$ to the locus where our assumption is not satisfied by imitating the discussion of 
Subsection 2.3.  
We denote by $Obs$ the rank $3d_a-d_b-1$ bundle on $ Mp(F_{3},\dd) $  so obtained.

Let us now turn to the construction of $\widetilde{Mp}_{0,2}(F_3,\dd)$. Since $F_3$ is non-nef, the boundary components of $\widetilde{Mp}_{0,2}(F_3,\dd)$
behave in a more complicated way than the $F_0$ case. Therefore, it is unclear to us whether there exists a simple toric construction 
like $\widetilde{Mp}_{0,2}(F_0,\dd)$. But we proceed under the assumption that the coordinates $u_{(i_a,i_b)}\;(0\leq i_a\leq d_a,\; 0\leq i_b\leq d_b,\;(i_a,i_b)\neq (0,0),(d_a,d_b))$
used in the $F_0$ case still work in the $F_3$ case. If we set one $u_{(i_a,i_b)}$ to zero, we expect that the following chain of two polynomial maps appears:
\begin{eqnarray}
( \sum_{j=0}^{i_{a}}{\bf a}_{j}s_{1}^{i_{a}-j}t_{1}^{j} , \sum_{j=0}^{i_{b}}{\bf b}_{j}s_{1}^{i_{b}-j}t_{1}^{j})
\cup ( \sum_{j=0}^{d_{a}-i_{a}}{\bf a}_{j}s_{2}^{d_{a}-i_{a}-j}t_{2}^{j} , \sum_{j=0}^{d_{b}-i_{b}}{\bf b}_{j}s_{2}^{d_{b}-i_{b}-j}t_{2}^{j}).
\label{chainf31}
\end{eqnarray} 
In the case where $d_a>0$, we have to deal with the behavior of the ${\bf b}_j$'s carefully. If $1\leq i_a\leq d_a-1$, $b_{0j}$ must be zero for all $j=0,1,\cdots,d_b$.
But if $i_a=0$, the first polynomial map becomes a polynomial map of degree $(0,i_b)$. Hence $b_{0h}$ can take  arbitrary values if $0\leq h\leq i_b-1$. 
On the other hand, $b_{0i_b}$ have to be zero since the second map is a polynomial map of degree $(d_a, d_b-i_b)$. If we set $i_a$ to $d_a$, we come 
across the same exotic behavior with the roles of of the first and the second polynomial maps interchanged. Let us compare the dimension of 
this boundary locus with dimension of $ Mp_{0,2}(F_3,\dd) $. If $d_a>0$, $\dim_{\bf C}(Mp_{0,2}(F_3,\dd) )$ equals $2d_a+d_b$. In the case where 
$1\leq i_a\leq d_a-i$, the dimension of the boundary locus is given by $2d_a+d_b-1$.  But if $i_a=0$ (resp. $i_a=d_a$), the dimension of the boundary 
locus becomes $2d_a+d_b+j-1$ (resp. $2d_a+2d_b-j-1$). Therefore, we are confronted with the singular phenomena that the dimension of the 
boundary locus exceeds the dimension of  $ Mp_{0,2}(F_3,\dd) $.  In such cases, we have to consider the rank of obstruction together with the dimension. 
As was computed before,   $\dim_{\bf C}(Mp_{0,2}(F_3,\dd) )-\mbox{rank}(Obs)=1-d_a+2d_b$. We can also define the obstruction of the 
chain of two polynomial maps in (\ref{chainf31}). If $1\leq i_a\leq d_a-1$, the obstruction is given by,
\begin{eqnarray}
H^1(CP^1,\varphi^{*}_1 {\cal O}_{F_{3}}(-3a+b) ) \oplus  {\cal O}_{F_{3}}(-3a+b) \oplus H^1(CP^1,\varphi^{*}_2 {\cal O}_{F_{3}}(-3a+b) ),
\label{bdh0}
\end{eqnarray}
where $\varphi_1$ (resp.  $\varphi_2$) is the rational map induced from the first (resp. the second) polynomial map in (\ref{chainf31}).   
Its rank equals,
\begin{eqnarray}
3i_a-(i_b-1)-1+1+3(d_a-i_a)-(d_b-i_b)-1=3d_a-d_b-1.
\label{bdh1}
\end{eqnarray}
Therefore, dimension of the boundary locus minus the rank of obstruction becomes $-d_a+2d_b$, which 
is less than $1-d_a+2d_b$ by $1$.
If $i_a=0$, the obstruction arises only from the second polynomial map, and its rank equals 
$3d_a-(d_b-i_b)-1$. Hence the dimension minus the rank turns out to be,
\begin{eqnarray}
2d_a+d_b+i_b-1-(3d_a-(d_b-i_b)-1)=-d_a+2d_b,
\label{bdh12}
\end{eqnarray} 
which is also less than the expected dimension of $Mp_{0,2}(F_3,\dd)$ by $1$.
We come to the same conclusion in the $i_a=d_a$ case. In this way, 
we can conclude that the expected dimension of the locus $u_{(i_a,i_b)}=0$ behaves well 
in the $F_3$ case. In general, we have to consider the locus:
\begin{eqnarray}
&&u_{(i_{a,1},i_{b,1})}=u_{(i_{a,2},i_{b,2})}=\cdots=u_{(i_{a,l},i_{b,l})}=0,\no\\
&&((0,0)=(i_{a,0},i_{b,0})<(i_{a,1},i_{b,1})<(i_{a,2},i_{a,2})<\cdots
<(i_{a,l},i_{a,l})=(d_a,d_b)). 
\label{seqf3}
\end{eqnarray} 
In the same way as in the $F_{0}$ case, we can associate a chain of $l$ polynomial maps 
to a point in this locus:
\begin{eqnarray}
\mathop{\cup}_{k=1}^{l}(\sum_{h=0}^{i_{a,k}-i_{a,k-1}}{\bf a}_{i_{a,k-1}+h}
(s_k)^h(t_k)^{i_{a,k}-i_{a,k-1}-h},\sum_{h=0}^{i_{b,k}-i_{b,k-1}}{\bf b}_{i_{b,k-1}+h}
(s_k)^h(t_k)^{i_{b,k}-i_{b,k-1}-h}).
\end{eqnarray}
But we have to impose the following conditions on ${\bf b}_j$'s:
\begin{itemize}
\item[(a)] If $i_{a,k}-i_{a,k-1}>0$, $b_{0h}$ $(h=i_{b,k-1},i_{b,k-1}+1,\cdots,i_{b,k})$
is $0$.
\item[(b)] If $i_{a,k}-i_{a,k-1}=0$, $b_{0h}$ $(h=i_{b,k-1}+1,i_{b,k-1}+2,\cdots,i_{b,k}-1)$
can take arbitrary value.
\item[(c)] If $i_{a,k+1}-i_{a,k}=0$ and $i_{a,k}-i_{a,k-1}=0$, $b_{0i_{b,k}}$ 
can take arbitrary value. (Otherwise, it is $0$.)
\end{itemize}
We can also define the obstruction of this chain of polynomial maps. Hence we can 
extend the bundle $Obs$ as a sheaf on the whole $\widetilde{Mp}_{0,2}(F_3,\dd)$.
\subsubsection{Virtual Structure  Constants and the Localization Computation}
In this section, we define and compute an analogue of Gromov-Witten invariants of $F_{3}$:
\begin{eqnarray}
\langle{\cal O}_{\alpha}{\cal O}_{\beta}\rangle_{0,2}=
\int_{[\overline{M}_{0,2}(F_3,\dd)]_{vir.}}ev_{1}^{*}(\alpha)\wedge ev_{2}^{*}(\beta).
\label{gwf3}
\end{eqnarray}
In (\ref{gwf3}), $ [\overline{M}_{0,2}(F_3,\dd)]_{vir.}$ is the virtual fundamental class of the 
moduli space $ \overline{M}_{0,2}(F_3,\dd)$, which means automatic insertion of the top Chern class of 
obstruction sheaf. As in the case of $F_0$, $\alpha$ and $\beta$ are elements of the classical 
cohomology ring $H^{*}(F_3,C)$ and $ev_{i}: \overline{M}_{0,2}(F_3,\dd)\rightarrow F_3$ is the 
evaluation map at the $i$-th marked point. To define an intersection number of 
$\widetilde{Mp}_{0,2}(F_3,\dd)$, which we expect to have geometrical meaning parallel
to $\langle{\cal O}_{\alpha}{\cal O}_{\beta}\rangle_{0,2}$, we introduce the heuristic 
notation to represent a point of $ \widetilde{Mp}_{0,2}(F_3,\dd)$:
\begin{equation}
[({\bf a}_{0},{\bf a}_{1},\cdots,{\bf a}_{d_a},{\bf b}_{0},{\bf b}_{1},\cdots,{\bf b}_{d_b},
u_{(1,0)},\cdots,u_{(d_a-1,d_b)})]\in \widetilde{Mp}_{0,2}(F_3,\dd).
\label{nof3}
\end{equation}
This is not rigorous in the sense that we haven't specified the equivalence relations which should 
come from the ${\bf C}^{\times}$ actions, but it is sufficient for our present purpose. Of course, 
the ${\bf b}_j$'s must obey the conditions (a), (b) and (c). With this notation, 
we define evaluation maps $ev_1$ and $ev_2$ from $ \widetilde{Mp}_{0,2}(F_3,\dd)$ to 
$F_3$ as follows: 
\begin{eqnarray}
&&ev_{1}( [({\bf a}_{0},{\bf a}_{1},\cdots,{\bf a}_{d_a},{\bf b}_{0},{\bf b}_{1},\cdots,{\bf b}_{d_b},
u_{(1,0)},\cdots,u_{(d_a-1,d_b)})])=[({\bf a}_{0},{\bf b}_{0})],\no\\
&&ev_{2}( [({\bf a}_{0},{\bf a}_{1},\cdots,{\bf a}_{d_a},{\bf b}_{0},{\bf b}_{1},\cdots,{\bf b}_{d_b},
u_{(1,0)},\cdots,u_{(d_a-1,d_b)})])=[({\bf a}_{d_a},{\bf b}_{d_b})].
\label{evf3}
\end{eqnarray}
We also define the virtual fundamental class $[\widetilde{Mp}_{0,2}(F_3,\dd)]_{vir.}$, which means automatic 
insertion of the top Chern class of the sheaf $Obs$ on $ \widetilde{Mp}_{0,2}(F_3,\dd)$. With this setup, we define an intersection number analogous to $\langle{\cal O}_{\alpha}{\cal O}_{\beta}\rangle_{0,2}$
as follows:
\begin{equation}
w({\cal O}_{\alpha}{\cal O}_{\beta})_{0,2}:=
\int_{[\widetilde{Mp}_{0,2}(F_3,\dd)]_{vir.}}ev_{1}^{*}(\alpha)\wedge ev_{2}^{*}(\beta).
\label{vscf3}
\end{equation}
Now, we compute $w({\cal O}_{\alpha}{\cal O}_{\beta})_{0,2}$ by using the localization theorem. As in the $F_0$ case, 
we consider the torus action flow:
\begin{eqnarray}
[( e^{\lambda_{0}t}{\bf a}_{0},  e^{\lambda_{1}t}{\bf a}_{1},\cdots,  e^{\lambda_{d_a}t}{\bf a}_{d_{a}},
 e^{\mu_{0}t}{\bf b}_{0},  e^{\mu_{1}t}{\bf b}_{1},\cdots,  e^{\mu_{d_b}t}{\bf b}_{d_{b}}, u_{(1,0)} ,\cdots,u_{(d_a-1,d_b)})].
\label{flowf3}
\end{eqnarray} 
In the same way as in the $F_0$ case, connected components of fixed point set under the above flow are classified 
by ordered partition $\sigma_{\dd}$, which is an element of the following set:
\begin{equation}
OP_{\dd}=\{\sigma_{\dd}=(\dd_{1},\dd_{2},\cdots,\dd_{l(\sigma_{\dd})})\;\;|\;\;
\sum_{j=1}^{l(\sigma_{d})}\dd_{j}=\dd\;\;,\;\;\dd_{j}=(d_{a,j},0)\;\mbox{or}\;\dd_{j}=(0,d_{b,j})\}.
\label{partf3} 
\end{equation}
Let $F_{\sigma_{\dd}}$ be a connected component of the fixed point set labeled by $\sigma_{\dd}$.
A point in $F_{\sigma_{\dd}}$ is represented by a chain of polynomial maps of length $l(\sigma_{\dd})$ whose 
$j$-th component is given by,
 \begin{eqnarray}
( {\bf a}_{j-1}(s_{j})^{|\dd_{j}|}+{\bf a}_{j}(t_{j})^{|\dd_{j}|} ,{\bf b}_{j-1})\; \mbox{or} \;({\bf a}_{j-1},{\bf b}_{j-1}(s_{j})^{|\dd_{j}|}+{\bf b}_{j}(t_{j})^{|\dd_{j}|}),
\label{chf3}
\end{eqnarray} 
respectively if $\dd_j=(d_{a,j},0)$ or $\dd_j=(0,d_{b,j})$. In (\ref{chf3}), we relabel the subscripts of ${\bf a}_{j}$'s and 
${\bf b}_j$'s in the same manner as in the $F_{0}$ case. We must be careful of the behavior of the ${\bf b}_j$'s because 
we have non-trivial restrictions imposed by the conditions (a), (b) and (c) in the previous sub-subsection.
If $\dd_j=(d_{a,j},0)$, the first entries of ${\bf b}_{j-1}$ and ${\bf b}_{j}$ should be $0$ because of the condition (a).
${\bf b}_j$ can take any value of ${\bf C}^{2}$ only if $\dd_j=(0,d_{b,j})$ and $\dd_{j+1}=(0,d_{b,j+1})$. (Precisely speaking, 
${\bf b}_{0}$ (resp. ${\bf b}_{l(\sigma_{\dd})}$) can take any value of ${\bf C}^{2}$ if $\dd_1=(0,d_{b,1})$ (resp. 
$\dd_{l(\sigma_{\dd})}=(0,d_{b,l(\sigma_{\dd})})$).)   
Therefore, it is set-theoretically given by a subset of,
\begin{eqnarray}
(F_{3})_{0} \times(F_{3})_{1} \times(F_{3})_{2} \times\cdots \times(F_{3})_{l(\sigma_{\dd})},\; ((F_{3})_j=\{[({\bf a}_j,{\bf b}_j)]\}),
\label{prof3}
\end{eqnarray} 
defined by the following conditions:
\begin{eqnarray}
&&{\bf b}_{j-1}={\bf b}_{j}=(0,1)\;\;\;\mbox{if}\; \dd_{j}=(d_{a,j},0),\no\\
&&{\bf a}_{j-1}={\bf a}_{j} \;\;\;\mbox{if}\; \dd_{j}=(0,d_{b,j}).
\label{idf3}
\end{eqnarray}
In (\ref{idf3}), we used the trivial fact that $[({\bf a},(0,b))]=[({\bf a},(0,1))]$. As usual, $F_{\sigma_{\dd}}$ should be regarded 
as an orbifold on which an abelian group $\oplus_{j=1}^{l(\sigma_{\dd})}\bigl({\bf Z}/(|\dd_{j}|{\bf Z})\bigr)$ acts.   
For later use, we introduce the inclusion map,
\begin{eqnarray}
i:F_{\sigma_{\dd}}\rightarrow (F_{3})_{0} \times(F_{3})_{1} \times(F_{3})_{2} \times\cdots \times(F_{3})_{l(\sigma_{\dd})},
\label{inf3}
\end{eqnarray}
and the projection map,
\begin{eqnarray}
\pi_{j}:(F_{3})_{0} \times(F_{3})_{1} \times(F_{3})_{2} \times\cdots \times(F_{3})_{l(\sigma_{\dd})}\rightarrow (F_{3})_j,
\;(j=0,1,\cdots,l(\sigma_{\dd})). 
\label{prf3}
\end{eqnarray}
Next, we determine the normal bundle $N_{F_{\sigma_{\dd}}}$ of $F_{\sigma_{\dd}}$ in $\widetilde{Mp}_{0,2}(F_3,\dd)$.
It consists of the degrees of freedom of deforming polynomial maps in $F_{3}$ and of resolving singularities of the image 
curve. Let $N_{\dd_j}$ be the direct summand of $N_{F_{\sigma_{\dd}}}$ coming from deformation of the polynomial map 
of degree $\dd_j$ and $N_{(\dd_j-1,\dd_j)}$ be the one coming from the resolution of the singularity between the polynomial 
maps of degree $\dd_{j-1}$ and $\dd_j$. Obviously, we have,
\begin{eqnarray}
N_{\sigma_{\dd}}=\bigl(\mathop{\oplus}_{j=1}^{l(\sigma_{\dd})}N_{\dd_j}\bigr)\oplus
\bigl(\mathop{\oplus}_{j=1}^{l(\sigma_{\dd})-1}N_{(\dd_{j-1},\dd_j)}\bigr).
\label{nf3}
\end{eqnarray} 
Following the $F_{0}$ case, we introduce the notation:
\begin{eqnarray}
{\cal O}_{F_{3}}(\frac{m}{\dd_{j}}):=\left\{\begin{array}{cc}  {\cal O}_{F_{3}}(\frac{m}{|\dd_{j}|}a) &\;\mbox{if}\;\dd_j=(d_{a,j},0),\\
  {\cal O}_{F_{3}}(\frac{m}{|\dd_{j}|}b) &\;\mbox{if}\;\dd_j=(0,d_{b,j}).\end{array}\right.
\end{eqnarray}  
Then $N_{\dd_j}$ and $N_{(\dd_{j-1},\dd_j)}$ are given as follows:
\begin{eqnarray}
N_{\dd_j}&=&\left\{\begin{array}{cc}i^{*}(\mathop{\oplus}_{i=1}^{|\dd_j|-1}
(\pi_{j-1}^{*}{\cal O}_{F_3}(\frac{-i}{\dd_j})\otimes\pi_{j}^{*}{\cal O}_{F_3}(\frac{i}{\dd_j})\otimes
\pi_{j-1}^{*}{\cal O}_{F_3}(a)\bigr)^{\oplus 2} &\mbox{if}\;\dd_j=(d_{a,j},0),\\
&\\
i^{*}\bigl(\mathop{\oplus}_{i=1}^{|\dd_j|-1}
(\pi_{j-1}^{*}{\cal O}_{F_3}(\frac{-i}{\dd_j})\otimes\pi_{j}^{*}{\cal O}_{F_3}(\frac{i}{\dd_j})\otimes
\pi_{j-1}^{*}{\cal O}_{F_3}(b-3a)\bigr)\oplus & \\ 
\oplus i^{*}\bigl(\mathop{\oplus}_{i=1}^{|\dd_j|-1}
(\pi_{j-1}^{*}{\cal O}_{F_3}(\frac{-i}{\dd_j})\otimes\pi_{j}^{*}{\cal O}_{F_3}(\frac{i}{\dd_j})\otimes
\pi_{j-1}^{*}{\cal O}_{F_3}(b)\bigr)
 &\mbox{if}\;\dd_j=(0,d_{b,j}),\end{array}\right.\no\\
N_{(\dd_{j-1},\dd_j)}&=&i^{*}\bigl(\pi^{*}_{j-1}{\cal O}_{F_3}(\frac{-1}{\dd_{j-1}})\oplus \pi^{*}_{j}{\cal O}_{F_3}(\frac{1}{\dd_{j-1}})\oplus 
\pi^{*}_{j}{\cal O}_{F_3}(\frac{1}{\dd_{j}})\oplus \pi^{*}_{j+1}{\cal O}_{F_3}(\frac{-1}{\dd_{j}})\bigr). 
\label{dj}
\end{eqnarray}
In the $F_3$ case, we also have to determine $Obs_{\sigma_{\dd}}$: restriction of $Obs$ to $F_{\sigma_{\dd}}$.
Let $Obs_{\dd_j}$ be the direct summand of $Obs_{\sigma_{\dd}}$ coming from the obstruction of deforming the polynomial 
map of degree $\dd_j$ and  $Obs_{(\dd_j-1,\dd_j)}$ be the one that arises from the effect of nodal singularities 
between the polynomial maps of degree $\dd_{j-1}$ and $\dd_j$. In the same way as in the $N_{F_{\sigma_{\dd}}}$
case, we have 
\begin{eqnarray}
Obs_{\sigma_{\dd}}=\bigl(\mathop{\oplus}_{j=1}^{l(\sigma_{\dd})}Obs_{\dd_j}\bigr)\oplus
\bigl(\mathop{\oplus}_{j=1}^{l(\sigma_{\dd})-1}Obs_{(\dd_{j-1},\dd_j)}\bigr).
\label{nif3}
\end{eqnarray} 
These direct summands turn out to be, 
\begin{eqnarray}
Obs_{\dd_j}&=&\left\{\begin{array}{cc}i^{*}(\mathop{\oplus}_{i=1}^{3|\dd_j|-1}
(\pi_{j-1}^{*}{\cal O}_{F_3}(\frac{-i}{\dd_j})\otimes\pi_{j}^{*}{\cal O}_{F_3}(\frac{i}{\dd_j})\otimes
\pi_{j-1}^{*}{\cal O}_{F_3}(-3a+b)\bigr)&\mbox{if}\;\dd_j=(d_{a,j},0),\\
0
 &\mbox{if}\;\dd_j=(0,d_{b,j}).\end{array}\right.\no\\
Obs_{(\dd_{j-1},\dd_j)}&=&\left\{\begin{array}{cc}
i^{*}\pi_{j}^{*}{\cal O}_{F_3}(-3a+b) &\mbox{if}\;\dd_{j-1}=(d_{a,j-1},0)\;\mbox{and}\;\dd_{j}=(d_{a,j},0),
\\
0
 &\mbox{otherwise}.\end{array}\right.
\label{dj2}
\end{eqnarray}
We then move on to the evaluation of the contribution from $F_{\sigma_{\dd}}$ to $w({\cal O}_{\alpha}{\cal O}_{\beta})_{0,\dd}$. 
We denote $\pi_{j}^{*}(z)$ 
(resp. $\pi_{j}^{*}(w)$) by $z_j$ (resp. $w_j$). As in the $F_{0}$ case, we define the following rational 
function to express the contributions from $N_{\dd_j}$ and $Obs_{\dd_j}$:
\begin{eqnarray}
G(\dd;z_{0},z_{1},w_0,w_1):=
\left\{\begin{array}{cc}
\frac{\prod_{j=1}^{3|\dd|-1}\bigl(\frac{-jz_0-(3|\dd|-j)z_1}{|\dd|}+w_0\bigr)}{\prod_{j=1}^{|\dd|-1}\bigl(\frac{jz_{0}+(|\dd|-j)z_{1}}{|\dd|}\bigr)^{2}},&\;\mbox{if}\;\dd=(d,0),\\
&\\
\frac{1}{\prod_{j=1}^{|\dd|-1}\bigl(( \frac{jw_{0}+(|\dd|-j)w_{1}}{|\dd|}-3z_{0} )(\frac{jw_{0}+(|\dd|-j)w_{1}}{|\dd|})\bigr)},
&\;\mbox{if}\;\dd=(0,d).
\end{array}\right.
\label{Gf3}
\end{eqnarray}
To express the contributions from $N_{(\dd_{j-1},\dd_j)}$ and $Obs_{(\dd_{j-1},\dd_j)}$, we introduce another rational function:
\begin{eqnarray}
H(\dd_1;\dd_2,z_0,z_1,z_2,w_0,w_1,w_2):=\left\{\begin{array}{cc}
\frac{(-3z_1+w_1)}{ \bigl(\frac{z_1-z_0}{|\dd_1|}+ 
\frac{z_1-z_2}{|\dd_2|}\bigr)},&\;\mbox{if}\;\dd_{1}=(d_{a,1},0)\;\mbox{and}\;\dd_{2}=(d_{a,2},0),\\
\frac{1}{ \bigl(\frac{z_1-z_0}{|\dd_1|}+ 
\frac{w_1-w_2}{|\dd_2|}\bigr)},&\;\mbox{if}\;\dd_{1}=(d_{a,1},0)\;\mbox{and}\;\dd_{2}=(0,d_{b,2}), \\
\frac{1}{ \bigl(\frac{w_1-w_0}{|\dd_1|}+ 
\frac{z_1-z_2}{|\dd_2|}\bigr)},&\;\mbox{if}\;\dd_{1}=(0,d_{b,1})\;\mbox{and}\;\dd_{2}=(d_{a,2},0),\\
\frac{1}{ \bigl(\frac{w_1-w_0}{|\dd_1|}+ 
\frac{w_1-w_2}{|\dd_2|}\bigr)},&\;\mbox{if}\;\dd_{1}=(0,d_{b,1})\;\mbox{and}\;\dd_{2}=(0,d_{b,2}).
\end{array}\right.
\label{Hf3}
\end{eqnarray}
With this setup, the contributions from $Obs$ and the normal bundle 
of $F_{\sigma_{\dd}}$ can be collected in the following integrand:
\begin{eqnarray}
K(\sigma_{\dd};z_{*},w_{*}):= \prod_{j=1}^{l(\sigma_{\dd})}G(\dd_j;z_{j-1},z_j,w_{j-1},w_j)
 \prod_{j=1}^{l(\sigma_{\dd})-1}H(\dd_j;\dd_{j+1},z_{j-1},z_j,z_{j+1},w_{j-1},w_j,w_{j+1}).
\label{Kf3}
\end{eqnarray}
Contributions from $ev_{1}^{*}(\alpha)$ and $ev_{2}^{*}(\beta)$ are given in the same way as in the $F_{0}$ case 
as follows: 
\begin{eqnarray}
ev_{1}^{*}(z^sw^t)=(z_{0})^s(w_{0})^t, \; ev_{2}^{*}(z^sw^t)=(z_{l(\sigma_{\dd})})^s(w_{l(\sigma_{\dd})})^t. 
\end{eqnarray} 
In integrating out $ ev_{1}^{*}(\alpha)ev_{2}^{*}(\beta)K(\sigma_{\dd};z_{*},w_{*}) $
over $F_{\sigma_{\dd}}$, we have to note the following three facts: 
\begin{itemize}
\item[(i)] Integration of the cohomology element $\alpha\in H^{*}(F_{3},C)$ is realized as the residue 
integral given in (\ref{res3}).
\item[(ii)] Looking at (\ref{prof3}) and (\ref{idf3}), we must identify:
\begin{eqnarray}
&&w_{j-1}= w_{j}=0\;\ ;\;\mbox{if}\; \dd_{j}=(d_{a,j},0),\no\\
&&z_{j-1}=z_{j} \;\;\;\mbox{if}\; \dd_{j}=(0,d_{b,j}).
\label{id2-2}
\end{eqnarray}
\item[(iii)]
$F_{\sigma_{\dd}}$ should be considered as an orbifold on which an abelian group 
$\oplus_{j=1}^{l(\sigma_{\dd})}\bigl({\bf Z}/(|\dd_{j}|{\bf Z})\bigr)$ acts.
\end{itemize}
Taking the facts (i) and (ii) into account, we define the following operation on 
rational functions $f$ in $z_*$ and$w_*$:
\begin{eqnarray}
Res_{(F_{3})_0}(f)&:=&\left\{\begin{array}{cc}
 \frac{1}{2\pi\sqrt{-1}}\oint_{C_{0}}\frac{dz_0}{(z_0)^2}(f|_{w_0=w_1}),
&\;\mbox{if}\; \dd_{1}=(d_{a,1},0),\\
 \bigl(\frac{1}{2\pi\sqrt{-1}}\oint_{C_{(0,3z_0)}}\frac{dw_0}{w_0(w_0-3z_0)}f\bigr)|_{z_0=z_{1}},
&\;\mbox{if}\; \dd_{1}=(0,d_{b,1}),
\end{array}\right.\no\\ 
Res_{(F_{3})_j}(f)&:=&\left\{\begin{array}{cc}
 \frac{1}{2\pi\sqrt{-1}}\oint_{C_{0}}\frac{dz_j}{(z_j)^2}(f|_{w_j=w_{j+1}}),
&\;\mbox{if}\; \dd_{j+1}=(d_{a,j+1},0),\\
 \bigl(\frac{1}{2\pi\sqrt{-1}}\oint_{C_{0}}\frac{dw_j}{w_j}f\bigr)|_{z_j=z_{j+1}},
&\;\mbox{if}\; \dd_{j}=(d_{a,j},0)\; \mbox{and}\; \dd_{j+1}=(0,d_{b,j+1}),\\
 \bigl(\frac{1}{2\pi\sqrt{-1}}\oint_{C_{(0,3z_j)}}\frac{dw_j}{w_j(w_j-3z_j)}f\bigr)|_{z_j=z_{j+1}},
&\;\mbox{if}\; \dd_{j}=(0,d_{b,j})\; \mbox{and}\; \dd_{j+1}=(0,d_{b,j+1}),
\end{array}\right.\no\\ 
&&\;\;\;\;\;\;\;\;\;\;\;\;\;\;\;\;\;\;\;\;\;\;\;\;\;\;\;\;\;\;\;\;\;\;\;\;\;\;\;\;\;\;\;\;\;\;\;\;\;\;\;\;\;
\;\;\;\;\;\;\;\;\;\;\;\;\;\;\;\;\;\;\;\;\;\;\;\;\;\;\;\;\;\;\;\;\;\;\;\;\;\;\;\;
(j=1,\cdots,l(\sigma_{\dd})-1),\no\\
Res_{(F_{3})_{l(\sigma_{\dd})}}(f)&:=&\left\{\begin{array}{cc}
\frac{1}{(2\pi\sqrt{-1})^2}\oint_{C_{0}}\frac{dz_{l(\sigma_{\dd})}}{(z_{l(\sigma_{\dd})})^2}\oint_{C_{0}}
\frac{dw_{l(\sigma_{\dd})}}{w_{l(\sigma_{\dd})}}f, 
&\;\mbox{if}\; \dd_{l(\sigma_{dd})}=(d_{a,l(\sigma_{\dd})},0),\\
 \frac{1}{(2\pi\sqrt{-1})^2}\oint_{C_{0}}\frac{dz_{l(\sigma_{\dd})}}{(z_{l(\sigma_{\dd})})^2}
\oint_{C_{(0,3z_{l(\sigma_{\dd})})}}\frac{dw_{l(\sigma_{\dd})}}
{w_{l(\sigma_{\dd})}(w_{l(\sigma_{\dd})}-3z_{l(\sigma_{\dd})})}f,
&\;\mbox{if}\; \dd_{l(\sigma_{\dd})}=(0,d_{b,l(\sigma_{\dd})}).
\end{array}\right.\no\\
\label{resf3}
\end{eqnarray}
Integration over $F_{\sigma_{\dd}}$ is done by successive use of the above operation and by dividing the result by the order 
of the abelian group $\oplus_{j=1}^{l(\sigma_{\dd})}\bigl({\bf Z}/(|\dd_{j}|{\bf Z})\bigr)$.
\begin{eqnarray}
&& Amp(\sigma_{\dd};\alpha,\beta) =\biggl(\prod_{j=1}^{l(\sigma_{\dd})}\frac{1}{|\dd_{j}|}\biggr) 
Res_{(F_{3})_{l(\sigma_{\dd})}}( Res_{(F_{3})_{l(\sigma_{\dd})-1}}(
\cdots Res_{(F_{3})_{0}}( ev_{1}^{*}(\alpha)ev_{2}^{*}(\beta)K(\sigma_{\dd};z_{*},w_{*}) )\cdots)).\no\\
\label{amp3}
\end{eqnarray}
Finally, we add up contributions from all the $F_{\sigma_{\dd}}$'s and obtain the formula:
\begin{eqnarray}
 w({\cal O}_{\alpha}{\cal O}_{\beta})_{0,\dd} =\sum_{\sigma_{\dd}\in OP_{\dd}}Amp(\sigma_{\dd};\alpha,\beta). 
\label{amp3f} 
\end{eqnarray}
\subsubsection{Numerical Results and the Mirror Computation}
In this section, we present the numerical results of $w({\cal O}_{\alpha}{\cal O}_{\beta})_{0,\dd}$
by using the formula (\ref{amp3f}). The topological selection rule for $w({\cal O}_{\alpha}{\cal O}_{\beta})_{0,\dd}$
is the same as the one for $\langle{\cal O}_{\alpha}{\cal O}_{\beta}\rangle_{0,\dd}$, as 
can be easily seen from dimension counting. Therefore, $w({\cal O}_{\alpha}{\cal O}_{\beta})_{0,\dd}$
is non-zero only if 
\begin{eqnarray}
\deg(\alpha)+\deg(\beta)=1-d_a+2d_b.
\label{sel3}
\end{eqnarray}
In (\ref{sel3}), $\deg(*)$ means total degree $s+t$ of a cohomology element $z^sw^t$. We write down below non-vanishing 
$w({\cal O}_{\alpha}{\cal O}_{\beta})_{0,\dd}$ up to $d_a=3$.
\begin{eqnarray}
&&w({\cal O}_{1}{\cal O}_{1})_{0,(1,0)}=5,\;\; w({\cal O}_{w}{\cal O}_{w^2})_{0,(0,1)}=3,\no\\
&&w({\cal O}_{1}{\cal O}_{w^2})_{0,(1,1)}=-6,\;\; w({\cal O}_{z}{\cal O}_{z})_{0,(1,1)}=1, \;\; 
w({\cal O}_{z}{\cal O}_{w})_{0,(1,1)}=-1,\no\\
&&w({\cal O}_{1}{\cal O}_{z})_{0,(2,1)}=-16,\;\; w({\cal O}_{1}{\cal O}_{w})_{0,(2,1)}=\frac{39}{2},\no\\
&&w({\cal O}_{1}{\cal O}_{1})_{0,(3,1)}=\frac{1901}{3},\no\\
&&w({\cal O}_{z}{\cal O}_{w^2})_{0,(2,2)}=15,\;\; w({\cal O}_{w}{\cal O}_{w^2})_{0,(2,2)}=-18,\no\\
&&w({\cal O}_{1}{\cal O}_{w^2})_{0,(3,2)}=-\frac{1035}{2},\;\; w({\cal O}_{z}{\cal O}_{z})_{0,(3,2)}=64,\;
w({\cal O}_{z}{\cal O}_{w})_{0,(3,2)}=-96,\;\; w({\cal O}_{w}{\cal O}_{w})_{0,(3,2)}=\frac{413}{3},\no\\
&&w({\cal O}_{w^2}{\cal O}_{w^2})_{0,(3,3)}=432.
\label{num3}
\end{eqnarray}
Then we compare these results with the B-model data used in the mirror computation of $F_{3}$ \cite{fj1}. 
In \cite{fj1}, we started from the so-called $I$-function of $F_{3}$,
\begin{eqnarray}
	I_{F_3}=e^{(zx_1 + wx_2)/\hbar}\sum_{\dd}\frac{\prod_{m=-\infty}^0(-3z+w+m\hbar)}{\prod_{m=-\infty}^{-3d_a+d_b}(-3z+w+m\hbar)\prod_{m=1}^{d_a}(z+m\hbar)^2\prod_{m=1}^{d_b}(w+m\hbar)}e^{d_ax_1+d_bx_2},
\label{If3}
\end{eqnarray}
and applied Birkhoff factorization with respect to the $\hbar$ parameter, to the connection matrix associated with $I_{F_3}$. 
This operation has been explained in Section 1. It resulted in the following two connection matrices: 
\begin{eqnarray}
B_{z}&:=&\left(\begin{array}{cccc}
-2{q_{1}}{q_{2}}-\frac {1035}{2}q_{1}^{3}q_{2}^{2}
& 1 + 135{q_{1}}^{2}{q_{2}}&-32{q_{1}}^{2}{q_{2}}
&\frac {5}{3}{q_{1}}+\frac {1901}{3}{q_{1}}^{3}{q_{2}}\\
10q_{1}^{2}q_{2}^{2}&-864q_{1}^{3}q_{2}^{2}-4{q_{1}}{q_{2}}& 
192q_{1}^{3}q_{2}^{2}+{q_{1}}q_{2}
&- \frac {32}{3}{q_{1}}^{2}{q_{2}}\\ 
 - 12q_{1}^{2}q_{2}^{2}& 1277q_{1}^{3}q_{2}^{2}+3{q_{1}}{q_{2}}& 
- 288q_{1}^{3}q_{2}^{2}-{q_{1}}{q_{2}}&  
\frac {1}{3}+ 13{q_{1}}^{2}{q_{2}}\\
 432{q_{1}}^{3}{q_{2}}^{3}+3{q_{1}}{q_{2}}^{2}& 
-126q_{1}^{2}q_{2}^{2}&30q_{1}^{2}q_{2}^{2}&  
-\frac{1035}{2}q_{1}^{3}q_{2}^{2}-2{q_{1}}
{q_{2}} 
\end{array}\right),\no\\
B_{w}&:=&\left(\begin{array}{cccc}
 - 345q_{1}^{3}q_{2}^{2}-2{q_{1}}{q_{2}}&  
\frac {135}{2}{q_{1}}^{2}{q_{2}}& 1-16
{q_{1}}^{2}{q_{2}}&\frac {1901}{9}{q_{1}}^{3}{q_{2}}\\
 10{q_{1}}^{2}{q_{2}}^{2}& -576q_{1}^{3}q_{2}^{2}-4{q_{1}}{q_{2}}&
128q_{1}^{3}q_{2}^{2}+{q_{1}}{q_{2}}&
\frac{1}{3}-\frac{16}{3}{q_{1}}^{2}{q_{2}}\\
{q_{2}}-12{q_{1}}^{2}{q_{2}}^{2}&\frac {2554}{3}q_{1}^{3}q_{2}^{2}+ 
3{q_{1}}{q_{2}}&-192q_{1}^{3}q_{2}^{2}-{q_{1}}{q_{2}}&  
1+\frac{13}{2}q_{1}^{2}{q_{2}}\\
 6{q_{1}}{q_{2}}^{2}+ 432{q_{1}}^{3}{q_{2}}^{3}& 
3{q_{2}}-126{q_{1}}^{2}{q_{2}}^{2}
&30{q_{1}}^{2}{q_{2}}^{2}&- 345q_{1}^{3}q_{2}^{2}-2{q_{1}}{q_{2}}
\end{array}\right),\no\\
\label{conf3a}
\end{eqnarray}
where $q_1=e^{x_1},\;q_2=e^{x_2}$.
In (\ref{conf3a}), we wrote down the results up to third order in $q_1$. To compare these matrices with (\ref{num3}), 
we multiply them by the classical intersection matrix of $F_3$:
\begin{eqnarray}
\eta:=\bordermatrix{&1&z&w&w^2\cr
1&0&0&0&3\cr
z&0&0&1&0\cr
w&0&1&3&0\cr
w^2&3&0&0&0\cr},
\label{intf3}
\end{eqnarray}
from the right. The results turns out to be,
\begin{eqnarray}
C_{z}&:=&\bordermatrix{&1&z&w&w^2\cr
  1&5{q_{1}} + 1901{q_{1}}^{3}{q_{2}}&  - 32{q_{1}}^{
2}{q_{2}}& 1 + 39{q_{1}}^{2}{q_{2}}&  
 - 6{q_{1}}{q_{2}} - { \frac {3105}{2}} {q_{1}
}^{3}{q_{2}}^{2}   \cr
  z& - 32{q_{1}}^{2}{q_{2}} &{q_{1}}{q_{2}} + 192
{q_{1}}^{3}{q_{2}}^{2}&  - {q_{1}}{q_{2}} - 288{q_{1}
}^{3}{q_{2}}^{2}& 30{q_{1}}^{2}{q_{2}}^{2}  \cr
 w& 1 + 39{q_{1}}^{2}{q_{2}}&  - {q_{1}}{q_{2}} - 
288{q_{1}}^{3}{q_{2}}^{2}& 413{q_{1}}^{3}{q_{2}}^{2}
&  - 36{q_{1}}^{2}{q_{2}}^{2} \cr
w^2&- 6{q_{1}}{q_{2}} -  \frac {3105}{2} {q_{1}}^{3}{q_{2}}^{2}&
 30{q_{1}}^{2}{q_{2}}^{2} 
&  - 36{q_{1}}^{2}{q_{2}}^{2} 
& 9{q_{1}}{q_{2}}^{2} + 1296{q_{1}}^{3}{q
_{2}}^{3} \cr  },\no\\
C_{w}&:=&\bordermatrix{ &1&z&w&w^2\cr
1&\frac {1901}{3}{q_{1}}^{3}{q_{2}}& 1
 - 16{q_{1}}^{2}{q_{2}} & 3 + { \frac {39}{2}} 
{q_{1}}^{2}{q_{2}} &  
 - 6{q_{1}}{q_{2}} - 1035{q_{1}}^{3}{q_{2}}^{2}  \cr
z&1 - 16{q_{1}}^{2}{q_{2}} &
 {q_{1}}{q_{2}} + 128{q_{1}}^{3}{q_{2}}^{2}&
  - {q_{1}}{q_{2}} - 192{q_{1}}^{3}{q_{2}}^{2}& 30{q_{1}}^{2}{q_{2}}^{2}  \cr
w&3 + { \frac {39}{2}} {q_{1}}^{2}{q_{2}}&
  - {q_{1}}{q_{2}} - 192{q_{1}}^{3}{q_{2}}^{2}& 
{ \frac {826}{3}} {q_{1}}^{3}{q_{2}}^{2}& 
3{q_{2}} - 36{q_{1}}^{2}{q_{2}}^{2}   \cr
w^2&- 6{q_{1}}{q_{2}} - 1035{q_{1}}^{3}{q_{2}}^{2}&
 30{q_{1}}^{2}{q_{2}}^{2}&
 3{q_{2}} - 36{q_{1}}^{2}{q_{2}}^{2}&
 18{q_{1}}{q_{2}}^{2} + 1296{q_{1}}^{3}{q_{2}}^{3} \cr}.\no\\
\end{eqnarray}
Let $(C_{z})_{\alpha\beta}(\dd)$ (resp. $ (C_{w})_{\alpha\beta}(\dd)$) be the coefficient of $q_1^{d_a}q_2^{d_b}$ in 
$(C_{z})_{\alpha\beta}$ (resp. $(C_{w})_{\alpha\beta}$). Then we notice that the following equalities
hold true up to the degrees we have computed:
\begin{eqnarray}
d_aw({\cal O}_{\alpha}{\cal O}_{\beta})_{0,\dd}=(C_{z})_{\alpha\beta}(\dd),\;
d_bw({\cal O}_{\alpha}{\cal O}_{\beta})_{0,\dd}=(C_{w})_{\alpha\beta}(\dd).
\label{conjf3}
\end{eqnarray}
Therefore, we confirmed Conjecture 2 for lower degrees.
If Conjecture 2 holds for arbitrary ${\bf d}$, we can construct B-model connection matrices $B_z$ and $B_w$ by using the data of 
$w({\cal O}_{\alpha}{\cal O}_{\beta})_{0,\dd}$'s. Hence we can execute the mirror computation of $F_3$ 
without using the $I$-function and Birkhoff factorization.

\subsection{Calabi-Yau Hypersurface in ${\bf P}(1,1,2,2,2)$}
Originally, ${\bf P}(1,1,2,2,2)$ is a weighted projective space:
\begin{eqnarray}
\{(x_0,x_1,x_2,x_3,x_4)\in{\bf C}^{5}\;|\;(x_0,x_1,x_2,x_3,x_4)\neq {\bf 0}\;\}/{\bf C}^{\times},
\label{wp2o}
\end{eqnarray} 
where the ${\ C}^{\times}$ action is given by,
\begin{eqnarray}
(x_0,x_1,x_2,x_3,x_4)\rightarrow (\mu x_0,\mu x_1,\mu^{2}x_2,\mu^{2} x_3,\mu^{2}x_4).
\label{wp2tro}
\end{eqnarray}
It has one K\"ahler form and a singular ${\bf P}^{2}=\{[(0,0,x_2,x_3.x_4)]\}$.
In this section, we use another space $WP_1$ instead of ${\bf P}(1,1,2,2,2)$. which 
is obtained from blowing up ${\bf P}(1,1,2,2,2)$ along the singular ${\bf P}^2$. It is
a smooth complex manifold and was used in \cite{mp}. Explicitly, it is given as follows:
\begin{eqnarray}    
WP_1=\{({\bf a},{\bf b})\;|\; {\bf a}=(a_0,a_1)\in{\bf C}^2,\; {\bf b}=(b_0,b_1,b_2,b_3)\in{\bf C}^4 
\;,{\bf a}\neq{\bf 0},\;{\bf b}\neq{\bf 0}\;\;\}/({\bf C}^{\times})^2, 
\label{towp1}
\end{eqnarray}
where the two ${\bf C}^{\times}$ actions are given by, 
\begin{eqnarray}
(a_0,a_1,b_0,b_1,b_2,b_3)\rightarrow (\mu a_0,\mu a_1, b_0,b_1,b_2,\mu^{-2}b_3), \;
(a_0,a_1,b_0,b_1,b_2,b_3)\rightarrow (a_0,a_1, \nu b_0,\nu b_1,\nu b_3,\nu b_4).
\label{cwp1}
\end{eqnarray}
From the above definition, we can see that $WP_1$ is nothing but the projective 
bundle $\pi:{\bf P}({\cal O}_{{\bf P}^1}\oplus{\cal O}_{{\bf P}^1}\oplus
{\cal O}_{{\bf P}^1}\oplus{\cal O}_{{\bf P}^1}(-2))\rightarrow {\bf P}^1$.
Let ${\cal O}_{WP_1}(a)$ be $\pi^{*}{\cal O}_{{\bf P}^{1}}(1)$ and ${\cal O}_{WP_1}(b)$
be the dual line bundle of the universal bundle of ${\bf P}({\cal O}_{{\bf P}^1}\oplus{\cal O}_{{\bf P}^1}\oplus
{\cal O}_{{\bf P}^1}\oplus{\cal O}_{{\bf P}^1}(-2))$.
The classical cohomology ring of $WP_{1}$ is generated by two K\"ahler 
forms, $z=c_1({\cal O}_{WP_1}(a))$ and $w=c_1({\cal O}_{WP_1}(b))$. They obey the relations:
\begin{eqnarray}
z^2=0,\;\;w^{3}(w-2z)=0.
\label{relwp1}
\end{eqnarray}
As in the previous examples, integration of $\alpha\in H^{*}(WP_1,{\bf C})$ over $WP_1$ 
can be realized as residue integral in $z$ and $w$ as follows:
\begin{equation}
\int_{WP_1}\alpha=\frac{1}{(2\pi\sqrt{-1})^2}\oint_{C_{0}}\frac{dz}{z^2}\oint_{C_{(0,2z)}}
\frac{dw}{w^3(w-2z)}\alpha.
\label{creswp1}
\end{equation}
In the r.h.s. of (\ref{creswp1}), $\alpha$ should be considered as a polynomial in $z$ and $w$. 
Since $c_{1}(WP_1)=4w$, the Calabi-Yau hypersurface $XP_1\subset WP_1$ is given by  
the zero locus of a holomorphic section of ${\cal O}_{WP_1}(4b)$. 
Let $i$ be the inclusion map of $XP_1$.
In this subsection,
we consider the K\"ahler sub-ring $H^{*}_{K}(XP_1,{\bf C})$, which is a sub-ring of 
$H^{*}(XP_1,{\bf C})$ generated by $i^{*}z$ and $i^{*}w$. We denote $i^{*}z$ and $i^{*}w$
by $z$ and $w$ for brevity. In this subsection, we consider the following intersection number on 
$\widetilde{Mp}_{0,2}(WP_1,\dd)$
\begin{eqnarray}
w({\cal O}_{\alpha}{\cal O}_{\beta})_{0,\dd}:=
\int_{ [\widetilde{Mp}_{0,2}(WP_1,(d_a,d_b))]_{vir.} }ev_{1}^{*}(\alpha)\wedge ev_{2}^{*}(\beta)\wedge c_{top}({\cal E}_{\dd}).
\label{gwwp1} 
\end{eqnarray}
In (\ref{gwwp1}), $[\widetilde{Mp}_{0,2}(WP_1,\dd)]_{vir.}$ is defined in the same way as in the $F_3$ case and 
${\cal E}_{\dd}$ is an orbi-bundle that corresponds to 
$R^{0}\pi_{*}ev_{3}^{*}({\cal O}_{WP_1}(4b))$ on $\overline{M}_{0,2}(WP_1,\dd)$.
It is constructed in the same way as in the discussions in Subsection 2.3. The structure of the moduli space 
$\widetilde{Mp}_{0,2}(WP_1,(d_a,d_b))$ is almost the same as $\widetilde{Mp}_{0,2}(F_3,(d_a,d_b))$ and an obstruction bundle similar 
to the $F_3$ case also appears. The process of the localization computation is also the same as in the $F_3$ case except that we have $c_{top}({\cal E}_{\dd})$ in this case. But this can be easily done by looking back at the computation in \cite{vs}. Therefore, 
we write down only the data to compute $w({\cal O}_{\alpha}{\cal O}_{\beta})_{0,(d_a,d_b)}$ numerically.
We introduce here two rational functions in $z_{*}$ and $w_{*}$ in the same way as the $F_3$ case:
\begin{eqnarray}
G(\dd;z_{0},z_{1},w_0,w_1):=
\left\{\begin{array}{cc}
4w_0\frac{\prod_{j=1}^{2|\dd|-1}\bigl(\frac{-jz_0-(2|\dd|-j)z_1}{|\dd|}+w_0\bigr)}{\prod_{j=1}^{|\dd|-1}\bigl
(\frac{jz_{0}+(|\dd|-j)z_{1}}{|\dd|}\bigr)^{2}},&\;\mbox{if}\;\dd=(d,0),\\
&\\
\frac{\prod_{j=0}^{4|\dd|}\bigl(\frac{jw_0+(4|\dd|-j)w_1}{|\dd|}\bigr) }{\prod_{j=1}^{|\dd|-1}\bigl((\frac{jw_{0}+(|\dd|-j)w_{1}}{|\dd|})^3( \frac{jw_{0}+(|\dd|-j)w_{1}}{|\dd|}-2z_{0} )\bigr)},
&\;\mbox{if}\;\dd=(0,d),
\end{array}\right. 
\label{Gwp1}
\end{eqnarray}
\begin{eqnarray}
H(\dd_1;\dd_2,z_0,z_1,z_2,w_0,w_1,w_2):=\left\{\begin{array}{cc}
\frac{(-2z_1+w_1)}{4w_1 \bigl(\frac{z_1-z_0}{|\dd_1|}+ 
\frac{z_1-z_2}{|\dd_2|}\bigr)},&\;\mbox{if}\;\dd_{1}=(d_{a,1},0)\;\mbox{and}\;\dd_{1}=(d_{a,2},0),\\
\frac{1}{ 4w_1\bigl(\frac{z_1-z_0}{|\dd_1|}+ 
\frac{w_1-w_2}{|\dd_2|}\bigr)},&\;\mbox{if}\;\dd_{1}=(d_{a,1},0)\;\mbox{and}\;\dd_{2}=(0,d_{b,2}), \\
\frac{1}{ 4w_1\bigl(\frac{w_1-w_0}{|\dd_1|}+ 
\frac{z_1-z_2}{|\dd_2|}\bigr)},&\;\mbox{if}\;\dd_{1}=(0,d_{b,1})\;\mbox{and}\;\dd_{2}=(d_{a,2},0),\\
\frac{1}{ 4w_1\bigl(\frac{w_1-w_0}{|\dd_1|}+ 
\frac{w_1-w_2}{|\dd_2|}\bigr)},&\;\mbox{if}\;\dd_{1}=(0,d_{b,1})\;\mbox{and}\;\dd_{2}=(0,d_{b,2}).
\end{array}\right.
\label{Hwp1}
\end{eqnarray}
Then the integrand associated with $\sigma_{\dd}\in OP_{\dd}$ is given by, 
\begin{eqnarray}
K(\sigma_{\dd};z_{*},w_{*}):= \prod_{j=1}^{l(\sigma_{\dd})}G(\dd_j;z_{j-1},z_j,w_{j-1},w_j)
 \prod_{j=1}^{l(\sigma_{\dd})-1}H(\dd_j;\dd_{j+1},z_{j-1},z_j,z_{j+1},w_{j-1},w_j,w_{j+1}).
\label{Kwp1}
\end{eqnarray}
The integration rule of $K(\sigma_{\dd};z_{*},w_{*})$ is almost the same as the $F_{3}$ case,
\begin{eqnarray}
Res_{(WP_1)_0}(f)&:=&\left\{\begin{array}{cc}
 \frac{1}{2\pi\sqrt{-1}}\oint_{C_{0}}\frac{dz_0}{(z_0)^2}(f|_{w_0=w_1}),
&\;\mbox{if}\; \dd_{1}=(d_{a,1},0),\\
 \bigl(\frac{1}{2\pi\sqrt{-1}}\oint_{C_{(0,2z_0)}}\frac{dw_0}{(w_0)^3(w_0-2z_0)}f\bigr)|_{z_0=z_{1}},
&\;\mbox{if}\; \dd_{1}=(0,d_{b,1}),
\end{array}\right.\no\\ 
Res_{(WP_1)_j}(f)&:=&\left\{\begin{array}{cc}
 \frac{1}{2\pi\sqrt{-1}}\oint_{C_{0}}\frac{dz_j}{(z_j)^2}(f|_{w_j=w_{j+1}}),
&\;\mbox{if}\; \dd_{j+1}=(d_{a,j+1},0),\\
 \bigl(\frac{1}{2\pi\sqrt{-1}}\oint_{C_{(0,2z_j)}}\frac{dw_j}{(w_j)^3}f\bigr)|_{z_j=z_{j+1}},
&\;\mbox{if}\; \dd_{j}=(d_{a,j},0)\; \mbox{and}\; \dd_{j+1}=(0,d_{b,j+1}),\\
 \bigl(\frac{1}{2\pi\sqrt{-1}}\oint_{C_{(0,2z_j)}}\frac{dw_j}{(w_j)^3(w_j-2z_j)}f\bigr)|_{z_j=z_{j+1}},
&\;\mbox{if}\; \dd_{j}=(0,d_{b,j})\; \mbox{and}\; \dd_{j+1}=(0,d_{b,j+1}),
\end{array}\right.\no\\ 
&&\;\;\;\;\;\;\;\;\;\;\;\;\;\;\;\;\;\;\;\;\;\;\;\;\;\;\;\;\;\;\;\;\;\;\;\;\;\;\;\;\;\;\;\;\;\;\;\;\;\;\;\;\;
\;\;\;\;\;\;\;\;\;\;\;\;\;\;\;\;\;\;\;\;\;\;\;\;\;\;\;\;\;\;\;\;\;\;\;\;\;\;\;\;
(j=1,\cdots,l(\sigma_{\dd})-1),\no\\
Res_{(WP_1)_{l(\sigma_{\dd})}}(f)&:=&\left\{\begin{array}{cc}
\frac{1}{(2\pi\sqrt{-1})^2}\oint_{C_{0}}\frac{dz_{l(\sigma_{\dd})}}{(z_{l(\sigma_{\dd})})^2}\oint_{C_{(0,2z_{l(\sigma_{\dd})})}}
\frac{dw_{l(\sigma_{\dd})}}{(w_{l(\sigma_{\dd})})^3}f ,
&\;\mbox{if}\; \dd_{l(\sigma_{dd})}=(d_{a,l(\sigma_{\dd})},0),\\
 \frac{1}{(2\pi\sqrt{-1})^2}\oint_{C_{0}}\frac{dz_{l(\sigma_{\dd})}}{(z_{l(\sigma_{\dd})})^2}\oint_{C_{(0,2z_{l(\sigma_{\dd})})}}\frac{dw_{l(\sigma_{\dd})}}
{(w_{l(\sigma_{\dd})})^3(w_{l(\sigma_{\dd})}-2z_{l(\sigma_{\dd})})}f,
&\;\mbox{if}\;  \dd_{l(\sigma_{\dd})}=(0,d_{b,l(\sigma_{\dd})}),
\end{array}\right.\no\\
\label{reswp1}
\end{eqnarray}
except that we also take the residue at $w_{j}=2z_{j}$ in the fourth and sixth lines of (\ref{reswp1}). It seems a little bit 
unnatural from geometrical point of view, but we need to do it to obtain the correct numerical results.
The reason of this modification seems to be a problem which should be pursued further.
 With this setup 
, contributions from $\sigma_{\dd}\in OP_{\dd}$ to $w({\cal O}_{\alpha}{\cal O}_{\beta})_{0,\dd}$ are given by,
\begin{eqnarray}
&& Amp(\sigma_{\dd};\alpha,\beta) =\biggl(\prod_{j=1}^{l(\sigma_{\dd})}\frac{1}{|\dd_{j}|}\biggr) 
Res_{(WP_1)_{l(\sigma_{\dd})}}( Res_{(WP_1)_{l(\sigma_{\dd})-1}}(
\cdots Res_{(WP_1)_{0}}( ev_{1}^{*}(\alpha)ev_{2}^{*}(\beta)K(\sigma_{\dd};z_{*},w_{*}) )\cdots)),\no\\
\label{amp3-2}
\end{eqnarray}
where $ev_1(z^sw^t)=(z_0)^s(w_0)^t$ (resp. $ev_2(z^sw^t)=(z_{l(\sigma_{\dd})})^s(w_{l(\sigma_{\dd})})^t$).
Finally,we obtain $w({\cal O}_{\alpha}{\cal O}_{\beta})_{0,\dd}$ as usual:
\begin{eqnarray}
 w({\cal O}_{\alpha}{\cal O}_{\beta})_{0,\dd} =\sum_{\sigma_{\dd}\in OP_{\dd}}Amp(\sigma_{\dd};\alpha,\beta). 
\end{eqnarray}
\subsubsection{Numerical Results and the Mirror Computation}
We present below numerical results of $w({\cal O}_{\alpha}{\cal O}_{\beta})_{0,\dd}$ up to $d_a+d_b\leq 3$ by using the
generating function:
\begin{eqnarray}
w({\cal O}_{\alpha}{\cal O}_{\beta})_{0}:=\left(\int_{XP_1}\alpha\beta z\right)x_1+
\left(\int_{XP_1}\alpha\beta w\right)x_2+\sum_{\dd>(0,0)} w({\cal O}_{\alpha}{\cal O}_{\beta})_{0,\dd}
e^{{d_a}x_1+{d_b}x_2}.
\label{genxp1}
\end{eqnarray}
In (\ref{genxp1}), the classical intersection number $\int_{XP_1}\alpha\beta\gamma$ is given by, 
\begin{eqnarray}
\int_{XP_1}\alpha\beta\gamma=\frac{1}{(2\pi\sqrt{-1})^2}\oint_{C_0}\frac{dz}{z^2}\oint_{C_{(0,2z)}}\frac{dw}{w^3(w-2z)}
4w\cdot\alpha\beta\gamma,
\label{clxp1}
\end{eqnarray}
where $\alpha$, $\beta$ and $\gamma$ in the r.h.s. are regarded as polynomials in $z$ and $w$.
\begin{eqnarray}
w({\cal O}_1{\cal O}_{w^2})_0&=&4x_{1}+8x_{2}+1024e^{x_{2}}+103872e^{2x_{2}}+\frac{46099456}{3}
e^{3x_{2}}+216576e^{x_{1}+2x_2}+\cdots,\no\\
w({\cal O}_1{\cal O}_{zw})_0&=&4x_{2}+416e^{x_{2}}-4e^{x_{1}}+39120
e^{2x_{2}}-6e^{2x_{1}}+192e^{x_{1}+x_{2}}+\frac{16567040}{3}e^{3x_{2}}-\frac{40}{3}e^{3x_{1}}+\no\\
&&+133920e^{x_{1}+2x_2}+192e^{2x_{1}+x_{2}}+\cdots,\no\\
w({\cal O}_z{\cal O}_z)_0&=&4e^{x_{1}}+10e^{2x_{1}}+832e^{x_{1}+x_{2}}+
\frac{88}{3}e^{3x_{1}}+199744e^{x_{1}+2x_2}+832e^{2x_{1}+x_{2}}+\cdots,\no\\
w({\cal O}_w{\cal O}_w)_0&=&4x_{1}+8x_2
+1664e^{x_{2}}+210880e^{2x_{2}}+\frac{108286976}{3}e^{3x_{2}}+486016e^{x_{1}+2x_2}+\cdots,\no\\
w({\cal O}_z{\cal O}_w)_0&=&4x_{2}+416e^{x_{2}}-4e^{x_{1}}+39120e^{2x_{2}}-6e^{2x_{1}}+832e^{x_{1}+x_{2}}+
\frac{16567040}{3}e^{3x_{2}}-\frac{40}{3}e^{3x_{1}}+\no\\
&&+375648e^{x_{1}+2x_2}+832e^{2x_{1}+x_{2}}+\cdots.
\label{dat1wp1}
\end{eqnarray}
Let $\eta_{\alpha\beta}$ be $\int_{XP_1}\alpha\beta$, i.e., the $(\alpha,\beta)$-element of classical intersection 
matrix of $XP_1$ and $\eta^{\alpha\beta}$ be the $(\alpha,\beta)$-element of the inverse of $(\eta_{\alpha\beta})$. 
One of our conjectures in this example is that $\eta^{z\alpha}w({\cal O}_1{\cal O}_{\alpha})_0$ and 
$\eta^{w\alpha}w({\cal O}_1{\cal O}_{\alpha})_0$ coincide with the mirror maps used in the standard mirror 
computation \cite{hkty}. Indeed, our numerical results: 
\begin{eqnarray}
t_1 &=&\frac{1}{4}w({\cal O}_1{\cal O}_{w^2})_0-\frac{1}{2}w({\cal O}_1{\cal O}_{zw})_0=\no\\
&=&{x_{1}} + 48e^{x_{2}} + 6408e^{2x_{2}} + 
1080448e^{3x_{2}} - 12816e^{x_1+2x_{2}} + 2e^{x_{1}} + 3e^{2x_{1}}- 96e^{x_{1}+x_{2}} 
+\no\\ 
&&+\frac {20}{3}e^{3x_{1}} - 96e^{2x_{1}+x_{2}}+\cdots,\no\\
t_2 &=&\frac{1}{4}w({\cal O}_1{\cal O}_{zw})_0=\no\\ 
&=&{x_{2}} + 104e^{x_{2}} - e^{x_{1}}
 + 9780e^{2x_{2}} - \frac {3}{2}e^{2x_{1}} 
+ 48e^{x_{1}+x_{2}} +\frac {4141760}{3}e^{3x_{2}} - \frac {10}{3}e^{3x_{1}}+ 
\no\\
&&+33480e^{x_{1}+2x_2} + 48e^{2x_{1}+x_{2}} +\cdots,
\label{dat2wp1}
\end{eqnarray}
give us the standard mirror maps in \cite{hkty}.
We then invert (\ref{dat2wp1}) and substitute $x_{i}=x_{i}(t_1,t_2)$ into $w({\cal O}_z{\cal O}_z)_0$, 
$w({\cal O}_z{\cal O}_w)_0$ and $w({\cal O}_w{\cal O}_w)_0$.
The results:
\begin{eqnarray}
w({\cal O}_z{\cal O}_z)_0|_{x_i=x_i(t_1,t_2)}&=&4t_{1} + 8t_{2} + 640e^{t_{2}}
 + 40448e^{2t_{2}} + 640e^{t_{1}+t_{2}} + 
\frac {7787008}{3} e^{3t_{2}} + 288896e^{t_{1}+2t_{2}}+\cdots, \no\\
w({\cal O}_z{\cal O}_w)_0|_{x_i=x_i(t_1,t_2)}&=&4e^{t_{1}} + 640e^{t_{1}+t_{2}} + 2e^{2t_{1}} + 72224
e^{t_{1}+2t_{2}} + \frac {4}{3} e^{3t_{1}}+\cdots \no\\
w({\cal O}_w{\cal O}_w)_0|_{x_i=x_i(t_1,t_2)}&=& 4{t_{2}} + 640e^{t_{1}+t_{2}} + 
144448e^{t_{1}+2t_{2}}+\cdots,\no\\
\label{dat3wp1}
\end{eqnarray}
give us the generating functions of 2-point Gromov-Witten invariants of $XP_1$. These results are evidences 
of Conjecture 1 in the case of Calabi-Yau hypersurface of $WP_1$.

\newpage
\section{Generalizations to Weighted Projective Space with One K\"ahler Form}
\subsection{K3 surface in ${\bf P}(1,1,1,3)$}
This subsection deals with results on the $j$-invariant of elliptic curves arising from our 
conjecture on mirror map. The $j$-invariant is a modular function of $\tau$: the flat coordinate 
of the moduli space of complex structures of elliptic curves, and its Fourier expansion is given by, 
\begin{eqnarray}
j(q)&=&q^{-1}+744+196884q+21493760q^2+864299970q^3+20245856256q^4+333202640600q^5+\cdots,\no\\
     &=:&q^{-1}+\sum_{d=1}j_dq^{d-1},\no\\ 
&&(q=\exp( 2\pi\sqrt{-1}\tau )).
\label{j}
\end{eqnarray}
By inverting (\ref{j}), we can express $2\pi\sqrt{-1}\tau$ as a power series in $j^{-1}$:
\begin{eqnarray}
2\pi\sqrt{-1}\tau &=&-\log(j)+744j^{-1}+473652j^{-2}+451734080j^{-3}+510531007770j^{-4}+\frac{3169342733223744}{5}j^{-5}+\cdots,\no\\
                  &=:&-\log(j)+\sum_{d=1}^{\infty}w_{d}j^{-d}.
\end{eqnarray}
Let $WP_{2}$ be the weighted projective space ${\bf P}(1,1,1,3)$:
\begin{eqnarray}
WP_{2}:=\{{\bf a}=(a_{0},a_{1},a_{2},a_{3})\;|\;{\bf a}\neq {\bf 0}\;\}/{\bf C}^{\times},
\label{d1wp1}
\end{eqnarray}
where the ${\bf C}^{\times}$ action is given by,
\begin{equation}
(a_0,a_1,a_2,a_3)\rightarrow (\mu a_0,\mu a_1,\mu a_2,\mu^3 a_3).
\label{twp2}
\end{equation}
We denote by ${\cal O}_{WP_2}(1)$ the line bundle whose holomorphic section is generated by $a_{0}$, $a_{1}$ and $a_{2}$. 
Let $z$ be $c_{1}({\cal O}_{WP_2}(1))$. Then $H^{*}(WP_{2},{\bf C})$ is isomorphic to ${\bf C}[z]/(z^4)$ and integration of 
$\alpha\in H^{*}(WP_{2},{\bf C})$ can be realized as the following residue integral:
\begin{eqnarray}
\int_{WP_2}\alpha=\frac{1}{3}\cdot\frac{1}{2\pi\sqrt{-1}}\oint_{C_0}\frac{dz}{z^4}\alpha,
\label{intwp2}
\end{eqnarray}
where $\alpha$ on the r.h.s. is regarded as a polynomial in $z$. The factor $\frac{1}{3}$ comes from the fact that $WP_2$ is 
an orbifold with ${\bf Z}_{3}$ singularity $[(0,0,0,1)]$. It is well-known that the zero locus of a holomorphic section of $ 
{\cal O}_{WP_2}(6)$ is a K3 surface. Let $XP_2$ be this K3 surface. In \cite{by}, it was proved that the mirror map used in 
the mirror computation of $XP_2$ is given by:
\begin{eqnarray} 
t&=&x+\sum_{d=1}^{\infty}w_d e^{dx}\no\\
 &=&x+744e^x+473652e^{2x}+451734080e^{3x}+510531007770e^{4x}+\frac{3169342733223744}{5}e^{5x}+\cdots,
\end{eqnarray}
where $t$ is the flat coordinate of K\"ahler moduli space of $XP_2$ and $x$ is a standard complex deformation parameter of the mirror 
manifold of $XP_2$. At this stage, we consider the following intersection number on $\widetilde{Mp}_{0,2}(WP_2,d)$:   
\begin{eqnarray}
w({\cal O}_{\alpha}{\cal O}_{\beta})_{0,d}:=
\int_{\widetilde{Mp}_{0,2}(WP_2,d)}ev_{1}^{*}(\alpha)\wedge ev_{2}^{*}(\beta)\wedge c_{top}({\cal E}_{d}),
\label{gwwp2} 
\end{eqnarray}
where  ${\cal E}_{d}$ is a sheaf that corresponds to $R^{0}\pi_{*}ev_{3}^{*}({\cal O}_{WP_2}(6))$ on $\overline{M}_{0,2}(WP_2,d)$.
We briefly mention the structure of $\widetilde{Mp}_{0,2}(WP_2,d)$. For brevity, we write $(a_0,a_1,a_2,a_3)$ as $({\bf a},a_3)$.
Then a polynomial map from $CP^1$ to $WP_2$ of degree $d$ is written as follows:
\begin{eqnarray} 
(\sum_{j=0}^{d}{\bf a}_js^{j}t^{d-j},\sum_{j=0}^{3d}a_{3,j}s^{j}t^{3d-j}).
\label{polywp2}
\end{eqnarray}
Therefore, $Mp_{0,2}(WP_2,d)$ is constructed as follows:
\begin{eqnarray}
Mp_{0,2}(WP_2,d)=\{({\bf a}_0,{\bf a}_1,\cdots,{\bf a}_d,a_{3,0},a_{3,1},\cdots,a_{3,3d})\;|\;({\bf a}_0,a_{3,0}),({\bf a}_d,a_{3,3d})\neq
{\bf 0}\}/({\bf C}^{\times})^2,
\label{mwp2}
\end{eqnarray}
where the two ${\bf C}^{\times}$ actions are given by,
\begin{eqnarray}
({\bf a}_0,{\bf a}_1,\cdots,{\bf a}_d,a_{3,0},a_{3,1},\cdots,a_{3,3d})&\rightarrow&
(\mu{\bf a}_0,\mu{\bf a}_1,\cdots,\mu{\bf a}_d,\mu^3a_{3,0},\mu^3a_{3,1},\cdots,\mu^3a_{3,3d}),\no\\
({\bf a}_0,{\bf a}_1,\cdots,{\bf a}_d,a_{3,0},a_{3,1},\cdots,a_{3,3d})&\rightarrow&
({\bf a}_0,\nu{\bf a}_1,\nu^2{\bf a}_2,\cdots,\nu^d{\bf a}_d,a_{3,0},\nu a_{3,1},\nu^2 a_{3,1},\cdots,\nu^{3d}a_{3,3d}).
\end{eqnarray}
Additional divisors added to construct $\widetilde{Mp}_{0,2}(WP_2,d)$ are fundamentally the same as in the $CP^{N-1}$ case. Therefore, 
a point in $\widetilde{Mp}_{0,2}(WP_2,d)$ can be represented as,
\begin{eqnarray}
[(({\bf a}_0,{\bf a}_1,\cdots,{\bf a}_d,a_{3,0},a_{3,1},\cdots,a_{3,3d},u_1,u_2,\cdots,u_{d-1})],
\label{wmwp2}
\end{eqnarray}
where $[*]$ means  taking the equivalence class under the $({\bf C}^{\times})^{d+1}$ action. 
We then compute $w({\cal O}_{\alpha}{\cal O}_{\beta})_{0,d}$ by using localization under the torus action flow:
\begin{eqnarray}
&&[((e^{\lambda_0 t}{\bf a}_0,e^{\lambda_1 t}{\bf a}_1,\cdots,e^{\lambda_d t}{\bf a}_d,
e^{\lambda_{3,0}t}a_{3,0},e^{\lambda_{3,1}t}a_{3,1},\cdots,e^{\lambda_{3,3d}t}a_{3,3d},u_1,u_2,\cdots,u_{d-1})],\no\\
&&\hspace{7.5cm}(\lambda_{3,3j}=3\lambda_{j},\;j=0,1,2,\cdots,d).
\label{toruswp2}
\end{eqnarray}
As in the $CP^{N-1}$ case, the connected components of the fixed point set are labeled by ordered partitions of the positive integer $d$: 
\begin{equation}
OP_{d}=\{\sigma_{d}=(d_{1},d_{2},\cdots,d_{l(\sigma_{d})})\;\;|\;\;
\sum_{j=1}^{l(\sigma_{d})}d_{j}=d\;\;,\;\;d_{j}\in{\bf N}\}.
\label{part} 
\end{equation}
Let $F_{\sigma_d}$ be the connected component labeled by $\sigma_d$. As in the previous examples, it is given by an orbifold: 
\begin{eqnarray}
(WP_2)_{0}\times(WP_{2})_1\times(WP_{2})_2\times\cdots\times(WP_2)_{l(\sigma_{d})},
\label{orbiwp2}
\end{eqnarray}
on which $\mathop{\oplus}_{j=1}^{l(\sigma_{d})}\left({\bf Z}/(d_{j}{\bf Z})\right)$ acts. Now, we define two rational functions 
in $z_{*}$ to write down the integrand for $F_{\sigma_d}$:
\begin{eqnarray}
G(d;z_{0},z_{1}):=
\frac{\prod_{j=0}^{6d}\bigl(\frac{jz_0+(6d-j)z_1}{d}\bigr) }{\prod_{j=1}^{d-1}(\frac{jz_{0}+(d-j)z_{1}}{d})^3\prod_{j=1}^{3d-1}(\frac{jz_{0}+(3d-j)z_{1}}{3d})}.
\label{Gwp2}
\end{eqnarray}
\begin{eqnarray}
H(d_1,d_2;z_0,z_1,z_2):=
\frac{1}{6z_1 \bigl(\frac{z_1-z_0}{d_1}+ 
\frac{z_1-z_2}{d_2}\bigr)}
\label{Hwp2}
\end{eqnarray}
As in the previous cases, the integrand is given by, 
\begin{eqnarray}
K(\sigma_{d};z_{*}):= \prod_{j=1}^{l(\sigma_{d})}G(d_j;z_{j-1},z_j)
 \prod_{j=1}^{l(\sigma_{d})-1}H(d_j,d_{j+1};z_{j-1},z_j,z_{j+1}).
\label{Kwp2}
\end{eqnarray}
Looking back at (\ref{intwp2}), we introduce the following operation:
\begin{eqnarray}
Res_{(WP_2)_j}(f)&:=&
 \frac{1}{3}\cdot\frac{1}{2\pi\sqrt{-1}}\oint_{C_{0}}\frac{dz_j}{(z_j)^4}f.
\label{reswp2}
\end{eqnarray}
Then contribution from $F_{\sigma_d}$ to $w({\cal O}_{\alpha}{\cal O}_{\beta})_{0,d}$ is given by,  
\begin{eqnarray}
&& Amp(\sigma_{d};\alpha,\beta) =\biggl(\prod_{j=1}^{l(\sigma_{d})}\frac{1}{d_{j}}\biggr) 
Res_{(WP_2)_{l(\sigma_{d})}}( Res_{(WP_2)_{l(\sigma_{d})-1}}(
\cdots Res_{(WP_2)_{0}}( ev_{1}^{*}(\alpha)ev_{2}^{*}(\beta)K(\sigma_{d};z_{*})\cdots)),\no\\
\label{ampwp2}
\end{eqnarray}
where $ev_{1}^{*}(z^s)=(z_{0})^s$ (resp. $ev_{2}^{*}(z^s)=(z_{l(\sigma_{d})})^s$).
Finally, we add up the contributions from all the $F_{\sigma_{\dd}}$'s and obtain the formula:
\begin{eqnarray}
 w({\cal O}_{\alpha}{\cal O}_{\beta})_{0,d} =\sum_{\sigma_{d}\in OP_{d}}Amp(\sigma_{d};\alpha,\beta). 
\end{eqnarray}
Now, our conjecture in this example becomes,
\begin{conj}
\begin{equation}
w_{d}=\frac{1}{2}w({\cal O}_{1}{\cal O}_{z})_{0,d}.
\end{equation}
\end{conj}
We checked the above equality up to degree 5. As a by-product of this conjecture, we can represent the Fourier coefficient $j_d$ 
of the $j$-invariant in terms of the intersection number $w({\cal O}_{\alpha}{\cal O}_{\beta})_{0,d}$ as follows:
\begin{cor}
\begin{equation}
j_{d}=\sum_{\sigma_d\in OP_d}(-(d-1))^{l(\sigma_d)-1}\frac{1}{(l(\sigma_d))!}
\prod_{j=1}^{l(\sigma_{d})}\bigl(\frac{w({\cal O}_{1}{\cal O}_{z})_{0,d_j}}{2}\bigr).
\end{equation}
\end{cor}
The above equation easily follows from standard combinatorics of inversion of power series. 
\subsection{Calabi-Yau Hypersurface in ${\bf P}(1,1,2,2,6)$}
As our last example, we deal with the Calabi-Yau hypersurface in ${\bf P}(1,1,2,2,6)$, which was discussed in much previous work 
\cite{hkty} \cite{KKLMV}, \cite{kv}. As in the case of ${\bf P}(1,1,2,2,2)$, we use the following toric manifold: 
\begin{eqnarray}    
WP_3=\{({\bf a},{\bf b})\;|\; {\bf a}=(a_0,a_1)\in{\bf C}^2,\; {\bf b}=(b_0,b_1,b_2,b_3)\in{\bf C}^4 
\;,{\bf a}\neq{\bf 0},\;{\bf b}\neq{\bf 0}\;\;\}/({\bf C}^{\times})^2, 
\label{towp3}
\end{eqnarray}
where the two ${\bf C}^{\times}$ actions are given by, 
\begin{eqnarray}
(a_0,a_1,b_0,b_1,b_2,b_3)\rightarrow (\mu a_0,\mu a_1, b_0,b_1,b_2,\mu^{-2}b_3), \;
(a_0,a_1,b_0,b_1,b_2,b_3)\rightarrow (a_0,a_1, \nu b_0,\nu b_1,\nu^3 b_3,\nu b_4).
\label{cwp3}
\end{eqnarray}
It can be obtained by blowing up ${\bf P}(1,1,2,2,6)$ along the singular ${\bf P}(1,1,3)$ in ${\bf P}(1,1,2,2,6)$. Let ${\cal O}_{WP_3}(a)$ be 
a line bundle whose holomorphic section is generated by $a_0$ and $a_1$ and let ${\cal O}_{WP_3}(b)$ be 
a line bundle whose holomorphic section is generated by $b_0$ and $b_1$. We denote $c_1({\cal O}_{WP_3}(a))$ 
(resp. $c_1({\cal O}_{WP_3}(b))$) by $z$ (resp. $w$). Then we can consider the following intersection number on 
$\widetilde{Mp}_{0,2}(WP_3,\dd)$:
\begin{eqnarray}
w({\cal O}_{\alpha}{\cal O}_{\beta})_{0,(d_a,d_b)}:=
\int_{ [\widetilde{Mp}_{0,2}(WP_3,(d_a,d_b))]_{vir.} }ev_{1}^{*}(\alpha)\wedge ev_{2}^{*}(\beta)\wedge c_{top}({\cal E}_{\dd}).
\label{gwwp3} 
\end{eqnarray}
where ${\cal E}_{\dd}$ is an orbi-bundle on $\widetilde{Mp}_{0,2}(WP_3,\dd)$ that corresponds to 
$R^{0}\pi_{*}ev_{3}^{*}({\cal O}_{WP_3}(6b))$ on $\overline{M}_{0,2}(WP_3,\dd)$.
From (\ref{towp3}) and (\ref{cwp3}), we can see that $WP_3$ is a ${\bf P}(1,1,1,3)$ bundle over ${\bf P}^{1}$. Therefore, it is 
straightforward to compute $w({\cal O}_{\alpha}{\cal O}_{\beta})_{0,(d_a,d_b)}$ by combining the result of $WP_1$ with the 
one of $WP_2$. We leave the remaining computations to readers  as an exercise. We end by presenting numerical results 
of $w({\cal O}_{\alpha}{\cal O}_{\beta})_{0,\dd}$ in the form of generating function:
\begin{eqnarray}
w({\cal O}_1{\cal O}_{w^{2}}) &:=& 2{x_{1}} + 4{x_{2}} 
+ 3456e^{x_{2}} + 2335968e^{2x_{2}} + 2313054720e^{3x_{2}} + 4836096e^{x_{1}+2x_2}+\cdots,\no\\
w({\cal O}_1{\cal O}_{zw})  &:=& 2{x_{2}} + 1488e^{x_{2}} - 
2e^{x_{1}} + 947304e^{2x_{2}} - 3e^{2x_{1}} + 
480e^{x_{1}+x_{2}} + 903468160e^{3x_{2}}- \no\\
&&- \frac {20}{3} e^{3x_{1}} + 
2859408e^{x_{1}+2x_2} + 480e^{2x_{1}+x_{2}}+\cdots, \no\\
w({\cal O}_z{\cal O}_{z})  &:=& 2e^{x_{1}} + 5e^{2x_{1}} + 
2976e^{x_{1}+x_{2}} +  \frac {44}{3} e^{3x_{1}} + 
4896288e^{x_{1}+2x_2} + 2976e^{2x_{1}+x_{2}}+\cdots,\no\\
w({\cal O}_w{\cal O}_{w}) &:=& 2{x_{1}} + 4{x_{2}}
 + 5952e^{x_{2}} + 5089248e^{2x_{2}} + 5867470336e^{3x_{2}}
 + 12006720e^{x_{1}+2x_2}+\cdots,\no\\
w({\cal O}_z{\cal O}_{w})  &:=& 2{x_{2}} + 1488e^{x_{2}} - 
2e^{x_{1}} + 947304e^{2x_{2}} - 3e^{2x_{1}} + 
2976e^{x_{1}+x_{2}} + 903468160e^{3x_{2}}- \no\\
&& - \frac {20}{3} e^{3x_{1}} 
+ 9198000e^{x_{1}+2x_2} + 2976e^{2x_{1}+x_{2}}+\cdots. 
\label{datwp3}
\end{eqnarray} 
Of course, we can perform the mirror computation by using these results as in the ${\bf P}(1,1,2,2,2)$ case.


\newpage  
\begin{thebibliography}{99}
\bibitem{givc} T.Coates, A.B.Givental.
\newblock{\em Quantum Riemann-Roch, Lefschetz and Serre}
\newblock{ Ann. of Math. (2) 165 (2007), no. 1, 15--53.}
\bibitem{ckyz} T.-M. Chiang, A. Klemm, S.-T. Yau, E. Zaslow.
\newblock{\em Local Mirror Symmetry: Calculations and Interpretations}
\newblock{Adv.Theor.Math.Phys. 3 (1999) 495-565.}
\bibitem{fj1}B.Forbes, M.Jinzenji 
\newblock{\em $J$ functions, non-nef toric varieties 
and equivariant local mirror symmetry of curves}
\newblock{ Int. J. Mod. Phys. A 22 (2007), no. 13, 2327--2360.}
\bibitem{fj3} B.Forbes and M.Jinzenji,
\newblock{ \em Extending the Picard-Fuchs 
system of local mirror symmetry.}
\newblock{ J.Math.Phys. 46 (2005) 082302.} 
\bibitem{fj4} B.Forbes and M.Jinzenji,
\newblock{\em Prepotentials for local mirror symmetry via Calabi-Yau fourfolds}
\newblock{J. High Energy Phys. 2006, no. 3, 061, 42 pp.}
\bibitem{guest} M.A.Guest.
\newblock{\em From Quantum Cohomology to Integrable Systems (Oxford Graduate Texts in Mathematics) }
\newblock{Oxford University Press, (2008).  }
\bibitem{hkty} S.Hosono, A.Klemm, S.Theisen, S.-T.Yau.
\newblock{\em Mirror symmetry, mirror map and applications to Calabi-Yau hypersurfaces}
\newblock{ Comm. Math. Phys. 167 (1995), no. 2, 301--350.}
\bibitem{iri} H.Iritani. 
\newblock{\em Quantum D-modules and Generalized Mirror Transformations}
\newblock{ Topology 47 (2008), no. 4, 225--276.}
\bibitem{fano} M.Jinzenji. 
\newblock{\em Completion of the Conjecture: 
Quantum Cohomology of Fano Hypersurfaces}
\newblock{Mod.Phys.Lett. A15 (2000) 101-120.}
\bibitem{vs}M.Jinzenji.
\newblock{\em Virtual Structure Constants as Intersection Numbers of Moduli Space of Polynomial Maps with Two Marked Points }
\newblock{Letters in Mathematical Physics, Vol.86, No.2-3, 99-114 (2008) }
\bibitem{gene0}M.Jinzenji.
\newblock{\em On the Quantum Cohomology Rings of General Type Projective Hypersurfaces and Generalized Mirror Transformation}
\newblock{ Int.J.Mod.Phys. A15 (2000) 1557-1596.}
\bibitem{gene}M.Jinzenji. 
\newblock{\em Coordinate change of Gauss-Manin system and generalized mirror transformation}
\newblock{ Internat. J. Modern Phys. A 20 (2005), no. 10, 2131--2156.}
\bibitem{gm}M.Jinzenji.
\newblock{\em Gauss-Manin System and the Virtual Structure Constants}
\newblock{Int.J.Math. 13 (2002) 445-478.}
\bibitem{prmir}M. Jinzenji.
\newblock{\em Direct Proof of Mirror Theorem of Projective Hypersurfaces up to degree 3 Rational Curves}
\newblock{Preprint, arXiv:0902.3863.}
\bibitem{KKLMV}S. Kachru, A. Klemm, W. Lerche, P. Mayr, C. Vafa.
\newblock{\em  Nonperturbative Results on the Point Particle Limit of N=2 Heterotic String Compactifications }
\newblock{ Nucl.Phys.B459:537-558,1996.}
\bibitem{kv}S. Kachru, C. Vafa.
\newblock{\em Exact Results for N=2 Compactifications of Heterotic Strings }
\newblock{Nucl.Phys.B450:69-89,1995.}
\bibitem{by}B.Lian, S.-T. Yau.
\newblock{\em Arithmetic properties of mirror map and quantum coupling}
\newblock{ Comm. Math. Phys. 176 (1996), no. 1, 163--191.}
\bibitem{mp}  D.R.Morrison, M.R.Plesser.
\newblock{\em Summing the Instantons: Quantum Cohomology and Mirror Symmetry in Toric Varieties}
\newblock{Nucl.Phys. B440 (1995) 279-354}
\end{thebibliography}
\end{document}